\documentclass[submission, Phys]{SciPost}
\usepackage{makecell}
\usepackage{amsmath,amsfonts,amssymb,array}
\usepackage{latexsym, array, bm}
\usepackage{lineno}
\usepackage{pgfplots}
\pgfplotsset{compat=1.16}
\usepackage{enumerate}
\usepackage{float}
\usepackage{graphicx}
\usepackage{multirow}
\usepackage{color}
\usepackage{ltablex}
\usepackage{caption, booktabs}
\usepackage{geometry}
\usepackage{caption}
\usepackage{subcaption}
\usepackage{bm}
\usepackage{mathtools}
\usepackage{diagbox}
\usepackage{adjustbox}
\usepackage{mathtools}
\usepackage{longtable}
\usepackage{algorithmic}
\usepackage[algo2e,linesnumbered,ruled,vlined,boxed,noresetcount,algosection]{algorithm2e}

\usepackage{url}
\usepackage{fancyvrb}
\newcommand{\email}[1]{\url{#1}}

\makeatletter
\g@addto@macro\bfseries{\boldmath}
\makeatother
\newcommand{\calO}{\mathcal{O}}
\newcommand{\abscalO}{\ensuremath{|\mathcal{O}|}}
\newcommand{\calB}{\mathcal{B}}

\newcommand{\calR}{\mathcal{R}}
\newcommand{\calS}{\mathcal{S}}
\newcommand{\tmax}{\text{max}}
\newcommand{\tmin}{\text{min}}
\newcommand{\tnew}{\text{new}}

\newcommand{\tcand}{\text{cand}}

\newcommand{\equaref}[1]{Eq.~(\ref{#1})}

\newcommand{\figref}[1]{Fig.~\ref{#1}}

\newcommand{\algref}[1]{Algorithm~\ref{#1}}

\usepackage[normalem]{ulem}

\newcommand{\mini}{\mathop{\mbox{minimize}}}
\newcommand{\maxi}{\mathop{\mbox{maximize}}}

\newcommand{\st}{\mbox{subject to }}

\newcommand{\be}{\begin{equation}}
\newcommand{\ee}{\end{equation}}
\newcommand{\bea}{\begin{eqnarray}}
\newcommand{\eea}{\end{eqnarray}}

\newcommand{\bvec}{\left(\begin{array}{c}}
	\newcommand{\evec}{\end{array}\right)}
\newcommand{\bsub}{\begin{subequations}}
	\newcommand{\esub}{\end{subequations}}

\usepackage{lipsum}
\usepackage{xspace,amsmath,booktabs,graphicx}
\usepackage[T1]{fontenc} %

\newcommand{\MC}{\ensuremath{\text{MC}}\xspace}

\newcommand{\chisq}{\ensuremath{\chi^2}\xspace}

\renewcommand{\vec}[1]{\ensuremath{\mathbf{#1}}\xspace}
\newcommand{\p}{\vec{p}}
\newcommand{\w}{\vec{w}}
\newcommand{\x}{\vec{x}}

\newcommand{\phatw}{\ensuremath{\widehat{\p}_\w\xspace}}
\newcommand{\phatws}{\ensuremath{\widehat{\p}_{\w^{*}}\xspace}}

\newcommand{\professor}{\textsc{Professor}\xspace}
\newcommand{\apprentice}{\textsc{apprentice}\xspace}

\newcommand{\pythia}{\textsc{Pythia}\xspace}

\newcommand{\sherpa}{\textsc{Sherpa}\xspace}

\newcommand{\MeV}{\text{Me\kern -0.15ex V}\xspace}
\newcommand{\GeV}{\text{Ge\kern -0.15ex V}\xspace}
\newcommand{\TeV}{\text{Te\kern -0.15ex V}\xspace}

\usepackage{todonotes}

\begin{document}

\begin{center}{\Large \textbf{BROOD:   Bilevel  and Robust Optimization and
      Outlier Detection for Efficient Tuning of High-Energy Physics Event Generators}}\end{center}

\begin{center}
Wenjing Wang\textsuperscript{1},
Mohan Krishnamoorthy\textsuperscript{2},
Juliane M\"uller\textsuperscript{1*},
Stephen Mrenna\textsuperscript{3*},
Holger Schulz\textsuperscript{4},
Xiangyang Ju\textsuperscript{1},
Sven Leyffer\textsuperscript{2},
Zachary Marshall\textsuperscript{1}
\end{center}

\begin{center}
{\bf 1} Lawrence Berkeley National Laboratory, Berkeley, CA 94720
\\
{\bf 2} Argonne National Laboratory, Lemont, IL 60439
\\
{\bf 3} Fermi National Accelerator Laboratory, Batavia, IL 60510
\\
{\bf 4} Department of Computer Science, Durham University, South Road, Durham DH1 3LE, UK
\\
* julianemueller@lbl.gov mrenna@fnal.gov 
\end{center}

\begin{center}
	\today
\end{center}


\section*{Abstract}
The parameters in Monte Carlo (MC) event generators are tuned on experimental measurements by evaluating the goodness of
fit between the data and the MC predictions.    The relative importance of each measurement is adjusted manually in an often time-consuming, iterative process to meet different experimental needs.  In this work, we introduce several optimization formulations and algorithms with new decision criteria for streamlining and automating this process. These algorithms are designed for two formulations: bilevel optimization and robust optimization. Both formulations are applied to the datasets used in the ATLAS A14 tune and to the dedicated hadronization datasets generated by the \sherpa generator, respectively. The corresponding tuned generator parameters are compared using three metrics. We compare the quality of our automatic tunes to the published  ATLAS A14 tune. Moreover, we analyze the impact of a pre-processing step that excludes data that cannot be described by the physics models used in the MC event generators.

\vspace{10pt}
\noindent\rule{\textwidth}{1pt}
\tableofcontents\thispagestyle{fancy}
\noindent\rule{\textwidth}{1pt}
\vspace{10pt}

\section{Introduction and Motivation}

Monte Carlo (MC) event generators are simulation tools that  predict  the properties of high-energy particle collisions.   
Event generators are built from theoretical formulae and models that describe the probabilities for various sub-event  phenomena that occur in a high-energy collision.
They are developed by physicists as a bridge between particle physics perturbation theory, which is defined at very high energy scales, and the observed sub-atomic particles, which are low-energy states of the strongly-interacting full theory.    This bridge is essential for interpreting event collision data in terms of the fundamental quantities of the underlying theory.   See \cite{Buckley:2011ms} for an overview of the event generators used for physics analysis at the Large Hadron Collider (LHC).

The description of particle collisions requires an understanding of phenomena at many different energy scales.
At high energy scales (much larger than the masses of the sub-atomic particles), first principle predictions can be made in a perturbative framework based on a few universal parameters.
At intermediate energy scales, an approximate perturbation theory can be established that introduces less universal parameters.
At low energy, motivated, but subjective, models are introduced
to describe sub-atomic particle production.   These low-energy models
 introduce a large number of narrowly defined parameters.
To make predictions or inferences, one must have a handle on the preferred models and the values of the parameters needed to describe the data.
This process of adjusting the parameters of the MC simulations to match data is called {\it tuning}.

This tuning task is complicated by the fact that the phenomenological models cannot claim to be complete or scale-invariant.
When compared to a large set of collider data collected in different energy regimes, the MC-models do not describe the full range of event properties equally well.
Typically, the physicists demand a tune that describes a subset of the data very well, another subset moderately well, and a remainder that must only  be
described qualitatively.
This distribution of subsets may well vary from one group of physicists to another and
has led to the education of experts who subjectively select and weigh data to achieve some physics goal.
Two such exercises are the Monash tune \cite{Skands:2014pea} and the A14 tune \cite{ATL-PHYS-PUB-2014-021}, though others exist in the literature.
Both of these tunes are successful, in the sense that they have been useful in understanding a wide range of phenomena observed at particle colliders.
However, the current approach to tuning remains inefficient and biased.

This work introduces a framework that, once agreed upon, greatly reduces the subjective element of the tuning process and
replaces it with an automated way to select the data for parameter tuning.

\subsection{Notation and terminology}

The data used in the tuning process are in the form of observables, denoted by $\calO$, and the set of observables is denoted by $\calS_\calO$.
Observables are quantities constructed from the (directly or indirectly) measured sub-atomic particles produced in an event.
In this case, each observable is presented as a \emph{histogram} that shows the frequency that the observable is measured over a range of possible values.   The range can be one or many divisions of the interval from the minimum to the maximum value that the observable can obtain.   These divisions are called \emph{bins}.    In practice, the size of a bin is set by how well an observable can be measured. The number of bins of an observable $\calO$ is denoted as \abscalO. 
We use $\calR$ to denote the reference data in the  histograms, a subscript $b$ to denote a bin, $\calR_b$  to denote the data value in a bin, and $\Delta\calR_b$ to denote the corresponding measurement uncertainty.

The MC-model has parameters $\p$, a $d$-dimensional vector in the  space $\Omega$, $\p\in\Omega\subset\mathbb{R}^d$.
The  MC-based simulations are denoted by $\MC(\p)$ to emphasize that they depend on the physics  parameters $\p$. The histograms computed from the MC simulation have the same structure as the histograms obtained from the measurement data $\mathcal{R}$, with a prediction per bin $\MC_b(\p)$ and an uncertainty associated with each bin $\Delta\MC_b(\p)$.   The uncertainty on the MC simulation comes from the numerical methods used to calculate the predictions, and it typically scales as the inverse of the square root of the number of simulated events in a particular bin.

\begin{figure}[htp!]
    \centering
    \includegraphics[width=.48\textwidth]{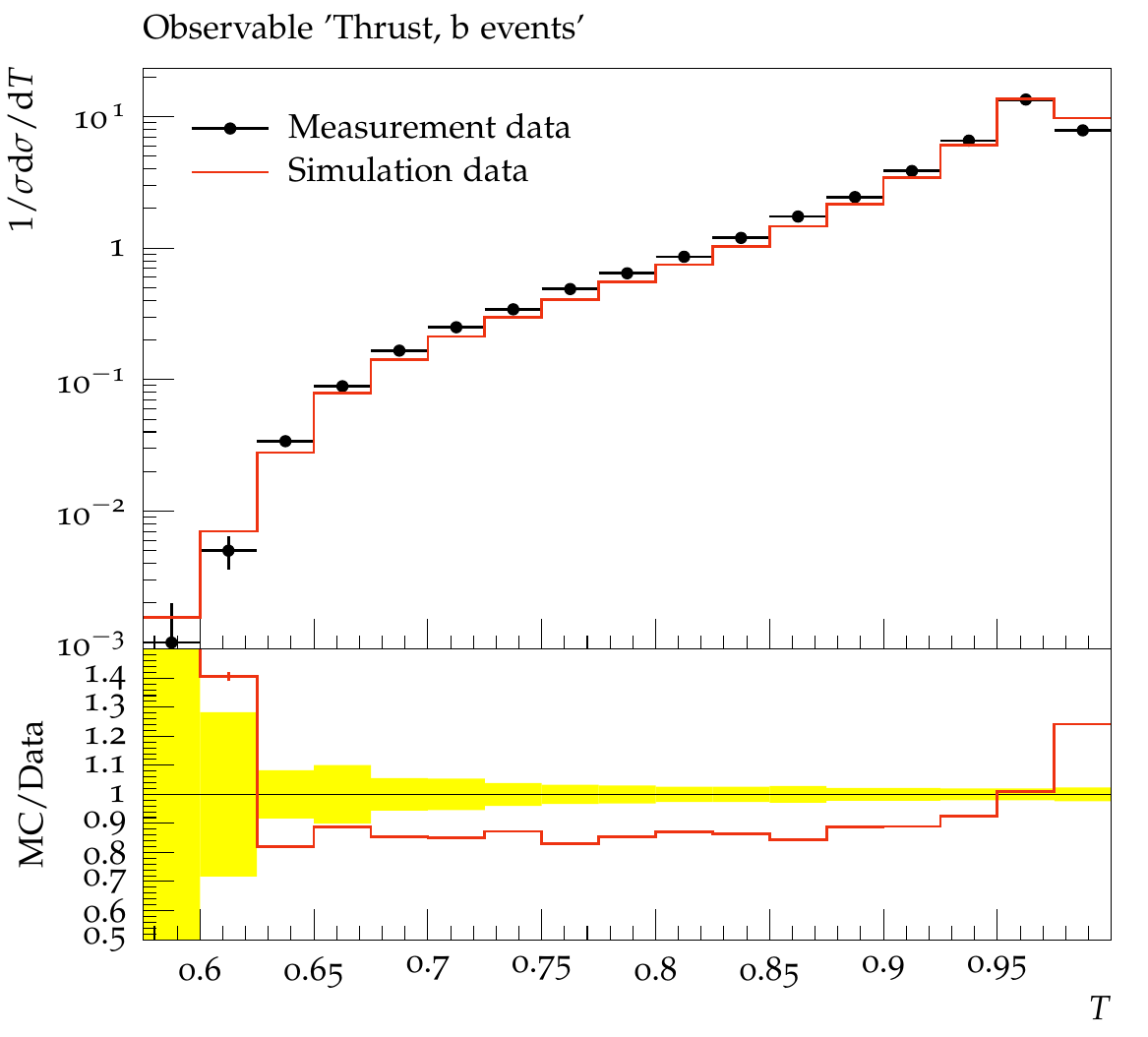}
    \caption{A histogram of a typical observable used in the tuning process.
    The top pane displays the measured (black) and predicted (red) data and their uncertainties.
    The bottom pane displays the ratio of predictions to measurements.   The yellow band displays the measurement uncertainty on the reference data: $\left[1 -\Delta\mathcal{R}_b/\mathcal{R}_b, 1 + \Delta\mathcal{R}_b/\mathcal{R}_b\right]$.}
    \label{fig:hist0}
\end{figure}

Figure~\ref{fig:hist0} shows a typical histogram. 
In this example, the observable, \textit{Thrust},  has 17 bins.   In the top pane, the black segments show the experimental data $\mathcal{R}$. The vertical error bars show the uncertainty associated with the data, i.e., $\Delta\mathcal{R}$. The red line shows the data obtained from the MC simulation $\MC(\p)$  with some parameter setting $\p$.
The bottom pane shows the ratio of $\MC(\p)$ to the data in each bin. The black horizontal line shows the reference ratio value one, to make the visual inspection easier. When the red line is above the black line, it means $\MC(\p)>\mathcal{R}$, and vice versa. The yellow region is defined by the range of the uncertainty on a measured value  (usually the 68\% confidence level on the reported value) relative to the measured value, i.e. $\left[1 -\Delta\mathcal{R}_b/\mathcal{R}_b, 1 + \Delta\mathcal{R}_b/\mathcal{R}_b\right]$.
A ``good'' tune is one where the red line falls within the yellow band.   In the example Figure~\ref{fig:hist0}, $\MC(\p)$ underpredicts the number of events with intermediate values of {\it Thrust} and overpredicts near the endpoints.

 \subsection{Mathematical formulation of the tuning problem}
 Our goal is to find an optimal set of physics parameters, $\p^*$, that  minimizes the difference between the experimental data and the simulated data from an MC event generator.
 This difference is defined as follows: 
\begin{equation}
\chisq_{\MC}(\p,\w) = \sum_{\calO\in\calS_\calO}w_\calO \sum_{b\in\calO}\frac{(\MC(\p)-\calR_b)^2}{\Delta \MC_b(\p)^2+\Delta\calR_b^2}, \label{eq:MC}
\end{equation}
where $w_\calO$ is the weight for an observable $\calO$ and 
$\w$ is a vector of weights, 
$\w=[w_1, \ldots, w_{|\calS_\calO|}]^T$.
In general, the number of bins can be different for different observables.  %
The weights $w_\calO\geq0$ reflect  how much an observable contributes to the tune, i.e., if $w_\calO=0$ for some $\calO$, then this observable will not influence the  tuning of $\p$.

The MC simulation is computationally expensive (the generation of 1 million events for a given set of parameters consumes about 800 CPU minutes on a typical computing cluster), severely
limiting the number of parameter choices $\p$ that can be used in the tuning.   %
To overcome these issues, we construct a parameterization of the MC simulation using an evolution of the \professor
framework~\cite{Buckley:2009bj}, named \apprentice. The code of \apprentice is available at \url{https://github.com/HEPonHPC/apprentice}.
The function in~\equaref{eq:MC} is not minimized directly.   Instead, during the optimization over $\p$, the MC simulation is replaced by a surrogate  model (here, a polynomial or a rational approximation to a number of MC simulations).
For each bin $b$ of each histogram, the central value and the corresponding uncertainty of the model prediction are parameterized independently as functions of the model parameters  $\p$.
This approach results in analytic expressions $f_b(\p)$ and $\Delta f_b(\p)$  that approximate the central value and the uncertainty, respectively,  and that can be evaluated in milliseconds. Thus, instead of~\equaref{eq:MC}, we minimize
\begin{equation}
\chisq(\p,\w) = \sum_{\calO\in\calS_\calO}w_\calO \sum_{b\in\calO}\frac{(f_b(\p)-\calR_b)^2}{\Delta f_b(\p)^2+\Delta\calR_b^2}. \label{eq:poly}
\end{equation}

Because $f_b(\p)$ and $\Delta f_b(\p)$ are given analytically and provide a prediction for any choice of $\p$ within
our domain of interest, we can minimize \equaref{eq:poly} efficiently using numerical methods.
\equaref{eq:poly} implicitly assumes that each bin $b$ is completely independent of all other bins.

In practice, the weights $w_\calO$ in~\equaref{eq:poly}  are adjusted manually, based on experience and physics intuition: the expert fixes the weights and minimizes the function in~\equaref{eq:poly}  over the parameters $\p$.   If the fit is unsatisfactory, a new set of weights is selected, and the optimization over $\p$ is repeated until the tuner is satisfied.\footnote{For the A14 tune, this selection was based on looking at hundreds of histograms such as the one shown in  \figref{fig:hist0}.}
 The selection of weights is time-consuming and different experts may have different  opinions about how well each observable is approximated by the model.
 Our goal is to automate the weight adjustment, yielding a less subjective and less time-consuming process to find the optimal physics parameters $\p$ that will then be used in the actual MC simulation.
 This problem was also considered in \cite{Bellm:2019owc},
 where weights are assigned based on how influential data is on constraining parameters.
 Also related to this work is that of \cite{Ilten:2016csi},
 which treats tuning as a black-box optimization problem within the framework of Bayesian optimization,
 but with no weighting of data.
 
 For convenience, we summarize our notation in Table~\ref{tab:notation}.

\begin{longtable}{p{.15\textwidth} | p{.85\textwidth}}
	\caption{Notation.}\label{tab:notation}   \\
	\hline
Notation & Definition \\\hline
	\endfirsthead
	\hline
	Notation & Definition \\\hline
	\endhead
	$\calO$ & observables that are constructed from data and MC-based simulations in the form of histograms\\
        \abscalO & the number of bins in an observable $\calO$\\
        $\calS_\calO$ & the set of observables used in the tune \\
    $|\calS_\calO|$ & the number of observables\\
        $\calR$ & the data in the histograms \\
        $b$ & a bin of a histogram  $\calO$ \\
        $\calR_b$  &  the data value in a bin \\
        $\Delta\calR_b$ & data uncertainty corresponding to the data value in a bin \\
        $\p$ & a $d$-dimensional vector of real-valued parameters \\
        $\MC(\p)$ & an  MC simulation that depends on the physics parameters $\p$ \\ 
        $\MC_b(\p)$ & the MC simulation in a bin $b$  \\
        $\Delta\MC_b(\p)$ & an uncertainty associated with the MC simulation in a bin $b$ \\
        $f_b(\p)$ & central value of the model prediction parameterized independently as a function of the model parameters  $\p$ \\
        $\Delta f_b(\p)$  & the uncertainty of the model prediction parameterized independently as a function of the model parameters  $\p$ \\        
        $\w$ &  an $|\calS_\calO|$-dimensional vector of real-valued weights \\
        $w_\calO$ & the weight given to a histogram in constructing a tune (if $w_\calO=0$ for some $\calO$,
        then this observable will not influence the  tuning of $\p$). \\
        $\phatw$ & optimal physics parameters for a given choice for the weights \\
        $\w^*$ &  an optimal set of weights for the observables \\
        $\phatws$ & the optimal set of simulation parameters corresponding to an optimal set of weights $\w^*$ for the observables \\
        $g$ & the outer objective function of $\mathbb{R}^{|\calS_\calO|\times d}\mapsto\mathbb{R}$  used in the bilevel optimization\\        
        $\mu$ & a hyperparameter that specifies the percentage of the observables used in the robust optimization\\
        $\chisq_\calO(\p)$ & the per-observable error averaged over all bins in the observable $\calO$\\
        $\p_\text{ideal}^\calO$ & the \textit{ideal} tune for an observable $\calO$, i.e., the parameters that minimize~\equaref{eq:obschi2} when using only observable $\calO$ for the tune\\ \hline	
 \end{longtable}

\section{Finding the Optimal Weights for Each Observable }
In this section, we describe two mathematical formulations for finding the optimal weights in \equaref{eq:poly} that determine how much influence each observable should have on the optimization over the physics parameters $\p$:
bilevel and robust optimization.

\subsection{Bilevel optimization formulation}\label{sec:bilevel}

We formulate a bilevel optimization problem as follows:
\begin{subequations}\label{eq:bilevel_whole}
\begin{align}
\min_{\w \in[0,1]^{|\calS_\calO|}, \hat{\p}_{\w}\in\Omega} &  g(\w, \widehat{\p}_{\w})\label{eq:outer_gen}\\
\text{s.t. } &  \sum_{\calO\in\calS_\calO} w_\calO=1\label{eq:sum1_bilevel}\\
& \widehat{\p}_\w \in \arg\min_{\p\in\Omega} \chisq(\p, \w)\label{eq:inner_bilevel}
\end{align}
\end{subequations}
where the function $g:\mathbb{R}^{|\calS_\calO|\times d}\mapsto\mathbb{R}$ describes a merit function 
to determine the goodness of weights (see below for the definitions we use in this work). 
The lower-level \equaref{eq:inner_bilevel} (same as~\equaref{eq:poly}) corresponds to finding the optimal parameters $\widehat{\p}_\w$ for a given set of weights $\w$, and the upper-level \equaref{eq:outer_gen} provides a measure of how good the weights are. The weights are normalized to sum to unity, see \equaref{eq:sum1_bilevel}, in order to prevent the trivial solution where  all weights are 0. 
Bilevel optimization problems have been studied extensively in the literature, see, e.g., \cite{ChenFlor:95,MarZhu:96,YeZhuZhu:97,colson2007overview,bard2013practical}.

In the following, we discuss two definitions of the outer objective function $g(\w,\phatw)$.  Other formulations are possible and our selection is driven by the goal to achieve reasonably
good agreement between the simulated and the observed data for all observables (rather than fitting a few observables extremely well and others poorly).

\subsubsection{Formulation 1: Portfolio to balance mean and variance of errors}\label{sec:bilevel-port}
The portfolio objective function is motivated by portfolio optimization in finance~\cite{portfolio}, where the  goal is to maximize the expected return while  minimizing  the risk. Translated to our problem, we want to minimize the expected error over all observables while also minimizing the variance over these errors. 

For a given set of weights $\w$, %
we obtain the ``$\w$-optimal'' parameters $\phatw$. For each observable $\calO$, an  error term is averaged over the number of bins in the observable ($\abscalO$):
\begin{equation} 
  e_{\cal O}( \phatw | \w) = \frac{1}{\abscalO}\sum_{b \in {\cal O}} \frac{\left(f_b(\phatw) - {\cal R}_b \right)^2}{\Delta f_b(\phatw)^2 + \Delta{\cal R}_b^2}, \ \calO\in\calS_\calO,
\end{equation}	
 where we write $e_{\cal O}( \phatw |\w)$ because the error value for each observable depends (implicitly) on the choice of the weights $\w$.    Thus, we obtain a set of $|\calS_\calO|$ average error values from which we compute the following statistics:

 \begin{subequations}
 \begin{equation} 
 \mu(\phatw|\w) = \frac{1}{|\calS_\calO|}\sum_{\calO\in\calS_\calO} e_\calO(\phatw|\w) \text{: average error over all observables,}
\end{equation}	
  
 \begin{equation} 
 \sigma^2(\phatw|\w) = \frac{1}{|\calS_\calO|}\sum_{\calO\in\calS_\calO} \left[e_\calO(\phatw|\w) -\mu(\phatw|\w)\right]^2 \text{: empirical variance of errors over all observables.}
\end{equation}
\end{subequations}
 
The portfolio objective function for the outer optimization then becomes 
 \begin{equation} 
g(\w,\phatw)=  \mu(\phatw|\w) + \sigma^2(\phatw|\w), \label{eq:portfolio}
\end{equation}	
which represents a simultaneous  minimization of the expected error \textit{and} the variance of the errors. For problems in which minimizing the variance is of higher priority, one can introduce a multiplier $\lambda$ before the variance term that reflects ``risk aversion''. In that case, if $\lambda$ is large,  we are more risk-averse, since reducing the variance associated with the errors will drive the minimization. If $\lambda$ is small, we are less risk-averse, and minimizing the mean of the errors is emphasized.

\subsubsection{Formulation 2:  Scoring of model fit and data uncertainty} \label{sec:bilevel-score}
We consider a second outer objective function formulation based on scoring schemes~(\cite[Eq.~(27)]{gneiting2007strictly}). 
The performance of a generic predictive model $P$  at a point $x$ is defined by a scoring rule, $S(P,x) = -\left(\frac{x-\mu_P}{\sigma_P}\right)^2 - \log \sigma_P^2$,
where $P$ has mean performance $\mu_P$ and variance $\sigma_P^2$.  A larger value for $S(P,x) $ signifies better model performance.
Thus, we minimize the negative  of $S(P,x)$:  
\begin{equation}\label{eq:score_oo}
    s(P,x) = - S(P,x) =  \left(\frac{x-\mu_P}{\sigma_P}\right)^2 + \log \sigma_P^2.
\end{equation}
For our application, $x$ corresponds to the simulation prediction $f_b(\p)$,  $\mu_P$ to our observation data $\calR_b$,  and the variance $\sigma^2_P$
to our data uncertainty $\Delta\calR_b$.
For each bin $b$ in an observable, we  calculate the score based on~\equaref{eq:score_oo}.
Then, we compute  the median (and mean) of the scores over all bins to obtain the median (average) performance for each observable.
In order to form the upper-level objective function, we sum up the median (mean) scores over all observables:

\begin{itemize}
\item{Outer objective based on median score}
\begin{subequations} 
\begin{equation}
    g(\w,\widehat\p_\w)=  \sum_{\calO\in\calS_\calO}\tilde{s}_\calO( \widehat{\p}_{\w}|\w), \end{equation}
\begin{equation}
\tilde{s}_\calO(\widehat{\p}_{\w}|\w)=\text{median of } \Bigg\{\left(\frac{f_b(\widehat{\p}_{\w})-\calR_b}{\Delta\calR_b}\right)^2 + \log(\Delta\calR_b^2), \forall {b\in\calO}\Bigg\}.\label{eq:medianscore}
\end{equation}
\end{subequations}

\item{Outer objective based on mean score}

  \begin{subequations}
\begin{equation}
    g(\w,\widehat\p_\w)=  \sum_{\calO\in\calS_\calO}\bar{s}_\calO( \widehat{\p}_{\w}|\w), \end{equation}
\begin{equation}
\bar{s}_\calO(\widehat{\p}_{\w}|\w)=
\frac{1}{\abscalO}\sum_{b \in {\cal O}} \Bigg\{\left(\frac{f_b(\widehat{\p}_{\w})-\calR_b}{\Delta\calR_b}\right)^2 + \log(\Delta\calR_b^2)\Bigg\}.
\label{eq:meanscore}
\end{equation}
\end{subequations}
\end{itemize}

In our numerical experiments, we analyze and compare both the performance of  the median score and the mean score. 
 Both the median and the mean score outer objective functions take into account the deviation of the prediction of $f_b(\widehat{\p}_\w)$ from $\calR_b$ and the uncertainty in the data $\Delta\calR_b$. Thus, if an observable has large uncertainties in the data or the model $f_b(\widehat{\p}_\w)$ does not approximate the data $\calR_b$ well, the score for this observable deteriorates. 
 Ideally, both terms $\left(\frac{f_b(\widehat{\p}_{\w})-\calR_b}{\Delta\calR_b}\right)^2$ and $\log(\Delta\calR_b^2)$ will be small.

 \subsubsection{Solving the bilevel optimization problem using surrogate models}
 Solving the inner optimization problem~(\ref{eq:inner_bilevel}) for each weight vector $\w$ is generally computationally non-trivial and its computational demand increases with the number of physics parameters $\p$ that have to be optimized and the number of observables present. Therefore, the goal is to try as few weights $\w$ as possible. We interpret the solution of the inner optimization problem as a black-box function evaluation of $g(\w, \hat{\p}_{\w})$ for $\w$. Given an initial set of  input-output data pairs $\left\{(\w_i, g(\w_i, \hat{\p}_{\w_i})\right\}_{i=1}^I$, we fit a surrogate model\footnote{This surrogate model for the weights is independent of the one used to evaluate the MC-based predictions.} (here a radial basis function~\cite{Powell1992}) that allows us to predict the values of $g(\w, \hat{\p}_{\w})$ at untried $\w$. In each iteration of the optimization algorithm, these predictions are used to select the most promising weight vector for which the inner optimization problem should be solved next. Promising weight vectors have either low predicted values of $g(\cdot)$ or are far away from already evaluated points~\cite{Muller2014, Mueller2017}.  Each time a new weight vector has been evaluated, the surrogate model is updated. This iterative process repeats until a stopping criterion has been met, e.g., a maximal number of weight vectors has been evaluated or a maximal CPU time has been reached. Details about the surrogate model algorithm are given in the online supplement Section~\ref{sec:surrogate}.
 
When solving the problem~\equaref{eq:bilevel_whole}, the optimization algorithm for the outer problem  selects a set of weights $\w$.
Given a choice for the weights,
we solve the inner optimization problem \equaref{eq:inner_bilevel} using \apprentice to obtain a set of optimal physics parameters $\phatw$. %
Given $\phatw$, we compute the corresponding function value of the outer objective function,  $g(\w, \phatw$). Based on this value, the outer optimization algorithm selects a new set of weights, which will be used to solve the inner optimization problem again. This leads to a new solution for \equaref{eq:inner_bilevel}, which in turn gives a new value for the outer objective function.
This process repeats until the outer optimization converges to an  optimal set of weights for the observables (denoted by $\w^*=[w_1^*, \ldots, w_{|\calS_\calO|}^*]^T$) and a corresponding optimal set of simulation parameters (denoted by $\phatws$).

\subsection{Robust optimization formulation}
\label{sec:robust}

As an alternative to the bilevel formulation, we developed a single-level robust optimization formulation for finding the optimal weights for~\equaref{eq:poly}. 
Robust optimization estimates the parameters ${\p}$ that
minimize
the largest deviation $\left(f_b(\p) - \mathcal{R}_b\right)^2$ over
all bins in an uncertainty set $\mathcal{U}_b$ of bin $b$:
\begin{align}
\label{eq:robo_f}
\begin{aligned}
\underset{\w\in[0,1],\p\in\Omega}{\mini}~& \sum_{\mathcal{O} \in \calS_\calO} \frac{w_\mathcal{O}}{|\calO|} \sum_{b \in \mathcal{O}} \underset{\mathcal{R}_b\in\mathcal{U}_b}{\maxi}\left(f_b(\p) - \mathcal{R}_b \right)^2.  \\
\end{aligned}
\end{align}
Assuming that the experiment and the MC simulation are described using
independent random variables with mean $\mathcal{R}_b$, the
uncertainty set $\mathcal{U}_b$ for each bin $b$ is described by the
interval $[\mathcal{R}_b-\Delta\calR_b-\Delta
f_b(\p),\mathcal{R}_b+\Delta\calR_b +\Delta f_b(\p)]$.

Introducing slack variables $\bm{t} = \left[ t_1,t_2,\dots,t_{|\calO|} \right]$, we rewrite \eqref{eq:robo_f} as:
\begin{subequations}\label{eq:robustopt_whole}
\begin{align}\label{eq:robo_fexp_obj}
\underset{\bm{t},\w\in[0,1],\p\in\Omega}{\mini}~& \sum_{\mathcal{O} \in \calS_\calO} \frac{w_\mathcal{O}}{|\calO|} \sum_{b \in \mathcal{O}} t_b \\
\begin{split}
\label{eq:robo_fexp_constr}
  \st &\\
t_b &\ge \left(f_b(\p) - (\mathcal{R}_b-\Delta\calR_b-\Delta f_b(\p)) \right)^2 \quad \forall b \in \mathcal{O}, \forall \mathcal{O} \in \calS_\calO\\
t_b &\ge \left(f_b(\p) - (\mathcal{R}_b+\Delta\calR_b+\Delta f_b(\p)) \right)^2 \quad \forall b \in \mathcal{O}, \forall \mathcal{O} \in \calS_\calO
\end{split}\\
\begin{split}
  \label{eq:robo_hyperparam_constr}
\sum_{\mathcal{O} \in \calS_\calO} \frac{w_\mathcal{O}}{|\calO|} &\ge 
 \frac{\mu}{100} \sum_{\mathcal{O} \in \calS_\calO} \frac{1}{|\calO|}
\end{split}
\end{align}
\end{subequations}
where the constraint \eqref{eq:robo_hyperparam_constr} is enforced to avoid the trivial solution of all weights being  zero. 
In this constraint, we bound the sum of the weights away from zero by a hyperparameter $\mu$ that  specifies the percentage of the observables that should be used in the  optimization.
Problem \eqref{eq:robustopt_whole} is attractive because it formulates the problem of finding optimal weights as a single-level optimization problem, which is easier to solve than the bilevel problem \equaref{eq:bilevel_whole}.

Selecting the best $\mu$ among all the 100 runs of robust optimization is determined using a cumulative density curve of the number of observables satisfying $\displaystyle \frac{\chisq_\calO (\p^*,\w)}{|\calO|} \le \tau$, where $\p^*$ is the optimal parameter obtained from the robust optimization run, $\w=\bm{1}$, $\tau \in \mathbb{R}^+$ and $\calO \in \calS_\calO$. 
Hence, in the plot of this curve (e.g., see Figure~\ref{fig:romuperformance}), the number of observables on the y-axis is monotonically increasing as $\tau$ increases on the x-axis. 
Then, the area between the cumulative density curve for each robust optimization run and the ideal cumulative density curve is computed. 
To build the ideal cumulative density curve, the $\p^*$ in $\displaystyle\frac{\chisq_\calO (\p^*,\w)}{|\calO|} \le \tau$ is obtained by considering only observable $\calO$ in \equaref{eq:poly}.
The best run is then chosen to be the one whose area to the ideal cumulative density curve is the smallest.
An example plot of the cumulative density curve and an illustration of the procedure to find the best run is included in Section~\ref{sec:ROmuselection} of the online supplement.

\section{Data Pre-processing: Filtering  Observables or Bins}

We also investigate the question of how to detect and exclude observables or bins whose data $\calR_b$ cannot be explained by the MC simulation model. 
One special choice of weight for an observable is $w_\calO=0$, which corresponds to excluding 
(filtering out) the observable $\calO$ from our parameter tune.
This is driven by a significant discrepancy between the simulation and data.  Such discrepancies can arise for at least two reasons:  (1) a mistake has been made in the experimental analysis; and/or (2) the observable is out of the domain of predictions that can be made reliably with the simulation. For our studies, we assume that the source of discrepancies is from (2). 
Because the simulation is a metamodel constructed from many smaller models, it is difficult to make
{\it a priori} statements about the domain of its predictions.  Important physics may be missing from the metamodel and/or a
model can describe the mean behavior but not the rarer fluctuations around the mean. 
The simulation should be able to describe the physics, but some observables worsen the description. Thus, it is quite reasonable to exclude these  observables.

In our discussion to this point, we have assumed that each {\it observable} has a given weight.   However, in those situations where the model can describe the mean behavior, it can be beneficial to filter out individual bins $b$ of the observable.   In the observables considered in this study, and typical of the high energy physics phenomenon, the models can have difficulties in describing the rise and/or fall of a distribution (consider the example
in Figure~\ref{fig:hist0} where there is a rise from the first  to the second bin and a fall from the penultimate to the last bin and the corresponding predicted data are far away from the measured, indicated by the red line in the lower pane.).

\subsection{Filtering of observables by outlier detection} \label{sec:obsfilter}

Using the surrogate model $f_b(\p)$ to approximate the expensive MC simulation, we can efficiently minimize the per-observable-$\chisq$ function:
\begin{equation}
\chisq_\calO(\p) =  \frac{1}{|\calO|}\sum_{b\in\calO}\frac{(f_b(\p)-\calR_b)^2}{\Delta f_b(\p)^2+\Delta\calR_b^2}, \calO \in\calS_\calO \label{eq:obschi2}
\end{equation}
for each observable $\calO \in\calS_\calO$, separately. %
$\chisq_\calO(\p)$ represents the average per-bin error for the observable and the best possible fit of the model for this single observable.
If we used  only one observable for the tune, the parameters $\p_\text{ideal}^\calO$ that minimize~\equaref{eq:obschi2} would represent the \textit{ideal} tune. The corresponding ideal objective function value $\chisq_\calO(\p_\text{ideal}^\calO)$ is the best possible result  for each individual observable $\calO$. Because the ideal parameter values will be different for each observable, we will not be able to obtain one parameter set that minimizes~\equaref{eq:obschi2} for all observables simultaneously. Therefore,  we obtain a set $\mathcal{X}$ of length $|\calS_\calO|$ of ideal objective function values of~\equaref{eq:obschi2}: $\mathcal{X}:=\{\chisq_{1} (\p_\text{ideal}^1),\chisq_{2} (\p_\text{ideal}^2), \ldots,\chisq_{|\calS_\calO|} (\p_\text{ideal}^{|\calS_\calO|})\} = \{\chi_{i}\}_{i=1}^{|\calS_\calO|}$.
If  the ideal error is large for some observables, it means that the model is not able to fit the data of these observables well at all, even with the freedom of not having to fit any other observables. Therefore, the inclusion of these data in optimizing~\equaref{eq:poly} may negatively impact the
overall optimization because large errors might drive the optimization.\footnote{We address later the fidelity of the
  surrogate model.}

To address this issue, we use the distribution of the values in $\mathcal{X}$ to identify outliers (observables with  values for ~\equaref{eq:obschi2} ``that deviate so much from other observations as to arouse suspicions that it was generated by a different mechanism.''\cite{hawkins1980identification}). We exclude the outlier  observables from the optimization of~\equaref{eq:poly} by setting their corresponding weights to zero, $w_\calO = 0$.

There are multiple methods that can be used for outlier detection, such as scatter plots \cite{croarkin2006nist}, Z-score \cite[Section 1.3.5.17]{croarkin2006nist}, interquartile range \cite{upton1996understanding}, generalized extreme studentized deviate~\cite{rosner1983percentage}, Grubb's test \cite{grubbs1969procedures, stefansky1972rejecting}, Dixon's Q test \cite{dixon1953processing}, Thompson tau test \cite{thompson1985note}, Pierce’s Criterion \cite{dardis2004peirce}, and Tietjen-Moore test \cite{tietjen1972some}, to name a few. 
We obtained reasonable results using the
Z-score.   For the set  $\mathcal{X} = \{\chi_{i}\}_{i=1}^{|\calS_\calO|}$, the Z-score of an observation $\chi_{i}$ is defined as $z_{i} = (\chi_{i} - m)/s$ where $m$ is the mean of the observation set $\mathcal{X}$ and $s$ is the standard deviation.
We calculate the Z-score for each data point $i$ in $\mathcal{X}$ and define an outlier as $z_i\ge 3$.
In other words, any ideal fit with a residual outside of 3 standard deviations is classified as an outlier.

The benefit of performing the outlier detection is that the computational cost of minimizing~\equaref{eq:poly} is reduced. %
In addition, the optimization will not be biased by observables that the underlying model cannot describe well.
\subsection{Filtering of  bins by hypothesis testing} \label{sec:binfilter}
We explore a second and more refined approach that allows us to  identify and exclude bin data
from the optimization of~\equaref{eq:poly}. Instead of eliminating whole observables, we identify a subset of bins for each observable that cannot be approximated well by the MC simulator model and we exclude only those bins from the optimization~\cite{ETZIONI1988416}.
The motivation behind this idea is driven by physics. In many cases, the physics observables are constructed such that one or both extrema of the observables' bins in the histogram represent high or low energy regimes. At low energy, a phenomenological model may be entirely inadequate. At high energy, genuine perturbative corrections not included in our calculations
may be needed. It would be non-physical to adjust the model parameters to explain these extremes.

To this end, we use the $\chisq$ test, which is a hypothesis test performed when the test statistic is $\chisq$-distributed under the null hypothesis~\cite{cochran1952}. Note that the $\chisq$ test statistic is different from the $\chisq(\p,\w)$ objective function introduced earlier. %
We first compute the $\chisq$ test statistic for a subset $\calB$ of the bins  in an observable $\calO$ using the computationally cheap approximation model $f_b(\p)$:
\begin{equation}\label{eq:subsetbinchi2}
\chisq_\calB(\p) =  \frac{1}{|\calB|}\sum_{b\in\calB\subset\calO}\frac{(f_b(\p)-\calR_b)^2}{\Delta f_b(\p)^2+\Delta\calR_b^2}.
\end{equation}
For this statistic, we hypothesize that: 
\begin{itemize}
	\item[] Null hypothesis $H_0$: There is no significant difference between $f_b(\p)$ and $\mathcal{R}_b$, in other words,  the data $\mathcal{R}_b$ can be appropriately described by $f_b(\p)$. 
	\item[] Alternate hypothesis $H_1$:  $H_0$ is rejected, i.e., there is a significant difference between $f_b(\p)$ and $\mathcal{R}_b$.	
\end{itemize}

In \eqref{eq:subsetbinchi2}, we have a sample of size $|\calB|$ based on which we compute the $\chi^2$ test statistic. However, the degrees of freedom of the $\chisq$ distribution is not $ |\calB|$ because the samples $f_b(\p), b\in\calB\subset\calO$ are not independent and they are related to each other through the parameters $\p$. 
Due to this relationship,  the number of degrees of freedom is reduced (see~\cite{alex2019goodnessoffit} for a similar argument). 
Hence the resulting degrees of freedom of the $\chisq$ distribution for the set $\calB$ is given by 
\begin{equation}\label{eq:degoffreedom}
\rho_\calB =  |\calB| - d,
\end{equation}
where $d$ is the dimension of $\p$.

We now choose a value for the significance level $\alpha$ (commonly used values are  0.01, 0.05, or 0.1). From a $\chi^2$ distribution table, we then obtain the critical value $\chi^2_{c,\calB}$ for bins in $\calB$ as a function of the significance level $\alpha$ and degrees of freedom $\rho_\calB$. 
More formally, we say that if the probability $P_{H_0} (T \le \chi^2_{c,\calB}) = \alpha$, then under $H_0: T \sim \chi^2 (\rho_\calB)$. Let us assume a random variable $Z \sim \chi^2 (\rho_\calB)$, then $P(Z \le \chi^2_{c,\calB}) = \alpha$. Thus, to find $\chi^2_{c,\calB}$, we need to compute the inverse of the cumulative distribution  function (CDF) of the $\chi^2$ distribution with $\rho_\calB$ degrees of freedom and at level $\alpha$.
Then we compare the test statistic with the critical value to decide whether $H_0$ is accepted or not, i.e., if $\chi^2_\calB \le \chi^2_{c,\calB}$, we keep the bin subset $\calB$; otherwise, we cannot keep this bin subset. 

We mainly intend to  exclude bins at the extremes of the observables, and hence we require that the bins we keep are contiguous. 
For some observables all  bins may  pass the $\chisq$ test, for others, all bins may be excluded, or  a subset of contiguous bins  is kept.

The problem is then to find the largest contiguous subset of bins $\calB$ such that $\chi^2_\calB \le \chi^2_{c,\calB}$. This is equivalent to solving the mixed-integer program
\begin{equation}\label{eq:hypotestmip}
\begin{aligned}
\max_{s,e\in\{1,2,\dots,\abscalO\}} & \quad e-s \\
\text{s.t. } &\chi^2_\calB \le \chi^2_{c,\calB}, \quad \calB = \{s,\dots,e\},
\end{aligned}
\end{equation}
where $s$ and $e$ are the start and end indices of contiguous bins in observable $\calO$. This problem can also be solved using a polynomial-time algorithm based on the maximum sub-array problem~\cite{10.1145/358234.381162}.  This algorithm is described in Section~\ref{sec:polyalghyptest} in the online supplement. 
In some cases, the bins to keep may not be unique, i.e., there may be multiple ranges of \{$s,\dots,e$\} that are of the same maximum length and satisfy the null hypothesis (or satisfy the constraint in~\equaref{eq:hypotestmip}). In practice, this is not a problem, since selecting any one of these bin subsets does not change the outcome of the filtering or the optimization in~\equaref{eq:poly}.

\section{Numerical Experiments and  Comparison of Different Tunes}

In this section, we describe the setup of our numerical experiments, the datasets we use in our study, and the results. Additional information can be found in the online supplement.
\subsection{Setup of the numerical experiments}

\begin{table}
\caption{Optimization methods used in this study.}
\label{tab:methods}
\begin{adjustbox}{width=\textwidth}
\begin{tabular}{l|l|l}
\hline
  Name & Methodology & Reference \\\hline
	 ``Bilevel-portfolio''& bilevel optimization with portfolio outer objective function & Section~\ref{sec:bilevel-port}. \\
	 ``Bilevel-medianscore''&  bilevel optimization with median score outer objective function & Section~\ref{sec:bilevel-score}. \\
	 ``Bilevel-meanscore''& bilevel optimization with mean score outer objective function & Section~\ref{sec:bilevel-score}. \\
	 ``Robust optimization''& single level robust optimization approach & Section~\ref{sec:robust}. \\
	 ``Expert''& weight adjustment done by the expert (only for the A14 dataset, see Section~\ref{sec:A14})&  ~\cite{ATL-PHYS-PUB-2014-021} \\
	 ``All-weights-equal''& no optimization is used and all observable weights are set to 1 & \\\hline
\end{tabular}
\end{adjustbox}
\end{table}

We compare the results of using the  methods shown in Table~\ref{tab:methods} for  adjusting the weights of the observables in our datasets.
The performance of each method is evaluated
with and without data pre-processing (observable-filtering and bin-filtering approaches, see Sections~\ref{sec:obsfilter} and~\ref{sec:binfilter}),
and when using a cubic polynomial (results presented in the online supplement) versus a rational approximation for $f_b(\p)$ in \apprentice. 
We found relatively good performance using the degrees 3 and 1 for the numerator and denominator polynomial, respectively, for the rational approximation.

For the bilevel optimization formulation (see \equaref{eq:bilevel_whole}), we made the following choices: The initial experimental design for the outer optimization has $|\calS_\calO| + 1$ points, where $|\calS_\calO|$ is the number of observables (number of weights to be adjusted) included. The total number of allowed outer objective  function evaluations (number of weight vectors tried) is 1000.  Because the inner optimization function is multimodal, we use 100 multi-starts with \apprentice to solve it. The bilevel optimization with each method (portfolio, meanscore, medianscore) is repeated three times with different random seeds and we report the results of the best run.

For the robust optimization formulation (\equaref{eq:robustopt_whole}), a total of 100 random values of $\mu \in (0,100]$ are used when evaluating \equaref{eq:robo_hyperparam_constr} and, for each $\mu$, the algorithm is run once. 
The best run amongst these is returned as the best $\mu$ for the robust optimization. The procedure to select the best $\mu$ is described in Section~\ref{sec:robust}.

\subsection{Comparison metrics and optimal tuning parameters}
\label{sec:metrics}
There are many ways to assess the quality of a tune.   In many cases, the domain experts visually inspect
a potentially large number of histograms (see, e.g., Figure~\ref{fig:hist0}) to make a judgment.   As an
objective measure, we propose three metrics, each represented as a single number for each tuning method,  that can be used  to compare the quality of the  model fits obtained by the different methods in a more objective fashion: 
\begin{enumerate}
\setlength\itemsep{0em}
    \item \textit{Weighted $\chi^2$}: the sum over all  $\chi^2$ at the best $\phatws$,  \[ \sum_{\calO\in\calS_\calO}w^*_\calO \sum_{b\in\calO}\frac{(f_b(\phatws)-\calR_b)^2}{\Delta f_b(\phatws)^2+\Delta\calR_b^2}\]
    where $w^*_\calO$, the weight of observable $\calO$, is scaled such that $w^*_\calO \in [0,1]$ and $\sum_{\calO\in\calS_\calO} w^*_\calO=1$.
    \item \textit{A-optimality}: \[
    \operatorname{Tr}\left(\mathbf{\Gamma}_{\mathrm{post}}(\phatws,\w^*)\right) = \sum_{i=1}^{d} \lambda_i \]
    \item  \textit{$\log$ D-optimality}: \[
        \log \operatorname{det}\left(\boldsymbol{\Gamma}_{\text {post }}(\phatws,\w^*)\right) = \sum_{i=1}^{d}\log\lambda_i,
    \]
\end{enumerate}
where $\lambda_i$ are the eigenvalues of $\mathbf{\Gamma}_{\mathrm{post}}$, $\mathbf{\Gamma}_{\mathrm{post}}$ is the weighted posterior covariance matrix in the Bayesian formulation of the inverse problem, $d$ is the dimension of $\phatws$.
To find $\mathbf{\Gamma}_{\mathrm{post}}$, we compute the optimal parameter point $\phatws$, which is also referred to as the maximum a posteriori probability estimate in the context of Bayesian inverse problems~\cite{pronzato2013design}.
Given the optimal parameters, we can find a linearization of the model as
\begin{equation*}
   \mathbf{F}_\calO(\phatws) = \left[\frac{\partial f_{1}\left(\phatws\right)}{\partial \mathbf{p}},\frac{\partial f_{2}\left(\phatws\right)}{\partial \mathbf{p}},\dots,\frac{\partial f_{|\calO|}\left(\phatws\right)}{\partial \mathbf{p}}\right]^{\top} 
\end{equation*}
for each observable $\calO$.
Then the weighted posterior can be approximated by a Gaussian $\mathcal{N}\left(\phatws, \mathbf{\Gamma}_{\mathrm{post}}\right)$. Here
\begin{equation}\label{eq:gammapost}
    \mathbf{\Gamma}_{\mathrm{post}}(\phatws,\w^*) = \left(\sum_{\calO\in\calS_\calO} w^*_\calO \mathbf{F}_{\calO}^{\top}(\phatws) \mathbf{\Gamma}_{\text {noise}}^{-1} \mathbf{F}_{\calO}(\phatws)\right)^{-1} 
\end{equation}
where $\mathbf{\Gamma}_{\text {noise}}[\calO] = \operatorname{diag}\left(\Delta f_1(\phatws)^2+\Delta\calR_1^2, \Delta f_2(\phatws)^2+\Delta\calR_2^2, \ldots, \Delta f_{|\calO|}(\phatws)^2+\Delta\calR_{|\calO|}^2\right)$ and $w_\calO^*$ is the weight of observable $\calO$ obtained from the methods and is scaled such that $w^*_\calO \in [0,1]$ and $\displaystyle \sum_{\calO\in\calS_\calO} w^*_\calO=1$.

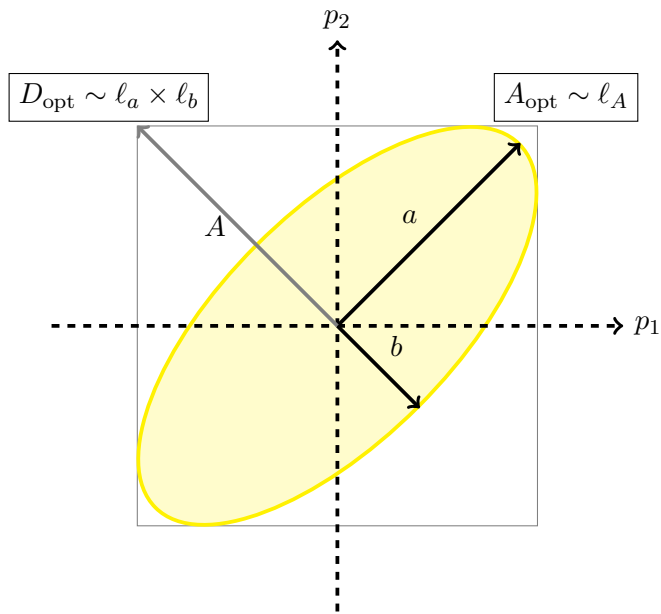
\begin{figure}[!h]
  \centering
  \begin{tikzpicture}[scale=0.76]
\draw [gray] (-3.5,-3.5) rectangle (3.5,3.5);    
\filldraw[yellow, fill opacity=0.2, rotate=45, line width = 0.5mm] (0, 0) ellipse (4.5cm and 2.0cm);
\draw[line width = 0.5mm, dashed, ->] (-5, 0) -- (5, 0) node[right]{$p_{1}$};
\draw[line width = 0.5mm, dashed, ->] (0, -5) -- (0, 5) node[above]{$p_{2}$};
\draw[line width = 0.5mm, ->] (0, 0) -- (3.2, 3.2) node[midway,above left]{$a$};
\draw[line width = 0.5mm, ->] (0, 0) -- (1.44, -1.44) node[midway,above right]{$b$};
\draw[gray,line width = 0.5mm, ->] (0, 0) -- (-3.5, 3.5) node[midway,black,below,left]{$A$};
\node[draw] at (-4,4) {$D_\textrm{opt} \sim \ell_a\times \ell_b$};
\node[draw] at (4,4) {$A_\textrm{opt} \sim \ell_A$};
\end{tikzpicture}
\caption{Graphical interpretation of optimal design criteria.
The D-criterion is the determinant of the covariance matrix or the volume of the joint confidence region.
The A-criterion is the trace of the covariance matrix or the dimension of a box that bounds the joint confidence region.
}
\label{fig:doecriteria}
\end{figure}

The $\mathbf{\Gamma}_{\mathrm{post}}(\phatws,\w^*)$ calculated at the optimal parameters and the optimal weights in \eqref{eq:gammapost} are used
here to describe the confidence region around the tuned parameters $\phatws$.  In order to summarize the multidimensional nature of $\mathbf{\Gamma}_{\mathrm{post}}$ into a scalar quantity, we use the A- and log D-optimality criteria. A graphical representation of the optimality criteria is shown in Figure~\ref{fig:doecriteria}.
The A-optimality criterion computes the trace of $\mathbf{\Gamma}_{\mathrm{post}}$, which is equivalent to the sum of its eigenvalues. This metric is proportional to the sum of the semiaxis lengths of the confidence ellipsoid of the parameters (lower is better), which corresponds to the average sum of the variances of the estimated parameters for the model~\cite{crestel2017optimal}.
The log D-optimality criterion computes the log of the determinant of $\mathbf{\Gamma}_{\mathrm{post}}$, which is equivalent to the sum of the log of the eigenvalues of $\mathbf{\Gamma}_{\mathrm{post}}$. This metric is proportional to the (log) volume of the confidence ellipsoid of the parameters (lower is better)~\cite{KURAM2013159}.
It can be interpreted in terms of Shannon information.

\subsection{The A14 dataset}
\label{sec:A14}
We chose the A14 tune~\cite{ATL-PHYS-PUB-2014-021} of the \pythia\footnote{To match the original study, we used version {\sc v8.186}.} event generator \cite{Sjostrand:2014zea} as one benchmark for developing and testing the methods proposed in this work.
This tune has been widely used for Large Hadron Collider (LHC) simulations, the methodology was reasonably well documented, and some of the simulation data were available to us.
The generator settings are available in Section~\ref{sec:gen_setup} of the online supplement.

The A14 dataset contains 406 observables (thus, 406 weights to optimize) and there are 10 tunable physics parameters $\p$. %
In our studies, we use the {\sc Rivet} \cite{Buckley:2010ar} package to compare our predictions to data.
The definition and ranges of the parameters $\p$ used to construct the original A14 tune
are shown in Table~\ref{tab:prm-ranges} in Section~\ref{sec:physicsparam} of the online supplement.

Because the coefficients of the cubic interpolation used in \cite{ATL-PHYS-PUB-2014-021} were not available to us, we start
by reproducing the hand-tuned parameter values with the NNPDF published in Table~3 in \cite{ATL-PHYS-PUB-2014-021}, which we refer to as {\it NNPDF}. In particular, we
use the weights given in Table~2 in \cite{ATL-PHYS-PUB-2014-021}, use their optimal parameter values as a starting point for the
$\chi^2$ minimization, and apply our optimizer to \equaref{eq:poly}. The resulting parameter values are reassuringly close
to the values reported in \cite{ATL-PHYS-PUB-2014-021} as showed in Table~\ref{tab:expert_comparison} where we label the original
parameters as NNPDF, and the re-optimized parameter values as Expert. We observe that most of the NNPDF parameter values lie 
within the confidence interval derived from eigentunes (see Section~\ref{sec:eigentune}) for the re-optimized Expert values.
Additionally, to check whether the parameter $\p$ reported in \cite{ATL-PHYS-PUB-2014-021} is within the confidence ellipsoid centered on the parameter $\phatw$ obtained from the $\chi^2$ minimization (i.e., Expert parameter values), we calculate $s \equiv\left\|\mathbf{L}^{T}(\mathbf{\p}-\phatw)\right\|_2$, where $\mathbf{L}$ is the Cholesky factor of $\mathbf{\Gamma}_{\mathrm{post}}(\phatw,\w)$ from \equaref{eq:gammapost} with weights $\w$ given in~\cite{ATL-PHYS-PUB-2014-021}. Since $s = 2.73\times 10^{-3}$  is less than one, we say that the parameter $\p$ is covered within the confidence ellipsoid centered on $\phatw$~\cite{pope2008algorithms}.

In the remainder of this paper, we use the Expert parameter values for comparison, rather than the NNPDF values, and we
refer to this tune as the {\it Expert} tune in our comparisons. This change 
provides a better comparison, because we found that the original NNPDF parameter values did not correspond to a minimizer of 
the $\chi^2$ optimization, \equaref{eq:poly}, and thus using the original values would unfairly disadvantage the NNPDF tune in
our comparisons. The main reason for this discrepancy is the fact that we use a better optimization routine,
and that we did not have access to the coefficients of the cubic interpolation used in \cite{ATL-PHYS-PUB-2014-021}.

\begin{table}[!htbp]
\centering
\caption{Parameter values for A14 published tune (left), and A14 corrected expert tune and corresponding eigentune confidence intervals (right).}
\begin{adjustbox}{width=\textwidth}
\begin{tabular}{l|l|l|l|l}
\hline
& 
\multicolumn{1}{c|}{A14 published expert tune}                                                & \multicolumn{3}{c}{A14 corrected expert tune}   \\\hline
         Parameter name                           & NNPDF & Expert & min    & max     \\\hline
\tt SigmaProcess:alphaSvalue            & 0.140        & 0.143       & 0 & 0.299  \\
\tt BeamRemnants:primordialKThard       & 1.88        & 1.904       & 1.899    & 1.908    \\
\tt SpaceShower:pT0Ref                  & 1.56        & 1.643       & 1.621   & 1.667    \\
\tt SpaceShower:pTmaxFudge              & 0.91        & 0.908       & 0.897   & 0.917    \\
\tt SpaceShower:pTdampFudge             & 1.05        & 1.046       & 1.041 & 1.052 \\
\tt SpaceShower:alphaSvalue             & 0.127       & 0.123       & 0.120 & 0.127  \\
\tt TimeShower:alphaSvalue              & 0.127       & 0.128       & 0 & 0.355  \\
\tt MultipartonInteractions:pT0Ref      & 2.09        & 2.149       & 1.055   & 3.442    \\
\tt MultipartonInteractions:alphaSvalue & 0.126       & 0.128       & 0 & 0.289  \\
\tt BeamRemnants:reconnectRange         & 1.71        & 1.792       & 1.784   & 1.802  
\\ \hline
\end{tabular}
\end{adjustbox}
\label{tab:expert_comparison}
\end{table}

The 10 tunable parameters are primarily related to the production of additional jets in the collisions, the distribution of energy within those jets,
and the kinematics (angles and momenta) of the jets.   They also relate to the sharing and spread of energy in the soft portion of the event,
the portion that is less dependent on the hard process (e.g., top-quark production or $Z$-boson production). 

The A14 observables are measurements of properties of proton-proton collisions at $\sqrt{s}=7$~TeV performed by the ATLAS collaboration.
These include event properties (e.g., the $Z$-boson transverse momentum, or the opening angles between the highest transverse momentum jets in the event) and properties of jets (e.g., the spread of energy within a jet, or the momentum of particles within a jet). In the publication~\cite{ATL-PHYS-PUB-2014-021}, the 406 observables are categorized into 10 groups (see Table~\ref{tab:weights}), namely {\it Track jet} properties (200 observables), {\it Jet shapes} (59 observables), {\it Dijet decorr} (9 observables),
{\it Multijets} (8 observables), $p_T^Z$ (fit range $<50$GeV, 20 observables), {\it Substructure} (36 observables), {\it $t\bar{t}$ gap} (4 observables),
{\it $t\bar{t}$ jet shapes} (20 observables), {\it Track-jet UE} (8 observables), and {\it Jet UE} (42 observables). 
The highest weights in \cite{ATL-PHYS-PUB-2014-021} are assigned to observables that relate to the production of additional high-momentum partons (the ratios of 3-jet to 2-jet events, and the fraction of top-quark production events that do not have an additional central jet). On the other hand, low weights are assigned to observables that measure the same physical phenomenon in several kinematic regimes.
The weighting of these observables ensures that the additional radiation and soft part of the events are consistent and well-modeled for all hard processes.
In addition, these parameters are difficult or impossible to constrain using data from $e^+e^-$ collision events, and they must be tuned using data from the LHC.

\subsection{The \sherpa dataset}
As a second benchmark, we tune a set of parameters for the \sherpa event generator \cite{Bothmann:2019yzt}.   To our knowledge, the default parameters were not optimized by weighting
data, and thus serve as an unbiased cross-check of our results.
In contrast to the A14 dataset used to tune \pythia, the data are confined to observables at $e^+e^-$ colliders.
The data includes event shapes and charged particle inclusive spectra from $Z$-boson decays,
differential and integrated jet rates, measurements of $B$-hadron fragmentation,
and the multiplicity of various hadrons \cite{Pfeifenschneider:1999rz,Abreu:1996na,Abe:2002iq,Amsler:2008zzb}.
Accordingly, the parameters are limited to those of the \sherpa hadronization model.

The \sherpa dataset contains 88 observables, hence 88 weights to optimize.
The number 88 is significantly less than the set of observables available in the {\sc Rivet} analyses (126) for the following reasons.
First, we reduce the number of observables to 114 by removing those that measure more than 3 jets,
since this is beyond the scope of the physics simulation.
Then, we apply a pre-filter step that removes distributions where {\it none} of the data bins fall within the envelope of
predictions from our surrogate model.   These all correspond to single-bin particle counts (such as
the number of $f_0$ mesons) that the \sherpa hadronization model either grossly under- or over-estimates. 
There are 13 tunable physics parameters whose definition and ranges are shown in
Table~\ref{tab:prm-ranges-sherpa} in Section~\ref{sec:physicsparam} of the online supplement.
These parameters are all part of the cluster model that produces physical particles from quarks and gluons.

\subsection{Data pre-processing: filtering out observables and bins}
In this subsection, we present the results of applying the filtering methods described in Sections~\ref{sec:obsfilter} and~\ref{sec:binfilter}.
First, we consider the outlier detection method described in Section~\ref{sec:obsfilter}.
We find that the filtering results differ based on the choice of surrogate function (cubic polynomial versus a rational approximation).
Based on the comparison of surrogate function predictions to the full MC simulations, we believe that the rational approximation yields
a more faithful representation.  Therefore, we present our main results using only the rational approximation.
The names of the outlier observables in the A14 and the \sherpa dataset using a cubic polynomial and a rational approximation, respectively, are shown in the online supplement in Sections~\ref{sec:A14outliers} and~\ref{sec:sherpaoutliers}.
Table~\ref{tab:outlier} shows a distribution of the $\chi^2_\calO$ values obtained for each observable $\calO$ from~\equaref{eq:obschi2} for A14  (left) and \sherpa (right) when using the rational approximation.
We find that the per-observable ideal parameters yield mostly small $\chi^2_\calO$ values (in $[0,1)$), but outliers  are present in both datasets. %
Using the rational approximation, 9 and 3 outlier observables are filtered from the A14 and \sherpa datasets, respectively.

\begin{table}[!htbp]
\centering
\caption{Distribution of the $\chi^2_\calO$ values for A14 (left) and \sherpa (right). $2.0438$ and $2.0177$ correspond to the values where the Z-score equals 3 (see Section~\ref{sec:obsfilter}). The observables with $\chi^2_\calO \geq 2.0438$ for A14 and $\chi^2_\calO \geq 2.0177$ for \sherpa are the outliers. There are 9 outliers (6+2+1) in A14 and 3 outliers (1+2+0) in \sherpa.}
\begin{tabular}{l|l|l|l}
\hline
\multicolumn{2}{c|}{A14}                                                & \multicolumn{2}{c}{\sherpa}                                             \\ \hline
\multicolumn{1}{c|}{$\chi^2_\calO$ range} & \multicolumn{1}{c|}{Number of observables} & \multicolumn{1}{c|}{$\chi^2_\calO$ range} & \multicolumn{1}{c}{Number of observables} \\ \hline
{[}0, 1)                     & 367                                        & {[}0, 1)                    & 82                                         \\
{[}1, 2.0438)    & 30                                         & {[}1, 2.0177)   & 3                                          \\
{[}2.0438, 3)    & \textbf{6}                                 & {[}2.0177, 3)   & \textbf{1}                                 \\
{[}3, 4)                     & \textbf{2}                                 & {[}3, 4)                    & \textbf{2}                                 \\
{[}4, 5)                     & \textbf{1}                                 & {[}4, 5)                    & \textbf{0}                                 \\ \hline
\end{tabular}
\label{tab:outlier}
\end{table}

\begin{figure}[!htbp]
	\begin{subfigure}[b]{\textwidth}
         \centering
         \includegraphics[width=\textwidth]{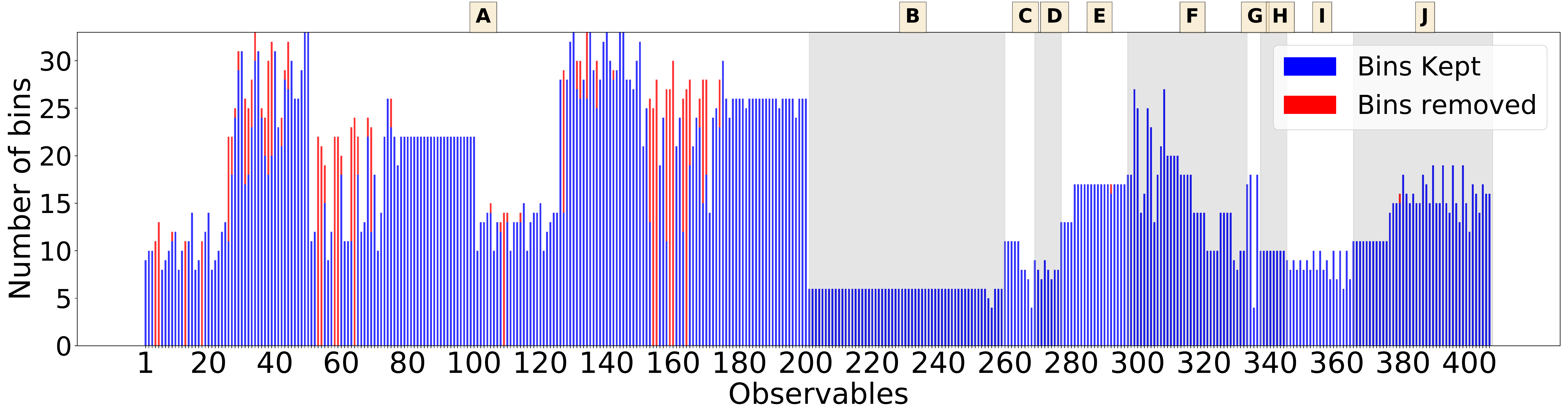}
         \caption{Bins kept and removed by the bin filter in all A14 observables organized by the observable group.    
         Group A is \textit{Track jet properties},
        group B is \textit{Jet shapes},
        group C is \textit{Dijet decorr},
        group D is \textit{Multijets},
        group E is $p^Z_T$,
        group F is \textit{Substructure},
        group G is \textit{$t \Bar{t}$  gap},
        group H is \textit{Track-jet UE},
        group I is \textit{$t \Bar{t}$  jet shapes}, and
        group J is \textit{Jet UE}.
        }
         \label{fig:binf-a14}
     \end{subfigure}
     \hfill
     \begin{subfigure}[b]{\textwidth}
         \centering
         \includegraphics[width=\textwidth]{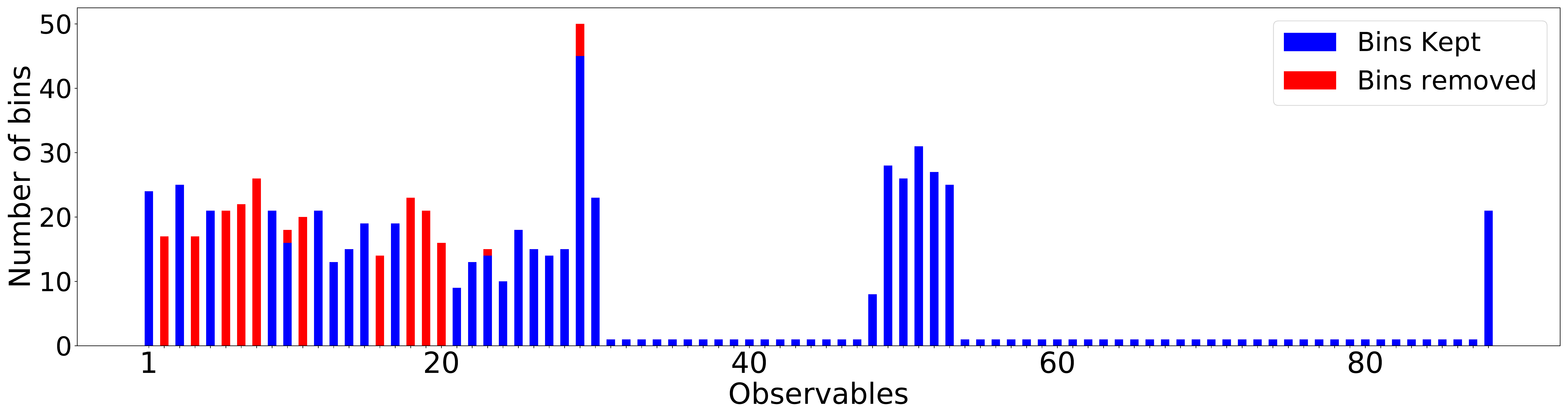}
         \caption{Bins kept and removed by the bin filter in all \sherpa observables.}
         \label{fig:binf-sherpa}
     \end{subfigure}
 	\caption{Illustration of the   bin filtering results.}
	\label{fig:binf}
\end{figure}

Figure~\ref{fig:binf} shows the outcomes of the bin-filtering approach described in Section~\ref{sec:binfilter} for each observable $\calO$ %
in A14 (top) and \sherpa (bottom) when using the rational approximation. In both datasets, multiple bins are removed. More specifically, most bins are removed in the \textit{Track jet properties} and $p^Z_T$ groups of the A14 dataset.   The patterns in the A14 plot result from the partitioning of the data.
For {\it Tracked jet properties} (labeled A), the observables are replicated for two values of jet cone size ($R = 0.4, 0.6$), explaining the similarities between
bins (1, 100) and (101, 200).   Furthermore, 4 types of observables are considered, and each is sliced into different ranges of transverse momentum and rapidity.

In the \sherpa dataset, all bins are removed from some observables whereas from two observables, we remove only two and five bins.  Additionally, since the number of degrees of freedom of the $\chi^2$ distribution is reduced by the number of parameters that the bins share in each observable (see Eq.~\eqref{eq:degoffreedom}), the bin filter is not applied to any observable with fewer than 10 and 13 bins in the A14 and the \sherpa datasets, respectively. The names of the observables from which the bins have been filtered and their $\chisq$ test statistic and critical $\chisq$ values are given in Sections~\ref{sec:binfiltA14} and~\ref{sec:binfiltSherpa} of the online supplement.   The single bin observables correspond to counts of a particular type of particle.

\subsection{Results for the A14 dataset}
In this section, we present a detailed analysis of our results with the A14 dataset.
\subsubsection{Comparison metric outcomes for the A14 dataset}
In this section, we consider the three metrics introduced in Section~\ref{sec:metrics} to compare various tunes.
For the A14 dataset, Tables~\ref{tab:comparison_full_31}-\ref{tab:comparison_binf_31} show the  results when using the rational approximation for the full data, the observable-filtered data, and the bin-filtered data, respectively. The results when using the cubic polynomial approximation are shown in the online supplement in Section~\ref{sec:a14_30_metric}. 
 Note that smaller numbers indicate better performance. We bold the smallest number of each metric for better visualization. For our comparison metrics, we take into account all observables and bins, respectively, but we do not use the filtered out observables and bins when determining the optimal parameters.

 Based on these results we can see that no method performs the best for all metrics in all cases.
 In fact, for the full dataset, the {\it Expert} tune has the best score for two of our three metrics.
 Nonetheless, the automated methods do produce comparable results in those cases.
The robust optimization consistently achieves the best performance under the Weighted $\chi^2$ criterion. %

The Bilevel-portfolio method performs the best under the A-optimality criteria, and the {\it Expert} tune performs the best under the D-optimality criteria for the observable-filtered datasets.
The Bilevel-portfolio method performs the best under the A- and D-optimality criteria for the bin-filtered datasets.
In  comparison to the results obtained with  the cubic polynomial approximation (see Section~\ref{sec:a14_30_metric} of the online supplement),   the rational approximation yields better results  for all methods under the Weighted $\chi^2$ criterion.

When comparing across Tables~\ref{tab:comparison_full_31}-\ref{tab:comparison_binf_31}, we see that in most cases, results with the observable-filtered data and bin-filtered data provide smaller values compared with those using the full dataset. We observe that by filtering out the observables and bins that cannot be well explained by the model, the quality of the model fits  can be improved.  

\begin{table}[!htbp]
	\caption{A14 results with the  \textit{full dataset} when using the \textit{rational} approximation. Lower numbers are better. The best results are in bold. %
	}
	\centering
	\begin{tabular}{l|rrr}
		\hline
		Method               & Weighted $\chi^2$ & A-optimality    & D-optimality (log) \\\hline
		Bilevel-meanscore & 0.1119	& 0.8513 & 	-63.6805
\\
		Bilevel-medscore  & 0.1320 & 	0.7673	 & -63.3846
 \\
		Bilevel-portfolio & 0.1224	 & 0.9425 & 	-61.1694
\\
		Expert        & 0.0965	 & \textbf{0.5705} & 	\textbf{-68.4091}
 \\
		All-weights-equal & 0.0815 & 	0.7673	 & -64.0008

 \\
		Robust optimization & \textbf{0.0402}	 & 1.0526 & 	-65.7547
	\\\hline
	\end{tabular}
\label{tab:comparison_full_31}
\end{table}

\begin{table}[!htbp]
	\caption{A14 results with the \textit{observable-filtered} data when using the \textit{rational} approximation. Lower numbers are better. The best results are indicated in bold. %
	}
	\centering
	\begin{tabular}{l|rrr}
		\hline
		Method         & Weighted $\chi^2$ & A-optimality & D-optimality (log) \\\hline
		Bilevel-meanscore & 0.0671 &	0.6793	 &	-65.1939
        \\
		Bilevel-medscore  &0.0823 &		0.7008	 &	-64.3410
 \\
		Bilevel-portfolio & 0.1372 &		\textbf{0.5130}	 &	-68.0382
\\
		Expert        & 0.0965	 &	0.5705	 &	\textbf{-68.4091}
 \\
		All-weights-equal      &   0.0843	 &	0.7110 &		-64.4546

\\
		Robust optimization            &         \textbf{0.0388}	 &	1.1086 &		-65.7182

		\\\hline
	\end{tabular}
	\label{tab:comparison_obsf_31}
\end{table}

\begin{table}[!htbp]
	\caption{A14 results with the \textit{bin-filtered} data when using the \textit{rational} approximation. Lower numbers are better. The best results are in bold. %
	1811 out of 7010 bins were filtered out and not used during the optimization. All numbers are computed over all 7010 bins.}
	\centering
	\begin{tabular}{l|rrr}
		\hline
		method         & Weighted $\chi^2$ & A-optimality & D-optimality (log) \\\hline
		Bilevel-meanscore   & 0.0677 &	0.7098	 &	-66.1484

 \\
		Bilevel-medianscore &0.1198	 &	0.7448	 &	-64.2754

\\
		Bilevel-portfolio   & 0.1464 &		\textbf{0.3747} &		\textbf{-70.5889}

        \\
		Expert        &0.0965 &		0.5705 &		-68.4091
\\
		All-weights-equal  & 0.0835	 &	0.7185	 &	-65.3507

  \\
		Robust optimization   &\textbf{0.0439} &		0.8332	 &	-67.4820
\\\hline
	\end{tabular}
	\label{tab:comparison_binf_31}
\end{table}

\subsubsection{Comparison of the cumulative distribution of bins at different variance levels}
\label{sec:cdfperfA14}

In this section, we introduce a new summarized graphical comparison of the results that is motivated by the bottom pane in the histogram plot of Figure~\ref{fig:hist0}.
We study the distribution of the $\chi^2$ values per bin obtained using
different tuning approaches.
For each   parameter set, we compute the
ratio  $\displaystyle r_b(\p)=\frac{\left(f_b(\p) -
    \calR_b\right)^2}{\Delta f_b(\p)^2+\Delta\calR_b^2}$
of the residual between the data and the prediction divided by the
variance per bin.   The $r_b$ values are sorted from the smallest to
largest, and the cumulative distribution is formed.

\begin{figure}[ht!]
  \centering
  \includegraphics[width=\textwidth]{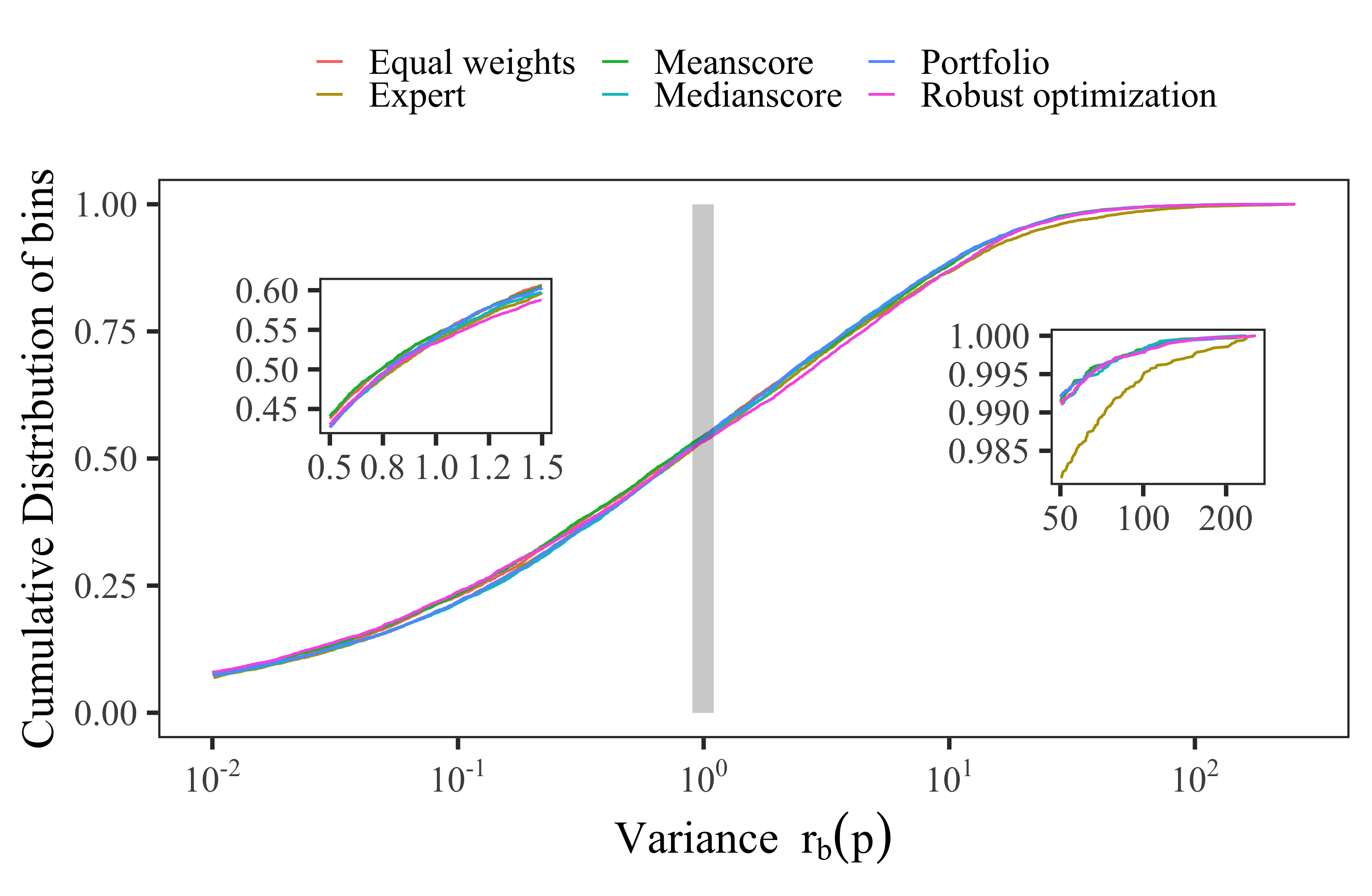}  
  \caption{
    Cumulative distribution of all bins (y-axis) in the A14 dataset at different bands of variance levels (x-axis) given by $r_b(\p)=\frac{\left(f_b(\p) - \calR_b\right)^2}{\Delta f_b(\p)^2+\Delta\calR_b^2}$.
  }
  \label{fig:chi2perf2binsAllA14}
\end{figure}

\begin{figure}[ht!]
  \centering
  \includegraphics[width=\textwidth]{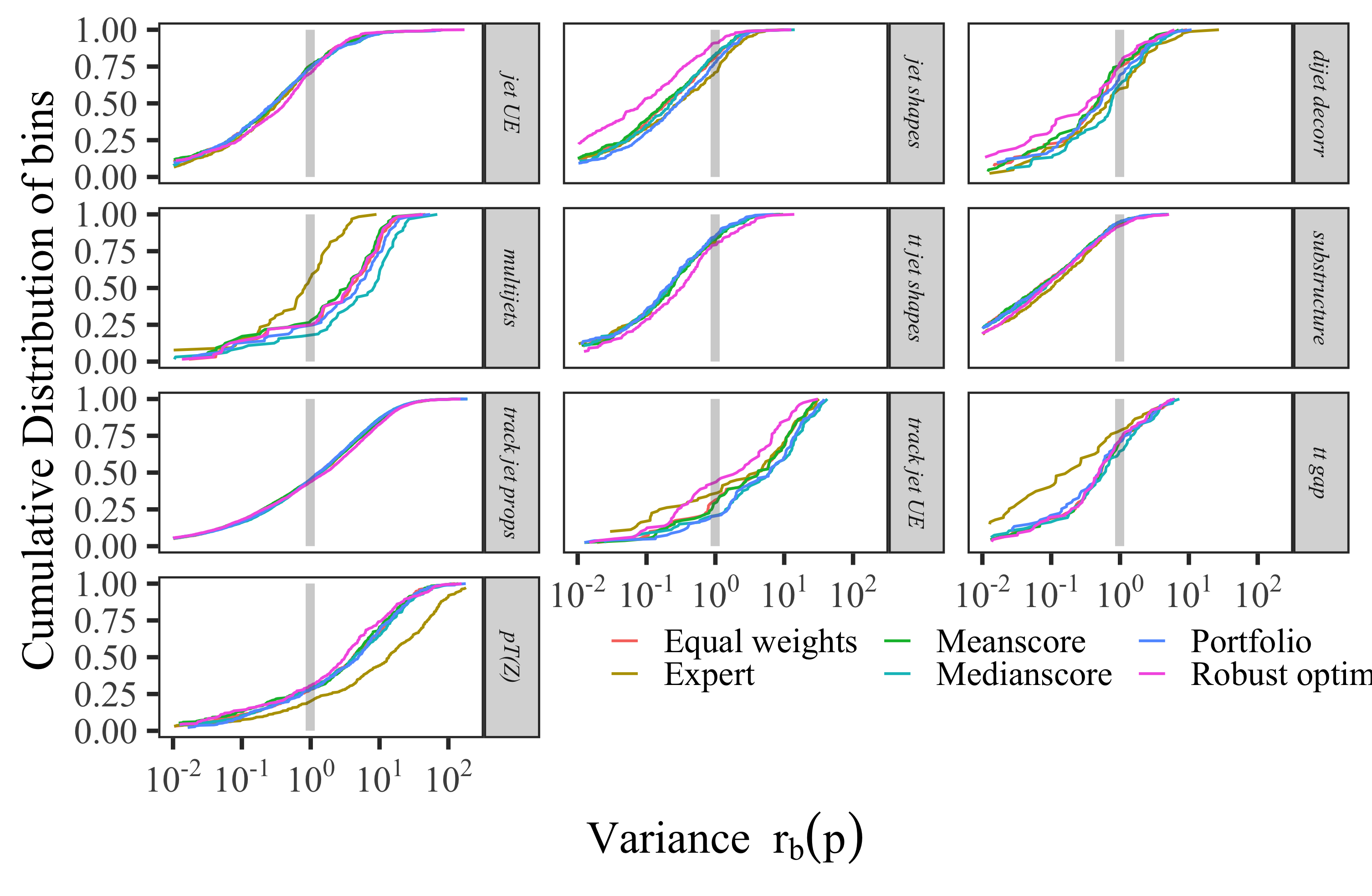}  
  \caption{
    Cumulative distribution of bins (y-axis) in each category of the A14 dataset at different bands of variance levels (x-axis) given by $r_b(\p)=\frac{\left(f_b(\p) - \calR_b\right)^2}{\Delta f_b(\p)^2+\Delta\calR_b^2}$.
  }
  \label{Fig:chi2perf2binsCat}  
\end{figure}

The cumulative distribution plot  for all bins in the A14 dataset is shown in Figure~\ref{fig:chi2perf2binsAllA14} and for the bins in each category  in Figure~\ref{Fig:chi2perf2binsCat}.
The more bins that reside on the bands of variance levels less than 1 the better, as this indicates smaller deviations of the model from the experimental data. 
When analyzing these results it is important to note that even though all the category plots have a scale between 0 and 1 on the y axis, the number of bins in one category of A14 is very different from the other. For e.g., more than 50\% of all bins in the A14 dataset belong to \textit{Track Jet Properties}. Hence, we see that the trend of the curves in the plot for \textit{Track Jet Properties} in Figure~\ref{Fig:chi2perf2binsCat} follows more closely to the trend of the curves when all A14 bins are considered as in Figure~\ref{fig:chi2perf2binsAllA14}.

It can be seen from Figure~\ref{fig:chi2perf2binsAllA14} that there is
a small difference among the approaches when all A14 bins are
considered.
Near the variance boundary, the difference between the
approaches is even smaller.  Additionally, at the variance boundary, all
approaches perform better than the \textit{Expert} tune. 
 Figure~\ref{Fig:chi2perf2binsCat} shows that these differences become more prominent when considering individual  categories of the A14 data. For instance, the parameters obtained from the robust optimization perform well for \textit{Jet shapes} and \textit{Track-jet UE}. 
We also see that near the variance boundary, the parameters obtained from the \textit{Expert} tune perform better for \textit{Multijets} and \textit{$t \Bar{t}$ gap} whereas the parameters obtained from the other approaches  perform better for \textit{Substructure}. 
These plots also show that there is a trade-off in fitting among the different approaches, which enables the physicist to use these results as guidance for selecting the most appropriate tuning method  depending on the categories that are of greater significance.

\subsubsection{Optimal parameter values for the A14 dataset with rational approximation}

The optimal parameter values for the A14 dataset when using the full dataset, the outlier-filtered dataset, and the bin-filtered dataset are shown 
in Tables~\ref{tab:optimal_prms_31}, \ref{tab:optimal_prms_obsfilter_31}, and  \ref{tab:optimal_prms_binfilter_31}, respectively. For a better visual comparison of the different solutions obtained with our methods, we illustrate the [0,1]-scaled optimal values in the online supplement Section~\ref{sec:pstar_a14_31}.
We have also computed the Euclidean distance between the \textit{Expert} tune and our tunes after normalizing the parameter values to [0,1].%

In Table~\ref{tab:optimal_prms_31}, we can see that there are differences between the optimal parameters obtained with different methods. In particular, the results of the Bilevel-meanscore method tend to be further away from the expert's solution than the other methods. The robust optimization and All-weights-equal results are very similar to each other as well as to the \textit{Expert}'s solution.

\begin{table}[!htbp]
	\centering
	\small
	\caption{Optimal parameter values for the A14 dataset obtained when using all observables in the optimization and the rational approximation. Euclidean distance is calculated based on the normalized parameter values.}\label{tab:optimal_prms_31}
	\begin{adjustbox}{width=\textwidth}
	\begin{tabular}{llllllll}
		\hline
		ID & Parameter name & Expert & Bil.-meanscore & Bil.-medianscore & Bil.-portfolio & Robust opt    & All-weights-equal \\\hline
1  & \tt SigmaProcess:alphaSvalue            & 0.143  & 0.138 & 0.133 & 0.136 & 0.139 & 0.137 \\
2  & \tt BeamRemnants:primordialKThard       & 1.904  & 1.855 & 1.723 & 1.796 & 1.883 & 1.851 \\
3  & \tt SpaceShower:pT0Ref                  & 1.643  & 1.532 & 1.184 & 1.322 & 1.588 & 1.493 \\
4  & \tt SpaceShower:pTmaxFudge              & 0.908  & 1.014 & 1.083 & 1.041 & 1.025 & 1.026 \\
5  & \tt SpaceShower:pTdampFudge             & 1.046  & 1.071 & 1.084 & 1.061 & 1.084 & 1.067 \\
6  & \tt SpaceShower:alphaSvalue             & 0.123 & 0.128 & 0.129 & 0.128 & 0.127 & 0.128 \\
7  & \tt TimeShower:alphaSvalue              & 0.128 & 0.130 & 0.129 & 0.128 & 0.132 & 0.129 \\
8  & \tt MultipartonInteractions:pT0Ref      & 2.149  & 2.033 & 1.883 & 1.937 & 2.052 & 2.076 \\
9  & \tt MultipartonInteractions:alphaSvalue & 0.128 & 0.124 & 0.118 & 0.120 & 0.126 & 0.125 \\
10 & \tt BeamRemnants:reconnectRange         & 1.792  & 2.082 & 1.914 & 1.987 & 2.602 & 1.980
 \\\hline\hline
& Euclidean distance from the expert solution & &	0.290 & 	0.664 & 	0.475 & 	0.268 & 	0.301

 \\\hline %
	\end{tabular}
	\end{adjustbox}
\end{table}

\begin{table}[htbp]
	\small
	\centering
	\caption{Optimal parameter values for A14 when using the rational approximation with all methods after outlier detection to  filter out observables that cannot be approximated well by the model. Euclidean distance is calculated based on the normalized parameter values. }\label{tab:optimal_prms_obsfilter_31}
	\begin{adjustbox}{width=\textwidth}
	\begin{tabular}{llllllll}
		\hline
	ID	& Parameter name & Expert & Bilevel-meanscore & Bilevel-medianscore & Bilevel-portfolio & Robust opt    & All-weights-equal \\\hline
1  & \tt SigmaProcess:alphaSvalue            & 0.143  & 0.140 & 0.138 & 0.141 & 0.138 & 0.139 \\
2  & \tt BeamRemnants:primordialKThard       & 1.904  & 1.865 & 1.839 & 1.861 & 1.879 & 1.843 \\
3  & \tt SpaceShower:pT0Ref                  & 1.643  & 1.574 & 1.603 & 1.593 & 1.614 & 1.550 \\
4  & \tt SpaceShower:pTmaxFudge              & 0.908  & 0.953 & 0.906 & 0.984 & 1.006 & 0.950 \\
5  & \tt SpaceShower:pTdampFudge             & 1.046  & 1.076 & 1.081 & 1.060 & 1.075 & 1.062 \\
6  & \tt SpaceShower:alphaSvalue             & 0.123 & 0.128 & 0.128 & 0.129 & 0.128 & 0.127 \\
7  & \tt TimeShower:alphaSvalue              & 0.128 & 0.123 & 0.123 & 0.118 & 0.132 & 0.124 \\
8  & \tt MultipartonInteractions:pT0Ref      & 2.149  & 2.064 & 2.017 & 2.095 & 2.022 & 2.039 \\
9  & \tt MultipartonInteractions:alphaSvalue & 0.128 & 0.126 & 0.125 & 0.129 & 0.125 & 0.126 \\
10 & \tt BeamRemnants:reconnectRange         & 1.792  & 1.852 & 1.903 & 1.801 & 2.719 & 1.937
\\\hline
		\hline
& Euclidean distance from the expert solution & & 0.227 & 	0.293 & 	0.273 & 	0.291 & 	0.254

 \\\hline %
	\end{tabular}
	\end{adjustbox}
\end{table}

\begin{table}[htbp]
	\small
	\centering
	\caption{Optimal parameter values obtained for A14 with the rational approximation with all methods after using the bin-filtering approach that excludes individual bins from the optimization 1811 bins out of 7010 total bins were filtered out. Euclidean distance is calculated based on the normalized parameter values.}\label{tab:optimal_prms_binfilter_31}
	\begin{adjustbox}{width=\textwidth}
	\begin{tabular}{llllllll}
		\hline
	ID	& Parameter name & Expert & Bilevel-meanscore & Bilevel-medianscore & Bilevel-portfolio & Robust opt    & All-weights-equal \\\hline
1  & \tt SigmaProcess:alphaSvalue            & 0.143  & 0.139 & 0.140 & 0.131 & 0.137 & 0.140 \\
2  & \tt BeamRemnants:primordialKThard       & 1.904  & 1.877 & 1.885 & 1.811 & 1.822 & 1.876 \\
3  & \tt SpaceShower:pT0Ref                  & 1.643  & 1.572 & 1.561 & 2.227 & 1.426 & 1.627 \\
4  & \tt SpaceShower:pTmaxFudge              & 0.908  & 0.964 & 0.968 & 0.869 & 0.948 & 0.943 \\
5  & \tt SpaceShower:pTdampFudge             & 1.046  & 1.056 & 1.053 & 1.481 & 1.053 & 1.068 \\
6  & \tt SpaceShower:alphaSvalue             & 0.123
 & 0.128 & 0.128 & 0.136 & 0.128 & 0.128 \\
7  & \tt TimeShower:alphaSvalue              & 0.128 & 0.128 & 0.129 & 0.126 & 0.136 & 0.130 \\
8  & \tt MultipartonInteractions:pT0Ref      & 2.149  & 2.028 & 2.175 & 2.338 & 1.931 & 2.080 \\
9  & \tt MultipartonInteractions:alphaSvalue & 0.128 & 0.124 & 0.128 & 0.135 & 0.120 & 0.126 \\
10 & \tt BeamRemnants:reconnectRange         & 1.792  & 2.047 & 1.854 & 1.820 & 2.404 & 2.001
 \\\hline 		\hline
& Euclidean distance from the expert solution & & 0.232 & 	0.179	 & 1.076 & 	0.426 & 	0.194

 \\\hline %
	\end{tabular}
	\end{adjustbox}
\end{table}

\subsubsection{Comparison of optimal weights for the A14 dataset with rational approximation}
We compare the optimal weights obtained by the different tuning methods in Table~\ref{tab:weights}. %
We normalize the weights obtained to match the scale of weights assigned by \textit{Expert} published in~\cite{ATL-PHYS-PUB-2014-021}. In each group, we report the average weight of observables in that group.
The {\it Expert} tune assigned the highest weights to the categories \textit{Multijets} and \textit{$t\bar{t}$ gap}. 
The robust optimization approach sets some of the weights for \textit{Track jet properties} to zero. The four \textit{Track-jet properties} classes of observables are nearly dependent  resulting
in redundant components of least-square residuals. Because the robust optimization
approach can be viewed as minimizing the maximum residual, it detects this redundancy of observables, and sets the weights accordingly to zero. We observe in
Figure~\ref{Fig:chi2perf2binsCat} that setting these weights to zero does not degrade the residuals of these observables, confirming that redundant information is presented.

  \begin{table}[!htbp]
		\centering
		\small
		\caption{Comparison of the optimal weights obtained by each method using the rational approximation. The observable grouping corresponds to the same grouping as in~\cite{ATL-PHYS-PUB-2014-021}.}\label{tab:weights}
		\begin{adjustbox}{width=\textwidth}
			\begin{tabular}{p{9cm}p{1.5cm}p{1.5cm}p{1.5cm}p{1.5cm}p{1.5cm}}
				\hline
				& expert & Bilevel-meanscore & Bilevel-medianscore  & Bilevel-portfolio & robustopt \\\hline
			\textbf{Track jet properties}                                        &     &       &       &       &       \\
			Charged jet multiplicity (50   distributions)                                           & 10  & 11.41 &	11.92 &	11.43 &	17.85
 \\
			Charged jet $z$ (50 distributions)                                                      & 10  & 11.01 &	10.00 &	10.28 &	0.00
 \\
			Charged jet $p^{rel}_T$ (50   distributions)                                            & 10  & 9.47 &	10.20 &	13.11 &	1.62
\\
			Charged jet $\rho_{ch}(r)$ (50   distributions)                                         & 10  & 10.63 &	12.72 &	12.19 &	0.00
 \\\hline
			\textbf{Jet shapes}                                                    &     &       &       &       &       \\
			Jet shape $\rho$ (59 distributions)                                                     & 10  & 12.46 &	8.49 &	9.69 &	17.85
 \\\hline
			\textbf{Dijet decorr}                                                  &     &       &       &       &       \\
			Decorrelation $\Delta \phi$ (Fit range:   $\Delta \phi>0.75$) (9 distributions)         & 20  &  18.82 &	10.32 &	18.50 &	15.87
\\\hline
			\textbf{Multijets}                                                     &     &       &       &       &       \\
			3-to-2 jet ratios (8 distributions)                                                     & 100 &  15.06 &	11.18 &	11.06 &	17.85
\\\hline
			\textbf{$p^Z_T$} (Fit range: $p^Z_T<50   \text{GeV}$)                  &     &       &       &       &       \\
			Z-boson $p_T$ (20 distributions)                                                        & 10  & 12.16 &	11.85 &	9.25 &	17.85
\\\hline
			\textbf{Substructure}                                                  &     &       &       &       &       \\
			Jet mass, $\sqrt{d_{12}},  \sqrt{d_{23}},  \tau_{21}, \tau_{23}$ (36 distributions)     & 5   & 10.71 &	12.75 &	14.23 &	17.85
\\\hline
			\textbf{$t \Bar{t}$ gap}                                               &     &       &       &       &       \\
			Gap fraction vs $Q_0$,  $Q_{\text{sum}}$ for $|y|<0.8$                                  & 100 & 24.56 &	5.05 &	1.97 &	17.85
 \\
			Gap fraction vs $Q_0$,  $Q_{\text{sum}}$ for $0.8<|y|<1.5$                              & 80  & 23.73 &	47.01 &	4.01 &	17.85
 \\
			Gap fraction vs $Q_0$,  $Q_{\text{sum}}$ for $1.5<|y|<2.1$                              & 40  & 2.39 &	14.20 &	7.35 &	17.85
 \\
			Gap fraction vs $Q_0$,  $Q_{\text{sum}}$ for $|y|<2.1$                                  & 10  & 5.47 &	19.00 &	12.82 &	17.85
 \\\hline
			\textbf{Track-jet UE}                                                  &     &       &       &       &       \\
			Transverse region $N_{ch}$ profiles (5   distributions)                                 & 10  &  13.01 &	24.18 &	7.46 &	17.85
\\
			Transverse region mean $p_T$ profiles for   $R=0.2, 0.4, 0.6$ (3 distributions)         & 10  &  7.91 &	16.89 &	9.68 &	17.85
\\\hline
			\textbf{$t \Bar{t}$ jet shapes}                                       &     &       &       &       &       \\
			Jet shapes $\rho(r),  \psi(r)$ (20 distributions)                                       & 5   & 10.44 &	11.47 &	10.29 &	15.17
 \\\hline
			\textbf{Jet UE}                                                        &     &       &       &       &       \\
			Transverse, trans-max, trans-min sum   $p_T$ incl. profiles (3 distributions)           & 20  &  12.11 &	5.32 &	10.51 &	17.85
\\
			Transverse, trans-max, trans-min $N_{ch}$   incl. profiles (3 distributions)            & 20  & 6.16 &	14.42 &	6.56 &	17.85
\\
			Transverse sum $E_T$ incl. profiles (2   distributions)                                 & 20  & 5.11 &	2.71 &	7.72 &	17.85
  \\
			Transverse sum $ET/$sum $p_T$ ratio incl.,   excl. profiles (2 distributions)            & 5   & 11.94 &	10.81 &	11.65 &	17.85
 \\
			Transverse mean $p_T$ incl. profiles (2   distributions)                                & 10  & 12.47 &	7.28 &	10.45 &	17.85
 \\
			Transverse, trans-max, trans-min sum   $p_T$ incl. distributions (15 distributions)     & 1   & 10.54 &	14.44 &	8.27 &	17.85
 \\
			Transverse,  trans-max, trans-min sum $N_{ch}$ incl.   distributions (15 distributions) & 1   & 11.62 &	10.33 &	11.48 &	17.85
\\\hline
			\end{tabular}
                      \end{adjustbox}
                      \label{tab:a14_weights}
	\end{table}

\subsubsection{Impact of data pre-processing by filtering on optimal results}
\label{sec:filterresults}
\begin{figure}[ht!]
  \centering
    \includegraphics[width=\textwidth]{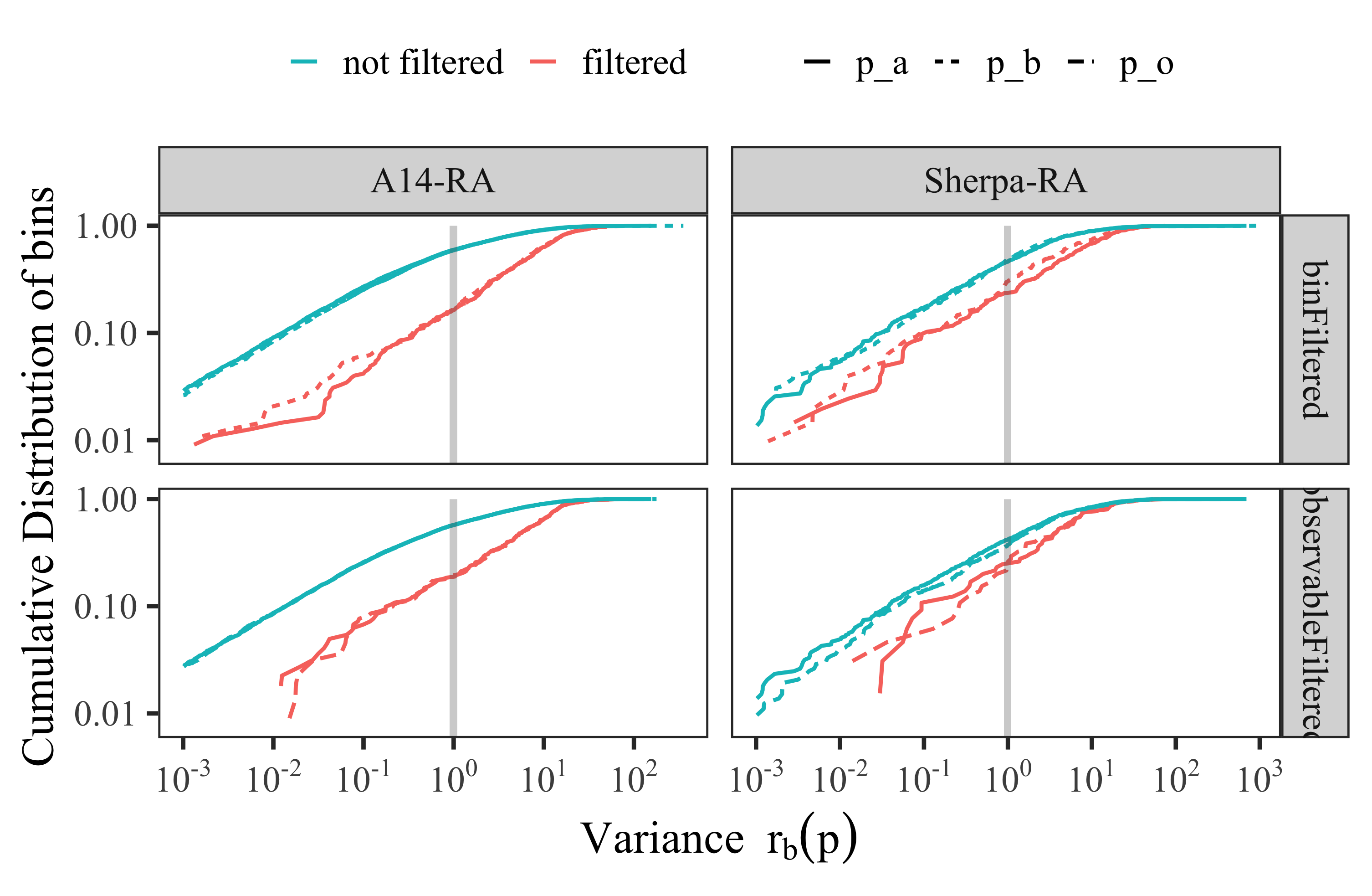}
	\caption{
		Cumulative distribution of bins remaining after filtering (\textit{not filtered}) and of those filtered out (\textit{filtered}) on the y-axis at different bands of variance levels on the x-axis. The variance level for each bin is calculated as $r_b(\p)=\frac{\left(f_b(\p) - \calR_b\right)^2}{\Delta f_b(\p)^2+\Delta\calR_b^2}$  with parameters $\p^*_a$, which is obtained when all bins were used, and parameters $\p^*_b$ and $\p^*_o$, which are obtained when only the bin filtered and observable filtered data are used, respectively. 
	}
	\label{fig:chi2perf2binFilter}
\end{figure}

In Table~\ref{tab:nbinatonevariance}, we show the number of filtered and unfiltered bins in the A14 and \sherpa datasets that lie  within a one $\sigma$ variance level. A large number of bins within a one $\sigma$ level indicates smaller deviations of the model from the experimental data. 
The cumulative distribution plot with the parameters obtained from the robust optimization approach for filtered and unfiltered data for the different categories is shown in Figure~\ref{fig:chi2perf2binFilter} (the  plots for the other methods are  shown in Section~\ref{sec:filterresonlinesup} of the online supplement).

From these results, we observe that there is no significant difference in the number of bins within the one $\sigma$ variance level between the optimal parameters $\p^*_a$ obtained when all bins were used for tuning and the optimal parameters $\p^*_b$ and $\p^*_o$ obtained when only the bin filtered and observable filtered bins are used for tuning, respectively.
Additionally, when comparing across Tables~\ref{tab:comparison_full_31}-\ref{tab:comparison_binf_31}, we see that in most cases, the results with the observable-filtered data and bin-filtered data provide smaller values in the proposed criteria compared with those using the full dataset. 
These observations indicate that the MC generator cannot explain the bins removed by the filtering approaches well and that the information contained in these bins does not add significant information to the tune.

\begin{table}[!htbp]
	\small
	\caption{
		\small Number of bins in the A14 and \sherpa datasets within the one $\sigma$ variance level. Larger numbers are better. The variance level for each bin is calculated as $r_b(\p)=\frac{\left(f_b(\p) - \calR_b\right)^2}{\Delta f_b(\p)^2+\Delta\calR_b^2}$.
		\textit{Test data type} specifies the data over which $r_b(\p)$ is calculated, where \textit{All} means that all bins are used, \textit{Not Filtered} refers to only the  bins that remain after filtering, and \textit{Filtered} refers to the bins that were filtered out by the respective filter specified in the \textit{Filtering Method} as well as the envelope filter. For each data type, the number of bins in the corresponding dataset is also specified.
		\textit{Parameters} specify the type of optimal parameters used in $r_b(\p)$ where $\p^*_a$ are the parameters obtained when all bins were used during tuning  whereas $\p^*_b$ and $\p^*_o$ are the parameters obtained when only the bin filtered and observable filtered date are used, respectively.
	 } \label{tab:nbinatonevariance}
	
		\resizebox{\columnwidth}{!}{%
	\centering
	\begin{tabular}{|c|c|c|c|c|c|c|c|} 	  
		\hline
		Dataset & Filtering method&\makecell{Test \\data type}&Parameters& \makecell{Robust  \\optimization} & Bilevel-meanscore &Bilevel-medianscore &Bilevel-portfolio\\
		\hline
		\multirow{12}{*}{\makecell{A14}}
		&\multirow{6}{*}{\makecell{Bin \\Filtered}}
		
		&\multirow{2}{*}{\makecell{All\\(\# 7010)}} 
		& $\p^*_a$&3730&3724&3687&3693
		\\\cline{4-8}
		&&& $\p^*_b$&3625&3775&3765&3573
		\\\cline{3-8}
		
		&&\multirow{2}{*}{\makecell{Not filtered\\(\# 5199)}} 
		& $\p^*_a$&3350&3317&3265&3273
		\\\cline{4-8}
		&&& $\p^*_b$&3248&3365&3342&3185
		\\\cline{3-8}
		
		&&\multirow{2}{*}{\makecell{Filtered\\(\# 1811)}} 
		& $\p^*_a$&380&407&422&420
		\\\cline{4-8}
		&&& $\p^*_b$&377&410&423&388
		\\\cline{2-8}
		
		&\multirow{6}{*}{\makecell{Observable \\Filtered}}
		
		&\multirow{2}{*}{\makecell{All\\(\# 7010)}} 
		& $\p^*_a$&3730&3724&3687&3693
		\\\cline{4-8}
		&&& $\p^*_o$&3732&3734&3695&3509
		\\\cline{3-8}
		
		&&\multirow{2}{*}{\makecell{Not filtered\\(\# 6707)}} 
		& $\p^*_a$&3675&3660&3624&3630
		\\\cline{4-8}
		&&& $\p^*_o$&3679&3672&3629&3444
		\\\cline{3-8}
		
		&&\multirow{2}{*}{\makecell{Filtered\\(\# 303)}} 
		& $\p^*_a$&55&64&63&63
		\\\cline{4-8}
		&&& $\p^*_o$&53&62&66&65
		\\\hline
		
		\multirow{12}{*}{\makecell{\sherpa}}
		&\multirow{6}{*}{\makecell{Bin \\Filtered}}
		
		&\multirow{2}{*}{\makecell{All\\(\# 792)}} 
		& $\p^*_a$&320&337&371&256
		\\\cline{4-8}
		&&& $\p^*_b$&343&328&345&243
		\\\cline{3-8}
		
		&&\multirow{2}{*}{\makecell{Not filtered\\(\# 588)}} 
		& $\p^*_a$&272&283&317&214
		\\\cline{4-8}
		&&& $\p^*_b$&282&270&292&200
		\\\cline{3-8}
		
		&&\multirow{2}{*}{\makecell{Filtered\\(\# 204)}} 
		& $\p^*_a$&48&54&54&42
		\\\cline{4-8}
		&&& $\p^*_b$&61&58&53&43
		\\\cline{2-8}
		
		&\multirow{6}{*}{\makecell{Observable \\Filtered}}
		
		&\multirow{2}{*}{\makecell{All\\(\# 792)}} 
		& $\p^*_a$&320&337&371&256
		\\\cline{4-8}
		&&& $\p^*_o$&286&348&386&252
		\\\cline{3-8}
		
		&&\multirow{2}{*}{\makecell{Not filtered\\(\# 727)}} 
		& $\p^*_a$&304&319&355&237
		\\\cline{4-8}
		&&& $\p^*_o$&271&331&370&235
		\\\cline{3-8}
		
		&&\multirow{2}{*}{\makecell{Filtered\\(\# 65)}} 
		& $\p^*_a$&16&18&16&19
		\\\cline{4-8}
		&&& $\p^*_o$&15&17&16&17
		\\\hline
		
	\end{tabular}
		}
\end{table}

\subsubsection{Comparison of rational approximation and the MC simulator}
\label{sec:mcvsra}
\begin{figure}[h!]
	\centering
	\includegraphics[width=\textwidth]{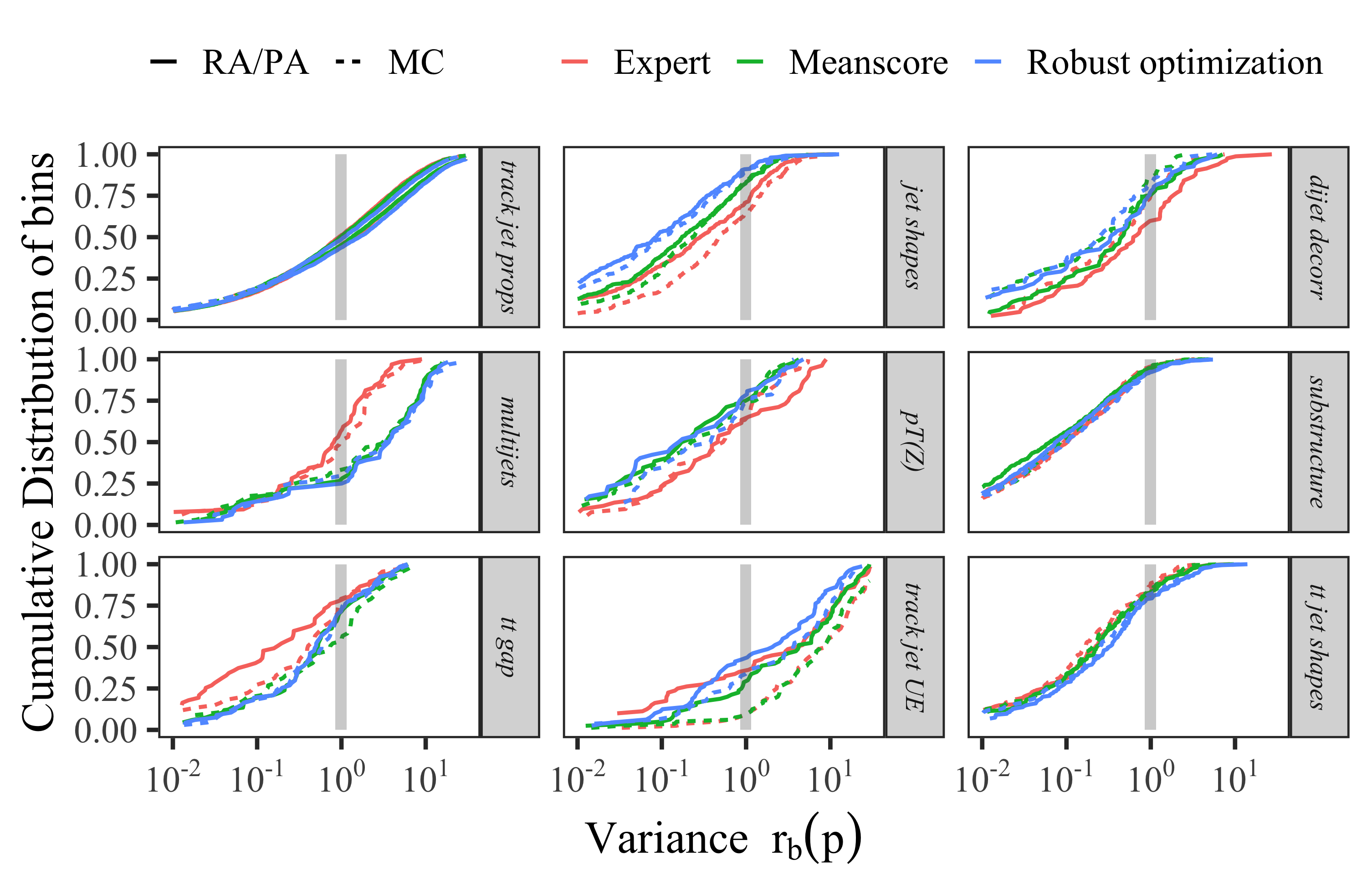}
	\caption{Cumulative distribution of bins (y-axis) in each category of the A14 dataset at different bands of variance levels (x-axis) computed with cubic polynomial approximation (PA) or rational approximation (RA) and the MC simulation.
	}
	\label{Fig:RAvsMCchi2perf2binsCat}
	
\end{figure}

Similar to the analysis conducted in Section~\ref{sec:cdfperfA14},  we compare the cumulative distribution of bins at different bands of variance levels computed using the approximation model as $r_b(\p)=\frac{\left(f_b(\p) - \calR_b\right)^2}{\Delta f_b(\p)^2+\Delta\calR_b^2}$ and the MC generator model as $\widetilde{r_b(\p)}=\frac{\left(\MC_b(\p) - \calR_b\right)^2}{\Delta \MC_b(\p)^2+\Delta\calR_b^2}$, where $\p$ are the parameters obtained from the tuning approaches. The more bins that are on the bands of variance levels less than one, the better. %
Figure~\ref{Fig:RAvsMCchi2perf2binsCat} shows the plot of this comparison for bins in each category of the A14 dataset.\footnote{The \textit{Jet UE} comparison is missing from this figure because the internal ATLAS analysis is not available to us.} 
To avoid making the plot too busy, we show the results using the parameters from three approaches. A similar plot showing the  results with parameters from the remaining approaches is given in Section~\ref{sec:mcvsra2} in the online supplement.

We observe in Figure~\ref{Fig:RAvsMCchi2perf2binsCat} that the \textit{Dijet decorr}, \textit{Jet shapes}, \textit{$p_T^Z$}, \textit{Track-jet UE}, and \textit{$t \Bar{t}$ gap} categories show differences in the performance between $r_b(\p)$ and $\widetilde{r_b(\p)}$ for each approach.
Additionally, for the \textit{robust optimization} and \textit{Bilevel-meanscore} approaches, this difference in the performance is not as wide as that of the \textit{Expert} (for e.g., see $p_T^Z$, \textit{Track-jet UE} categories). This suggests that 
(a) there are categories where the approximations are not able to capture the MC generator  perfectly, and 
(b) in general, the rational approximation is a better surrogate for the MC generator than the polynomial approximation, i.e., the rational approximation gives better predictions of the MC generator than the polynomial approximation.

\subsection{Results for the \sherpa dataset}
In this section, we present the detailed results for the \sherpa dataset.

\subsubsection{Comparison metric outcomes for the \sherpa dataset}

Tables~\ref{tab:comparison_full_sherpa_31}-\ref{tab:comparison_binf_sherpa_31} show the  results when using the rational approximation (results for the cubic polynomial approximation are in the online supplement Section~\ref{sec:sherpa_30_metric}). Smaller numbers indicate better performance. The smallest number of each metric is bold for better visualization. %
Similar to A14, we find that the robust optimization approach achieves the best performance in terms of the Weighted $\chi^2$ criterion. Assigning All-weights-equal to all observables yields the best results in terms of A- and D-optimality for the full and the bin-filtered dataset. The portfolio approach yields the best A- and D-optimality values when using the observable-filtered dataset. 

Compared with the results of A14, we see that the magnitudes of the numbers obtained for the \sherpa dataset for the Weighted $\chi^2$, A- and D-optimality criteria are much larger, indicating that we are not certain about the optima found by the methods.

\begin{table}[htbp]
	\caption{Results for the comparison metrics for the full \sherpa  dataset using the rational approximation.}%
	\centering
	\begin{tabular}{l|rrr}
		\hline
		method         & Weighted $\chi^2$ & A-optimality    & D-optimality (log) \\\hline
		Bilevel-meanscore & 0.2201 &	9.0147 &	-39.3957

 \\
		Bilevel-medscore  & 0.2249 &	43.2031 &	-25.7164
\\
		Bilevel-portfolio & 0.1510	& 11.9869 &	-35.7488

         \\
		All-weights-equal         & 0.2794	&\textbf{6.8428} &	\textbf{-42.0325}

      \\
		Robust optimization          & \textbf{0.0603}	& 55.8079 &	-22.0884

\\\hline 
	\end{tabular}
	\label{tab:comparison_full_sherpa_31}
\end{table}

\begin{table}[htbp]
	\caption{Results for comparison metrics for the  observable-filtered \sherpa dataset using the rational approximation. Three observables were filtered out and not used during the optimization. }
	\centering
	\begin{tabular}{l|rrr}
		\hline
		method         & Weighted $\chi^2$ & A-optimality & D-optimality (log) \\\hline
		Bilevel-meanscore & 0.3621 &	11.1570	& -36.5249

      \\
		Bilevel-medscore  & 0.2315 &	13.0679	& -35.3498
 \\
		Bilevel-portfolio & 0.4728 &	\textbf{8.5578} &	\textbf{-38.6042}

         \\
        All-weights-equal         & 0.4587	& 59.7043 &	-19.1257

\\
		Robust optimization            & \textbf{0.0509} &	32.9470	& -30.5536
 \\\hline        
	\end{tabular}
	\label{tab:comparison_obsf_sherpa_31}
\end{table}

\begin{table}[htbp]
	\caption{Results for the bin-filtered \sherpa dataset using the rational approximation.  204 out of 5426 bins were filtered out and not used during the optimization.}
	\centering
	\begin{tabular}{l|rrr}
		\hline
		method         & Weighted $\chi^2$ & A-optimality & D-optimality (log) \\\hline
        Bilevel-meanscore & 0.1406 &	16.5417 &	-33.3334

          \\
		Bilevel-medscore  & 0.1352	& 16.9715 &	-33.7009

       \\
		Bilevel-portfolio & 0.2792 &	15.2932	& -35.8314

 \\
        All-weights-equal         & 0.2105	& \textbf{8.9591} &	\textbf{-38.3039}

         \\
		Robust optimization            & \textbf{0.0869} &	17.4497	& -34.0525

     	\\\hline
	\end{tabular}
	\label{tab:comparison_binf_sherpa_31}
\end{table}

\subsubsection{Comparison of the cumulative distribution of bins at different variance levels}\label{sec:cdfperfSherpa}

\begin{figure}[h!]
	\centering
	\includegraphics[width=\textwidth]{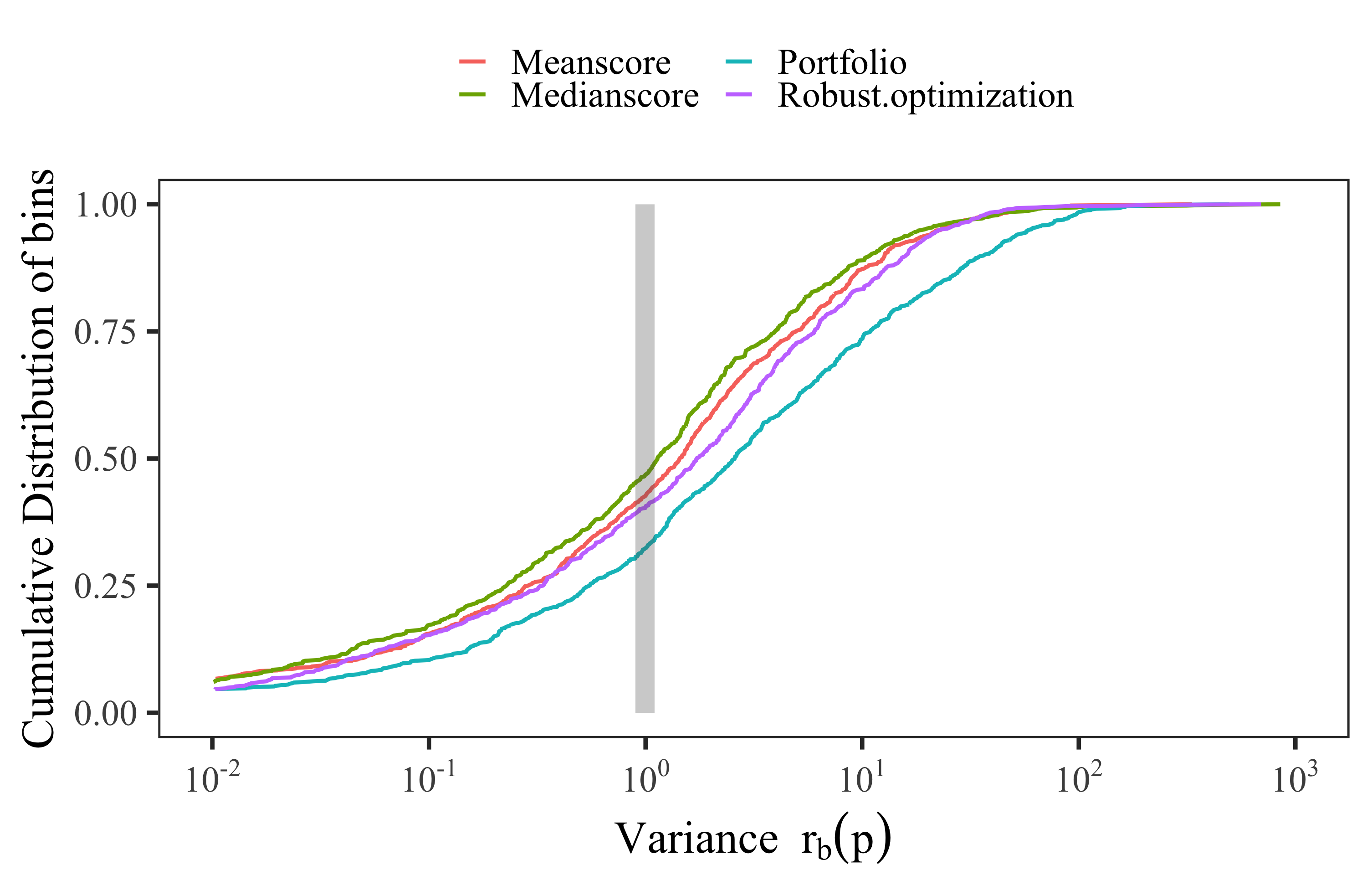}        
	\caption{Cumulative distribution of all bins (y-axis) in the \sherpa dataset at different bands of variance levels (x-axis) given by $r_b(\p)=\frac{\left(f_b(\p) - \calR_b\right)^2}{\Delta f_b(\p)^2+\Delta\calR_b^2}$.}
	\label{fig:chi2perf2binsAllSherpa}
\end{figure}

Similar to the analysis conducted in Section~\ref{sec:cdfperfA14}, we compare the cumulative distribution of bins at different bands of variance level computed using the optimal parameters $\p$ obtained from the tuning approaches. 
Figure~\ref{fig:chi2perf2binsAllSherpa} shows the plot of this comparison for all bins.
The results  show that fewer bins lie within the variance boundary of one when using the  parameters of  the bilevel-portfolio approach. On the other hand, the bilevel-medianscore approach  finds parameters that yield the  most bins at lower bands of variance levels.

\subsubsection{Optimal parameter values for the \sherpa dataset with rational approximation}

The optimal parameter values for the \sherpa dataset  are shown in Tables~\ref{tab:optimal_prms_sherpa_RA}, \ref{tab:optimal_prms_obsfilter_sherpa_RA} and \ref{tab:optimal_prms_binfilter_sherpa_RA}.
For a different visualization of the different solutions obtained with our methods,
we illustrate the [0,1]-scaled optimal parameters in the online supplement Section~\ref{sec:pstar_sherpa_31}.

We see that most of the parameters are on the boundaries of the parameter space (indicated in the table in bold), except for \texttt{KT\_0} and \texttt{STRANGE\_FRACTION}.
This observation indicates that we might need to change the size of the parameter domain  to avoid model extrapolation. 

Note that for the \sherpa dataset, we do not have an ``expert'' solution for benchmark comparison. Instead, we compare the solutions to the chosen reasonable default setting. 
The parameter range is constructed by multiplying the default value by $0.5$ and $1.5$ to obtain the lower and the upper bound of its range respectively, i.e., the default values lie  in the middle of the parameter range. We see that there are differences between the optimal parameters obtained with the different methods, in particular, bilevel-medianscore gives a very similar solution to the default setting.

\begin{table}[htp!]
	\centering
	\small
	\caption{Optimal parameter values obtained with all methods using rational approximation when no filtering was applied before optimization (88 observables). The parameter values on the boundaries of the parameter space are indicated in bold. }\label{tab:optimal_prms_sherpa_RA}
	\begin{adjustbox}{width=\textwidth}
	\begin{tabular}{llllllll}
		\hline
	ID	&	Parameter name	& Default & Bilevel-meanscore & Bilevel-medscore & Bilevel-portfolio & Robust opt & All-weights-equal   \\\hline
1  & \tt KT\_0                  & 1.00 & 0.888          & 0.789          & 0.919          & 0.909 & 0.872          \\
2  & \tt ALPHA\_G               & 1.25 & \textbf{0.626} & 1.500          & \textbf{0.626} & 1.874 & \textbf{0.626} \\
3  & \tt ALPHA\_L               & 2.50 & \textbf{3.749} & 1.890          & \textbf{3.749} & 1.252 & \textbf{3.749} \\
4  & \tt BETA\_L                & 0.10 & \textbf{0.150} & \textbf{0.050} & 0.087          & 0.150 & \textbf{0.150} \\
5  & \tt GAMMA\_L               & 0.50 & 0.274          & 0.339          & \textbf{0.750} & 0.683 & 0.293          \\
6  & \tt ALPHA\_H               & 2.50 & 3.400          & 2.897          & \textbf{1.251} & 2.841 & 3.440          \\
7  & \tt BETA\_H                & 0.75 & 0.827          & 0.536          & 0.783          & 0.540 & 0.795          \\
8  & \tt GAMMA\_H               & 0.10 & 0.148          & \textbf{0.050} & 0.082          & 0.150 & \textbf{0.150} \\
9  & \tt STRANGE\_FRACTION      & 0.50 & 0.517          & 0.498          & 0.583          & 0.508 & 0.546          \\
10 & \tt BARYON\_FRACTION       & 0.18 & 0.100          & 0.175          & 0.106          & 0.136 & \textbf{0.090} \\
11 & \tt P\_QS\_by\_P\_QQ\_norm & 0.48 & \textbf{0.720} & 0.419          & 0.572          & 0.613 & \textbf{0.720} \\
12 & \tt P\_SS\_by\_P\_QQ\_norm & 0.02 & \textbf{0.010} & 0.015          & \textbf{0.030} & 0.030 & \textbf{0.010} \\
13 & \tt P\_QQ1\_by\_P\_QQ0     & 1.00 & \textbf{1.499} & 1.206          & 0.948          & 1.190 & \textbf{1.499}
 \\\hline
\hline
& Euclidean distance from the default solution & & 1.513 & 	0.984 & 	1.244 & 	1.289 & 	1.531

 \\\hline %
					\end{tabular}
					\end{adjustbox}
\end{table}

\begin{table}[htbp]
	\small
	\centering
	\caption{Optimal parameter values for the \sherpa dataset obtained with all methods using rational approximation after observable filtering  (3 observables were filtered out). The parameter values on the boundaries of the parameter space are indicated in bold.}\label{tab:optimal_prms_obsfilter_sherpa_RA}
	\begin{adjustbox}{width=\textwidth}
	\begin{tabular}{llllllll}
		\hline
	ID	&	Parameter name	& Default & Bilevel-meanscore & Bilevel-medscore & Bilevel-portfolio & Robust opt & All-weights-equal   \\\hline
1  & \tt KT\_0                  & 1.00 & 0.867          & 0.744          & 0.952          & 0.876          & 0.886 \\
2  & \tt ALPHA\_G               & 1.25 & 0.775          & \textbf{0.626} & \textbf{0.626} & 0.626          & 0.957 \\
3  & \tt ALPHA\_L               & 2.50 & \textbf{3.749} & \textbf{1.252} & \textbf{3.749} & 1.252          & 2.424 \\
4  & \tt BETA\_L                & 0.10 & 0.109          & \textbf{0.050} & \textbf{0.050} & 0.150          & 0.113 \\
5  & \tt GAMMA\_L               & 0.50 & \textbf{0.250} & 0.437          & 0.413          & 0.750          & 0.460 \\
6  & \tt ALPHA\_H               & 2.50 & 3.053          & 2.318          & \textbf{1.251} & 2.826          & 3.132 \\
7  & \tt BETA\_H                & 0.75 & 0.827          & 0.625          & 0.750          & 0.375          & 0.969 \\
8  & \tt GAMMA\_H               & 0.10 & \textbf{0.050} & 0.134          & 0.094          & 0.050          & 0.131 \\
9  & \tt STRANGE\_FRACTION      & 0.50 & 0.479          & 0.580          & 0.651          & 0.506          & 0.511 \\
10 & \tt BARYON\_FRACTION       & 0.18 & \textbf{0.270} & 0.137          & \textbf{0.090} & 0.137          & 0.180 \\
11 & \tt P\_QS\_by\_P\_QQ\_norm & 0.48 & \textbf{0.720} & 0.469          & 0.495          & 0.470          & 0.601 \\
12 & \tt P\_SS\_by\_P\_QQ\_norm & 0.02 & \textbf{0.010} & \textbf{0.030} & \textbf{0.030} & \textbf{0.030} & 0.019 \\
13 & \tt P\_QQ1\_by\_P\_QQ0     & 1.00 & \textbf{0.500} & \textbf{1.499} & \textbf{1.499} & 1.499          & 0.958
 \\\hline
\hline
& Euclidean distance from the default solution & & 1.408 & 	1.249 & 	1.372 & 	1.446 & 	0.637

 \\\hline %
 \end{tabular}
	\end{adjustbox}
\end{table}

\begin{table}[htbp]
	\centering
	\small
	\caption{Optimal parameter values for the \sherpa dataset obtained with all methods using rational approximation after  bin-filtering  (7 bins were filtered out). The parameter values on the boundaries of the parameter space are indicated in bold.}\label{tab:optimal_prms_binfilter_sherpa_RA}
	\begin{adjustbox}{width=\textwidth}
	\begin{tabular}{llllllll}
		\hline
	ID	& Parameter name		& Default & Bilevel-meanscore & Bilevel-medscore & Bilevel-portfolio & Robust opt & All-weights-equal   \\\hline
1  & \tt KT\_0                  & 1.00 & 0.895          & 0.821          & 0.948          & 0.820 & 0.899          \\
2  & \tt ALPHA\_G               & 1.25 & 0.893          & 1.483          & \textbf{0.626} & 1.874 & \textbf{0.626} \\
3  & \tt ALPHA\_L               & 2.50 & \textbf{3.749} & 2.334          & 2.567          & 3.749 & \textbf{3.749} \\
4  & \tt BETA\_L                & 0.10 & 0.050          & \textbf{0.150} & 0.074          & 0.050 & 0.067          \\
5  & \tt GAMMA\_L               & 0.50 & 0.390          & \textbf{0.250} & \textbf{0.750} & 0.250 & 0.454          \\
6  & \tt ALPHA\_H               & 2.50 & \textbf{1.251} & 3.670          & \textbf{1.251} & 1.969 & \textbf{1.251} \\
7  & \tt BETA\_H                & 0.75 & 0.715          & 0.534          & 0.739          & 1.125 & 0.715          \\
8  & \tt GAMMA\_H               & 0.10 & 0.119          & 0.142          & 0.105          & 0.050 & 0.089          \\
9  & \tt STRANGE\_FRACTION      & 0.50 & 0.556          & 0.542          & 0.570          & 0.531 & 0.559          \\
10 & \tt BARYON\_FRACTION       & 0.18 & 0.122          & 0.120          & 0.124          & 0.138 & 0.124          \\
11 & \tt P\_QS\_by\_P\_QQ\_norm & 0.48 & 0.595          & \textbf{0.720} & 0.492          & 0.497 & 0.577          \\
12 & \tt P\_SS\_by\_P\_QQ\_norm & 0.02 & \textbf{0.030} & 0.030          & \textbf{0.030} & 0.030 & \textbf{0.030} \\
13 & \tt P\_QQ1\_by\_P\_QQ0     & 1.00 & \textbf{1.499} & \textbf{1.499} & \textbf{1.499} & 1.499 & \textbf{1.499}       
 \\\hline
\hline
		& Euclidean distance from the default solution & & 1.266 & 	1.377 & 	1.201	 & 1.462 & 	1.242

 \\\hline %
 \end{tabular}

	\end{adjustbox}
      \end{table}

      The distribution of weights from the different methods has a similar pattern as for the tunes based on the A14 dataset.
      These patterns are displayed in Fig.~\ref{fig:sherpa_weights} in the online supplement.
      Robust optimization selects only one of the event shape observables as relevant, while applying the same equal weight to
      most of the particle multiplicity (one bin) distributions.   The other methods have weights that are more widely distributed
      among the observables with a small number of weights far from the average.

\subsection{A note on computation times}

The bilevel optimization approaches of medianscore, meanscore, and portfolio are run on a 4-core, 32 GB RAM machine running at 1.1 GHz. For the results of robust optimization  presented in this paper, 100 values for $\mu$ are used that are run on 100 threads in parallel on a server with 64 Intel Xeon Gold CPU cores running at 2.30 GHz. There are two threads per core, but each run of robust optimization is run on a single thread. Additionally, this server is equipped with 1.5TB DDR4 2666 MHz of memory. A simple comparison to find the best $\mu$ takes one minute. 
The all-weights-equal approach is run on a 4-core, 32 GB RAM machine running at 1.1 GHz.

The time taken by all the tuning approaches for unfiltered (\textit{All data}) as well as for bin filtered and observable filtered A14 data is given in Table~\ref{tab:filterApproachPerfromanceA14}.
In the unfiltered data case, the bilevel optimization approaches of medianscore, meanscore, and portfolio  take approximately 14.5 hours and
each run (i.e., one $\mu$) of robust optimization takes an average of about 0.8 hours. Since all 100 values of $\mu$ were run in parallel, the total time to complete all 100 runs of robust optimization is approximately two hours.
In comparison, campaigns to tune weights by hand takes many weeks or months.
Given our results, we can see that the automated weight adjustment by optimization is significantly faster than hand-tuning. 
The all-weights-equal approach took less than 10 minutes, however, this approach leads to worse results.

The observable filtering method requires a single-tune to obtain the $\chisq$ values per observable which takes 1647 seconds (0.45 hours) for all observables in the A14 dataset, which is followed by applying the Z-score method to filter out outliers (see Section~\ref{sec:obsfilter}) and this takes about 10 seconds. 
Once the single-tune to obtain the $\chisq$ values per observable is performed, the bin-filtering method takes an additional 300 seconds to filter out the bins from the A14 dataset. 
Thus, the total pre-processing time required for observable filtering is 1657 seconds (0.46 hours) and for bin-filtering is 1947 seconds (0.54 hours).

From Table~\ref{tab:filterApproachPerfromanceA14}, we observe that the time taken to tune parameters in the observable filtered and bin filtered data case is significantly smaller than for the unfiltered data case. 
For the bilevel optimization approaches, the time required per iteration for the observable- and bin-filtered cases is 6\%, and 55\% less, respectively, and for each run of robust optimization, it is 9\% and 36\% less, respectively. Also, the overhead of performing observable and bin filtering is small compared to the time it takes to tune parameters. 
Since the results from Section~\ref{sec:filterresults} show that the bins filtered by bin and observable filtering do not add significant information to the tune, we can claim that using filtered data provides a significant improvement in compute-time performance for tuning parameters.

\begin{table}[h!]
	\caption{CPU time (in seconds) and time per iteration (in seconds) taken by all approaches when using  all, the observable-filtered, and the bin-filtered  A14 data. The robust optimization approach converges after  69, 105, and 83 iterations, respectively. The bilevel-medianscore, -meanscore, and -portfolio approaches are all run for 1000 iterations.
	}
	\label{tab:filterApproachPerfromanceA14}
	
	\centering
	\begin{tabular}{|c|r|r|r|r|r|r|}
		\hline
		
		\multirow{3}{*}{Method} & \multicolumn{2}{c|}{All data} & \multicolumn{2}{c|}{Bin filtered} & \multicolumn{2}{c|}{Observable filtered}\\
		\cline{2-7}
		& CPU time & \makecell{Time per \\iteration} & CPU time & \makecell{Time per \\iteration}& CPU time & \makecell{Time per \\iteration}
		\\\hline
		\makecell{Robust optimization} &  3035 & 44 & 2989 & 28 & 3327 & 40 \\
		Bilevel-medianscore &52326&52&23600&24&49057&49\\
		Bilevel-meanscore &52169&52&23600&24&49018&49\\
		Bilevel-portfolio &52366&52&23609&24&49084&49\\\hline
	\end{tabular}
\end{table}

\section{Eigentunes}
\label{sec:eigentune}
We use the eigentune approach to calculate confidence intervals for the optimal parameters. 
We note that the A- and D-optimality criteria provide the size of confidence ellipsoid around the optimal parameters. Here, we expand this information by scanning generator parameters along the principal axes of this ellipsoid.
Details of this method are described in \cite{Buckley:2009bj} and  a similar approach is used in estimating the uncertainties of predictions from the parton distribution functions~\cite{Pumplin:2001ct}. The interval defines a boundary beyond which the value of the objective function is larger than the objective function value at the minimum by a criterion. The criterion is normally chosen to be the number of degrees of freedom $n$, which is defined as the total number of bins of all observables minus the number of generator parameters, $d$, i.e., $ n = \sum_{\calO \in \calS_\calO}\abscalO - d$. %
However, to properly take into account the weights assigned to observables, we use the scaled effective sample size as the criteria, which is calculated as follows:
\begin{equation*}
    n = \gamma \times \left(\frac{(\sum_i w_i)^2}{\sum_i w_{i}^{2}} - d\right)
\end{equation*}
The weights are normalized so that the sum of weights associated with all observables equals one. $\gamma$ is iteratively tuned and chosen to be 0.01.
The interval would represent the uncertainties of the  parameters assuming that the objective function follows a $\chi^2$ distribution. Smaller intervals associated with the tuned parameters indicate that the parameters are better constrained by the experimental data.

Given the non-linearity of the objective function and parameter correlations, a reliable approach to find the 68\% confidence interval is to evaluate the objective function for all possible parameter values. However, this  poses a computational challenge. Instead, we project the multidimensional  parameter space into two directions defined by the eigenvectors $u_{1,2}$ associated with the largest and smallest eigenvalues of the covariance matrix of the parameters, which are calculated using the inverse of \equaref{eq:gammapost}. Then we find an offset $\alpha$ such that the sum of all $\chi^2$ satisfies
\begin{equation}
    \chi^2(\Vec{p}^\prime_{1,2}) = \chi^2(\Vec{p}^*) + n
    \label{eq:eigen}
\end{equation}
where $\Vec{p}^\prime_{1,2} = \Vec{p}^* \pm u_{1,2}\times\alpha$. For each eigenvector, we obtain  two vectors  $\Vec{p}^\prime$ from  Eq.~(\ref{eq:eigen}). Finally, the procedure results in a matrix of sizes of 4 times $d$. Each column represents a generator parameter; the minimum and maximum in each column  are used to define the eigentune as shown in Tables~\ref{tab:eigentune_a14_ra} and~\ref{tab:eigentune_sherpa_ra} for the A14 and the \sherpa dataset, respectively, using the rational approximation. The same surrogate model is used for all methods. It is possible that the determined intervals  go beyond the predefined parameter range. In this case, the MC predictions are extrapolated by the surrogate model. When the lower part of the interval goes negative, we force the value  to be zero.

For the A14 data, different optimization methods result in similar intervals for all  parameters. The beam remnants (e.g. {\tt BeamRemnants:reconnectRange}) and space-like showering parameters (e.g. {\tt SpaceShower:pT0Ref}) are better constrained; their intervals are within 1\% of their optimized parameters. However, the strong coupling constant $\alpha_s$ in hard scattering processes ({\tt SigmaProcess:alphaSvalue}) and time-like showering ({\tt TimeShower:alphaSvalue}) are less constrained. 

For the \sherpa data, different optimization methods produce quite different intervals. Overall, the bilevel-meanscore method results in relatively small intervals for all  parameters. The heavy quark fragmentation parameters (e.g. {\tt ALPHA\_H}) are well-constrained thanks to the $B$-hadron fragmentation measurements, but the light quark fragmentation parameters are not.

\begin{table}[!htbp]
\caption{Eigentune results for the A14 data using the rational approximation for different optimization methods.}
\label{tab:eigentune_a14_ra}
\begin{adjustbox}{width=\textwidth}
\begin{tabular}{l|rr|rr|rr|rr|rr}\hline
    Parameters & \multicolumn{2}{c|}{Expert}&  \multicolumn{2}{c|}{Bilevel-meanscore} & \multicolumn{2}{c|}{Bilevel-mediansocre} & \multicolumn{2}{c|}{Bilevel-portfolio} & \multicolumn{2}{c}{Robust optimization}  \\ 
     &  min & max &  min & max  &  min & max &  min & max &  min & max\\  \hline 
\tt SigmaProcess:alphaSvalue &  0.075 &  0.193 &  0.079 &  0.192 &  0.079 &  0.190 &  0.074 &  0.195 &  0.085 &  0.183 \\
\tt BeamRemnants:primordialKThard &  1.903 &  1.906 &  1.805 &  1.910 &  1.674 &  1.769 &  1.744 &  1.850 &  1.876 &  1.892 \\
\tt SpaceShower:pT0Ref &  1.636 &  1.653 &  1.516 &  1.547 &  1.142 &  1.228 &  1.298 &  1.344 &  1.586 &  1.591 \\
\tt SpaceShower:pTmaxFudge &  0.905 &  0.912 &  1.012 &  1.016 &  1.069 &  1.096 &  1.037 &  1.046 &  1.025 &  1.026 \\
\tt SpaceShower:pTdampFudge &  1.044 &  1.048 &  1.064 &  1.076 &  1.082 &  1.086 &  1.058 &  1.064 &  1.078 &  1.091 \\
\tt SpaceShower:alphaSvalue &  0.121 &  0.124 &  0.125 &  0.131 &  0.127 &  0.130 &  0.124 &  0.133 &  0.123 &  0.129 \\
\tt TimeShower:alphaSvalue &  0.043 &  0.197 &  0.044 &  0.192 &  0.039 &  0.213 &  0.030 &  0.213 &  0.051 &  0.198 \\
\tt MultipartonInteractions:pT0Ref &  1.665 &  2.543 &  1.649 &  2.562 &  1.780 &  1.979 &  1.160 &  2.829 &  1.461 &  2.528 \\
\tt MultipartonInteractions:alphaSvalue &  0.068 &  0.177 &  0.072 &  0.161 &  0.115 &  0.121 &  0.062 &  0.186 &  0.094 &  0.151 \\
\tt BeamRemnants:reconnectRange &  1.788 &  1.795 &  2.065 &  2.105 &  1.912 &  1.915 &  1.972 &  2.000 &  2.589 &  2.618 \\
\hline
\end{tabular}
\end{adjustbox}
\end{table}

\begin{table}[!htbp]
\caption{Eigentune results for the \sherpa data using the rational approximation for different optimization methods. Parameters with negative values are set to zero.}
\label{tab:eigentune_sherpa_ra}
\begin{adjustbox}{width=\textwidth}
\begin{tabular}{l|rr|rr|rr|rr}\hline
    Parameters & \multicolumn{2}{c|}{Bilevel-meanscore} & \multicolumn{2}{c|}{Bilevel-mediansocre} & \multicolumn{2}{c|}{Bilevel-portfolio} & \multicolumn{2}{c}{Robust optimization}  \\ 
     &  min & max &  min & max &  min & max &  min & max\\  \hline 

\tt KT\_0 &  0.815 &  0.970 &  0.688 &  0.957 &  0.524 &  1.254 &  0.491 &  1.273 \\
\tt ALPHA\_G &  0.438 &  0.792 &  1.325 &  1.604 &  0.571 &  0.691 &  1.597 &  2.115 \\
\tt ALPHA\_L &  3.683 &  3.824 &  1.309 &  2.863 &  3.525 &  3.939 &  0.291 &  2.088 \\
\tt BETA\_L & 0 &  0.460 &  0.043 &  0.062 & 0 &  0.440 & 0 &  0.387 \\
\tt GAMMA\_L &  0.175 &  0.362 &  0.330 &  0.352 &  0.688 &  0.823 &  0.220 &  1.087 \\
\tt ALPHA\_H &  3.245 &  3.537 &  2.843 &  2.988 &  1.200 &  1.311 &  2.289 &  3.475 \\
\tt BETA\_H &  0.747 &  0.898 &  0.484 &  0.585 &  0.623 &  0.972 &  0.350 &  0.759 \\
\tt GAMMA\_H &  0.059 &  0.249 & 0 &  0.080 &  0.013 &  0.133 & 0 &  0.469 \\
\tt STRANGE\_FRACTION &  0.496 &  0.556 &  0.395 &  0.595 &  0.415 &  0.706 &  0.440 &  0.567 \\
\tt BARYON\_FRACTION & 0 &  0.459 &  0.129 &  0.218 &  0.018 &  0.170 & 0 &  0.342 \\
\tt P\_QS\_by\_P\_QQ\_norm &  0.552 &  0.809 &  0.319 &  0.524 &  0.552 &  0.588 &  0.594 &  0.629 \\
\tt P\_SS\_by\_P\_QQ\_norm & 0. &  0.031 & 0. &  0.103 & 0 &  0.081 & 0 &  0.068 \\
\tt P\_QQ1\_by\_P\_QQ0 &  1.492 &  1.512 &  1.202 &  1.210 &  0.945 &  0.952 &  1.167 &  1.210 \\

\hline
\end{tabular}
\end{adjustbox}
\end{table}

\section{Discussion}

The results presented in the previous sections demonstrate that automated
tuning methods can produce better fits of the generator predictions to data.
Several figures of merit for comparing different tunes were considered. The automation of the process means
that tuning can be performed in less time and with less subjective bias.
In this section, we discuss 
the physics impact of various tuning results.  

\subsection{Implications of our results on physics}

Physics event generators are imperfect tools.    They contain a mixture of solid physics predictions, approximations, and
{\it ad hoc} models.   The approximations and models are expected to be incomplete, and thus are unlikely to describe
the full range of observables accessible by the experiment.   Despite this fact, for a certain choice of parameters, a model
may be able to describe part of the data.   This agreement would be accidental and would likely compromise predictions
of this model for different parts of the data.     The weighting of data by an expert is a primitive attempt to force the model to
agree with data in a region of interest to the physicist -- which, most of the time, corresponds to a region where a model
should be applied.  It is equivalent to adding a large systematic uncertainty to the
data that is de-emphasized by the weighting.

Here, we address whether the automated methods accomplish this weighting of data without explicit input from the
physicist.     First, we should state our expectations for a tune to the A14 dataset.
The features of the expert tune were previously discussed in Section.~2.2.1 of the A14 publication \cite{ATL-PHYS-PUB-2014-021}.
The A14 data is all of interest to the physicist, but some of those observables are expected {\it a priori} to
be described better by the event generator than others.   The parton shower and hadronization model are expected to describe well
{\it Tracked jet properties} and {\it Jet shapes}.   The description of jets is essential for all hadron collider analyses and is the {\it raison d'\^etre} for
event generators.
{\it $t\bar t$ jet shapes} emphasize the final state parton shower, and is critical to be described well when making
precision predictions that are sensitive to the top quark mass.
{\it Dijet decorr} and {\it $p_T^Z$} observables provide constraints on initial state parton shower and intrinsic transverse momentum parameters free
from most other parameters, and are generically important to be described well.
Additional properties, such as the number of jets produced in di-jet or $Z$ events or the production of jets at extreme angles,
are beyond the scope of the \pythia predictions.
{\it Track-jet UE} and {\it Jet UE} observables
are sensitive to \pythia's multi-parton-interaction model, which describes most of the particles
produced in a high-energy collision.   The addition of {\it Multijets} observables is biasing the parton shower to describe a
next-to-leading order observable, while the leading-logarithm parton shower includes only an approximation to the full
result.   Experience shows that this biasing provides a globally better description of many observables
of interest to the physicist with little effort and without significantly impacting other predictions.   This feature was built into the \textit{Expert} tune by applying
a large weight to this dataset.
Finally, adding the
{\it $t\bar t$ gap} category is asking for the description of an exclusive observable, which has very strong requirements in its construction,
whereas the \pythia prediction here is valid for more inclusive observables.   Including this data in the tune is a
very specific physics requirement that may be beyond the scope of the \pythia approximations.

\subsection{Observables with improved descriptions}

	\begin{figure}[htbp]
	\begin{subfigure}[b]{0.5\textwidth}
         \centering
         \includegraphics[width=\textwidth]{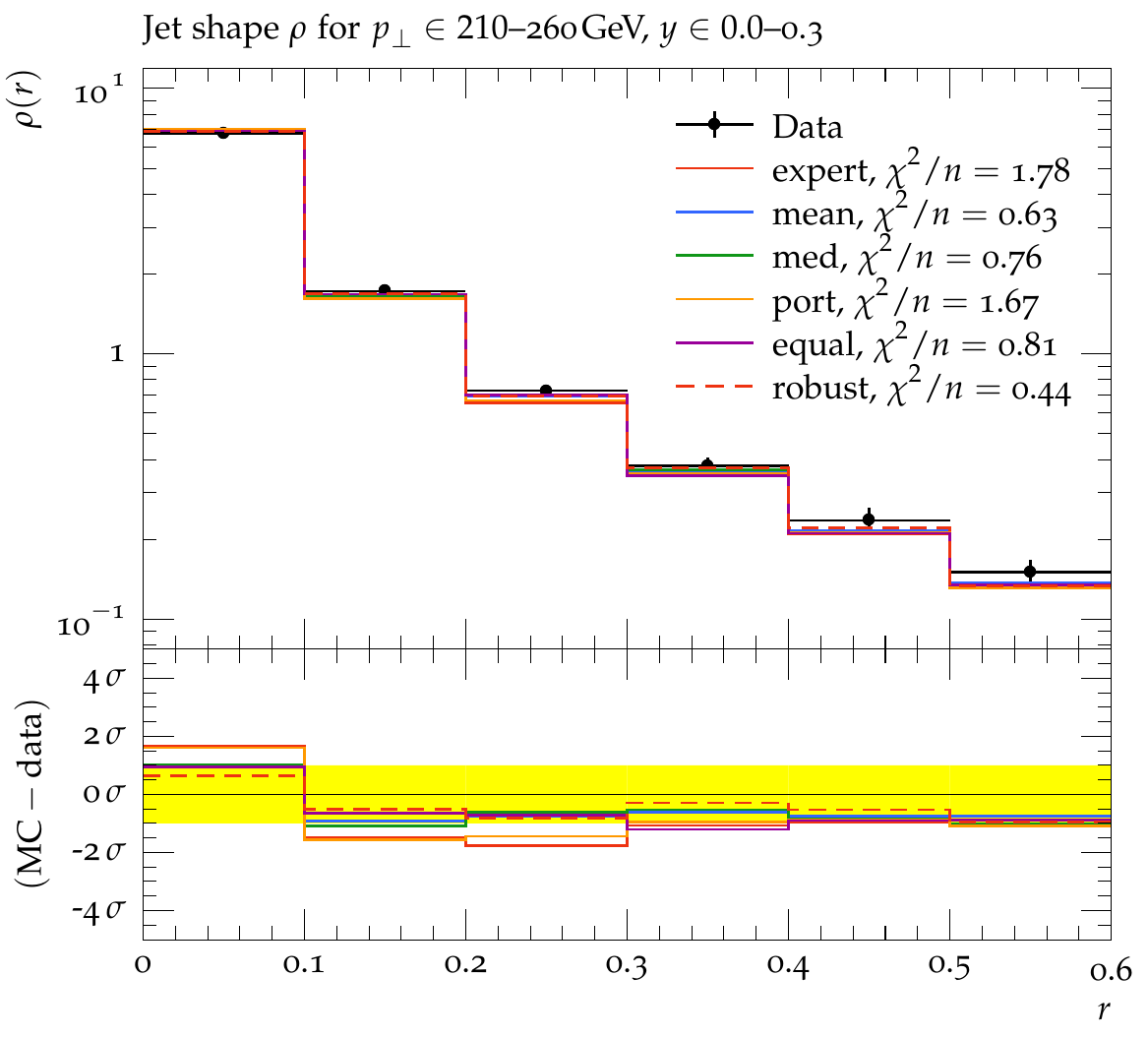}
         \caption{\it  Jet shapes}
         \label{fig:improved-a14-ex1}
     \end{subfigure}
     \hfill
     \begin{subfigure}[b]{0.5\textwidth}
         \centering
         \includegraphics[width=\textwidth]{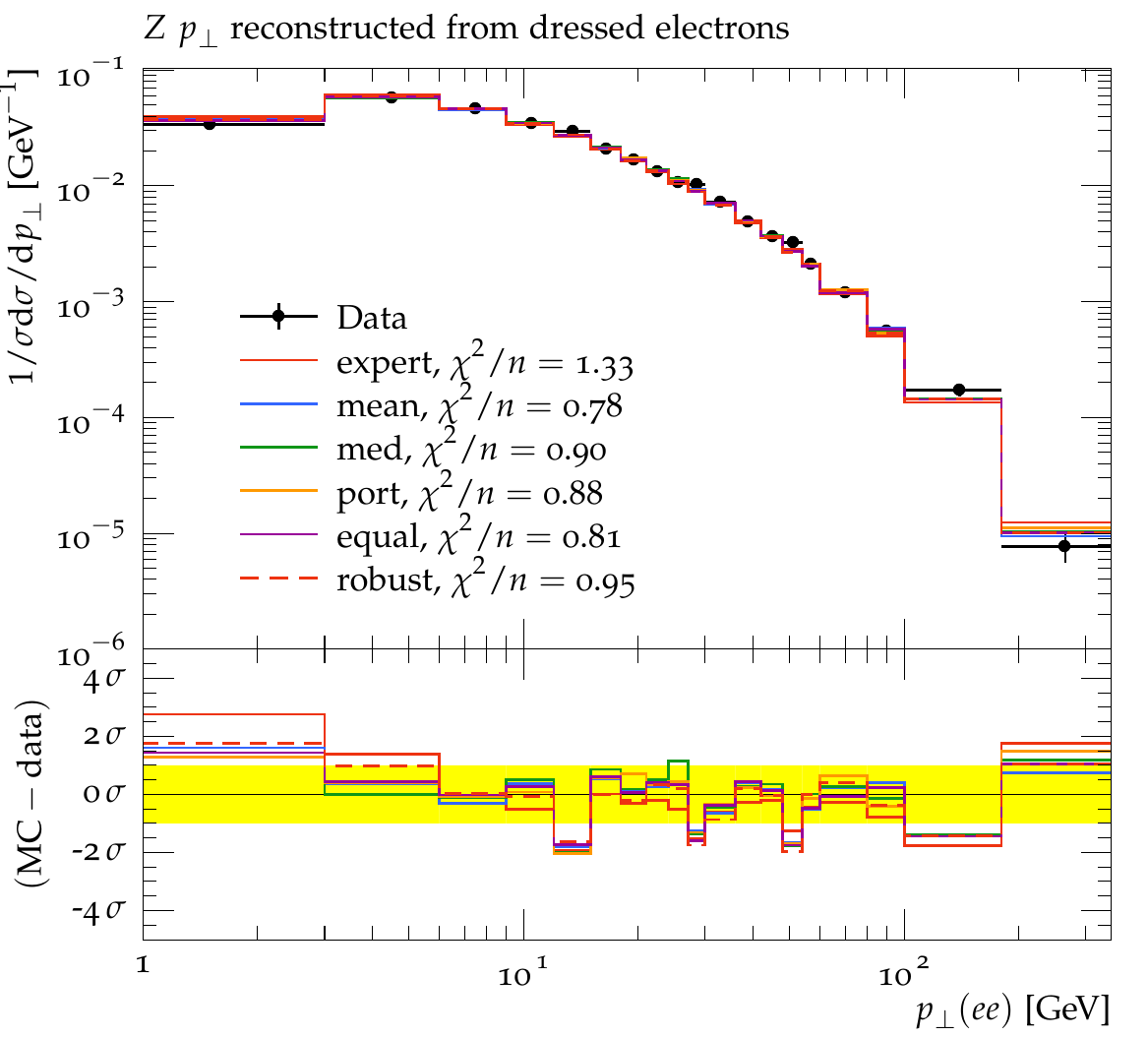}
         \caption{\it $p^Z_T$}
         \label{fig:improved-a14-ex2}
     \end{subfigure}
          \begin{subfigure}[b]{0.5\textwidth}
         \centering
         \includegraphics[width=\textwidth]{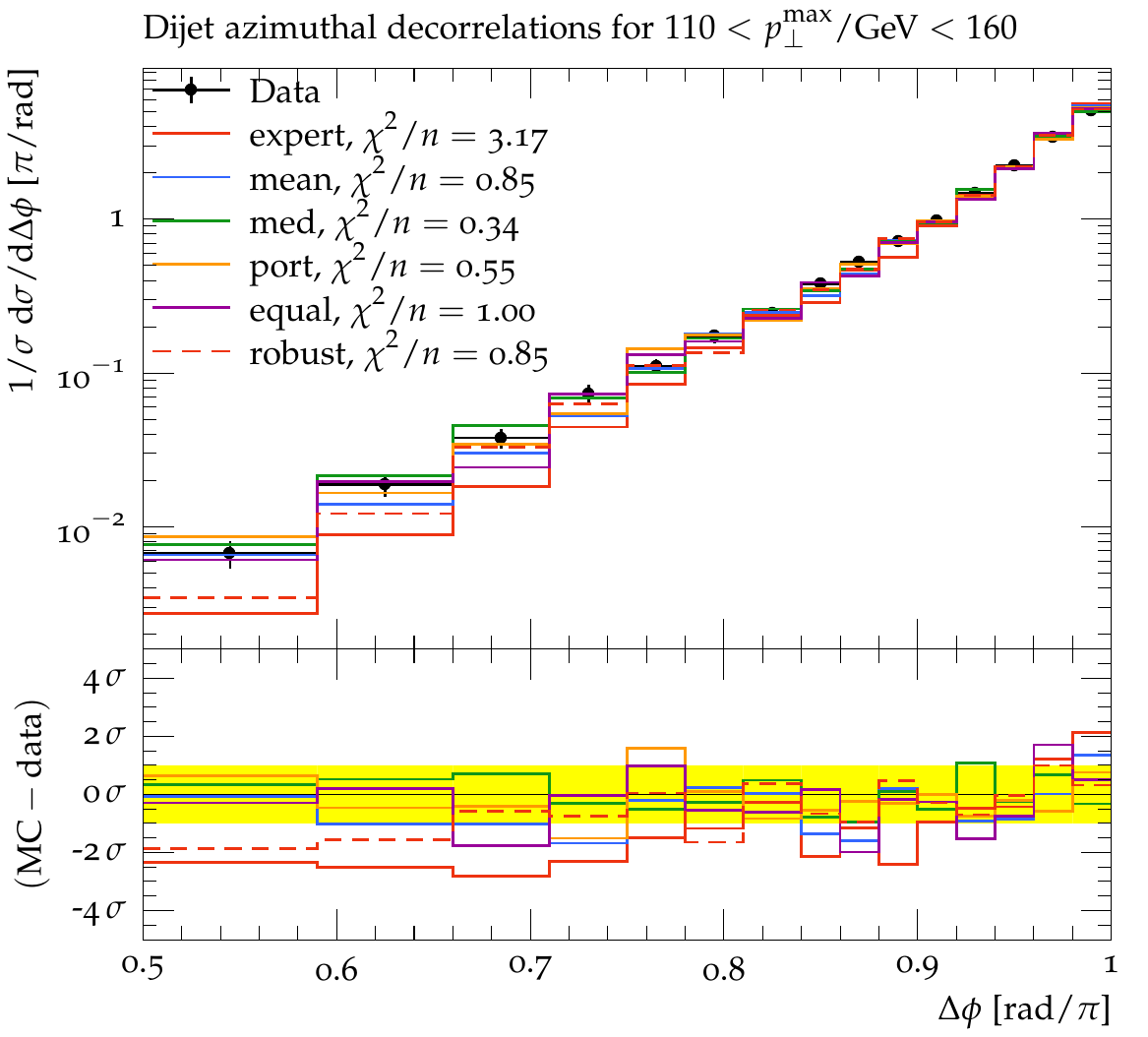}
         \caption{\it Dijet decorr}
         \label{fig:improved-a14-ex3}
     \end{subfigure}
 	\caption{Examples of A14 observables and their $\chi^2_\calO$ values for which the automated tuning leads to better fits than the expert's hand tuning.}
\end{figure}

Examples of observable predictions with a lower $\chi^2_\calO$ value than the {\it
  expert} tune are displayed in Figures~\ref{fig:improved-a14-ex1}-\ref{fig:improved-a14-ex3}.
These reflect an improvement in a class of observables and are
indicative of all the comparisons between predictions and data.

All of our methods produce a better description of the data than the
expert tune for the category {\it Jet shapes}, though the expert
prediction is mainly differing in only the first bin.
This observable is expected to be described well, in general, since it
lies in a physics regime compatible with the \pythia approximations.

The predictions for the {\it $p_T^Z$} and {\it Dijet decorr}
categories are also improved.   We note that the weights found for these analyses are not substantially
different than for the expert tune, but that other categories have their weights reduced (see Table~\ref{tab:weights} for reference).
This implies some tension between
these observables and the \textit{Multijets} category (to be discussed below).

The comparisons between predictions and data shown in our figures are
based on runs of the MC event generator for the parameter
values derived using the surrogate model.
Before continuing, we should comment on the differences in Figure~\ref{Fig:RAvsMCchi2perf2binsCat} (and in Figure~\ref{Fig:RAvsMCchi2perf2binsCat2} in Section~\ref{sec:filterresonlinesup} of the online supplement) between the surrogate model (RA) and
explicit runs of the event generator (MC) at the output tuned
parameters.
The surrogate model would be unreliable if the output tune parameters were outside or near the boundary of the parameter
range used to derive the inputs for the surrogate.   
A comparison of the parameter values relative to the expert tune and
Figure \ref{fig:prms_all_31} shows the distribution of parameter
values normalized to the sampling range:  $\mathrm{r}_\mathrm{param} = \frac{\p - \p_\mathrm{min}}{\p_\mathrm{max}-\p_\mathrm{min}}$.
All of the
central values for the parameters are well within the sampling range.
Only the parameters {\tt   SpaceShower:pTdampFudge} and {\tt BeamRemnants:reconnectRange}
come near the boundaries.  For the former, the minimum sampling value
was 1.0, and the tuning results only indicate that this parameter
should be near 1.0.   For the latter, the maximum sampling value was
chosen quite large so that all results appear to be close to the
minimum value.

Furthermore, the most noticeable differences between the RA surrogate
predictions and MC occur for
rather small values of the variance between the data and predictions.   These values have a negligible impact on the full $\chi^2$, and
are within the expected range of validity of the surrogate model.

\subsection{Observables with worse descriptions}

The predictions for {\it Track jet properties} and {\it Substructure} are not significantly improved,
but also not degraded.
Most of the observables in these categories were
designed to tune and test the multi-parton interaction model, and thus it is no surprise that they are
described well.

\begin{figure}[htbp]
	\begin{subfigure}[b]{0.5\textwidth}
         \centering
         \includegraphics[width=\textwidth]{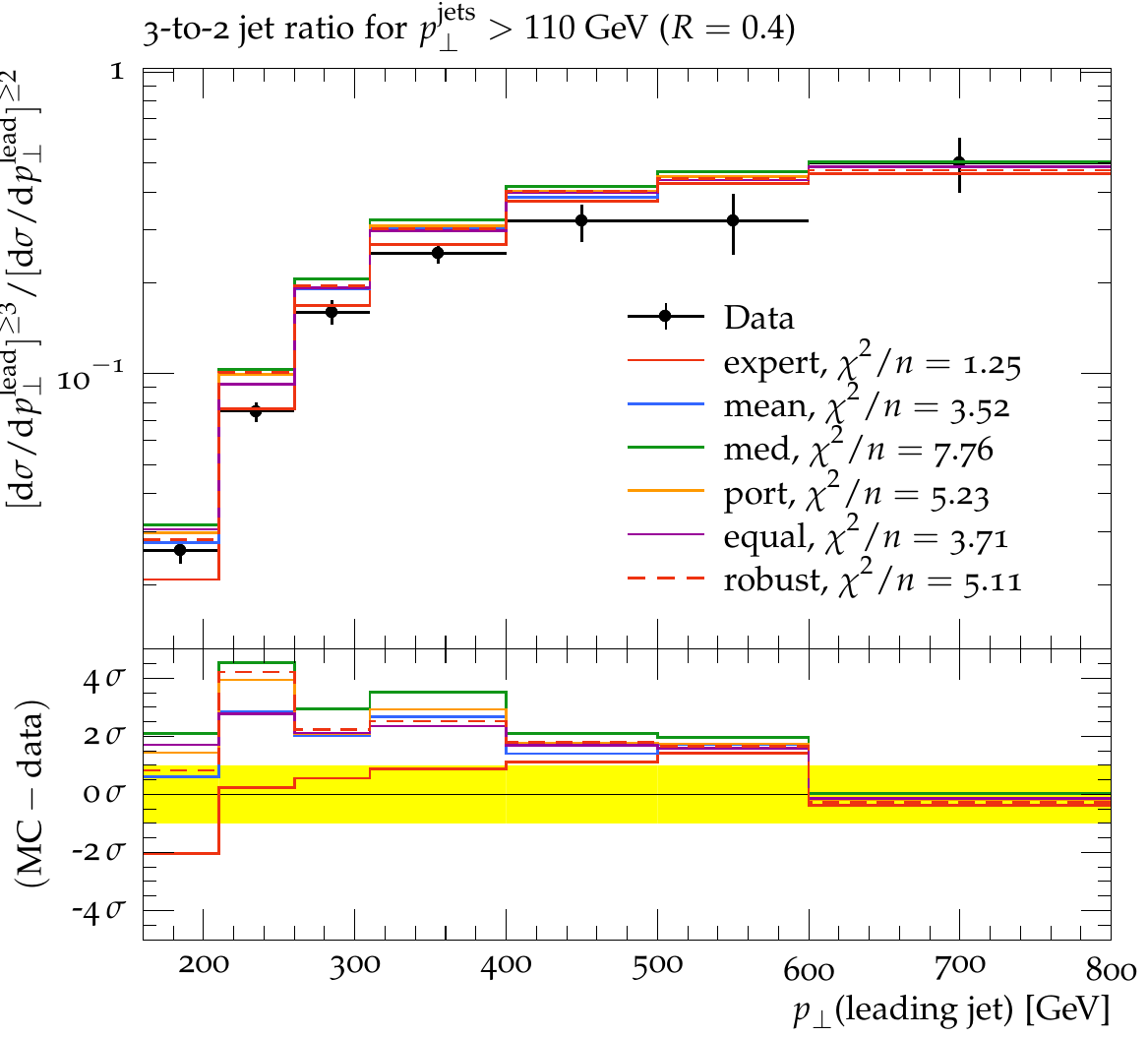}
         \caption{\it Multijets}
         \label{fig:worse-a14-ex1}
     \end{subfigure}
     \hfill
          \begin{subfigure}[b]{0.5\textwidth}
         \centering
         \includegraphics[width=\textwidth]{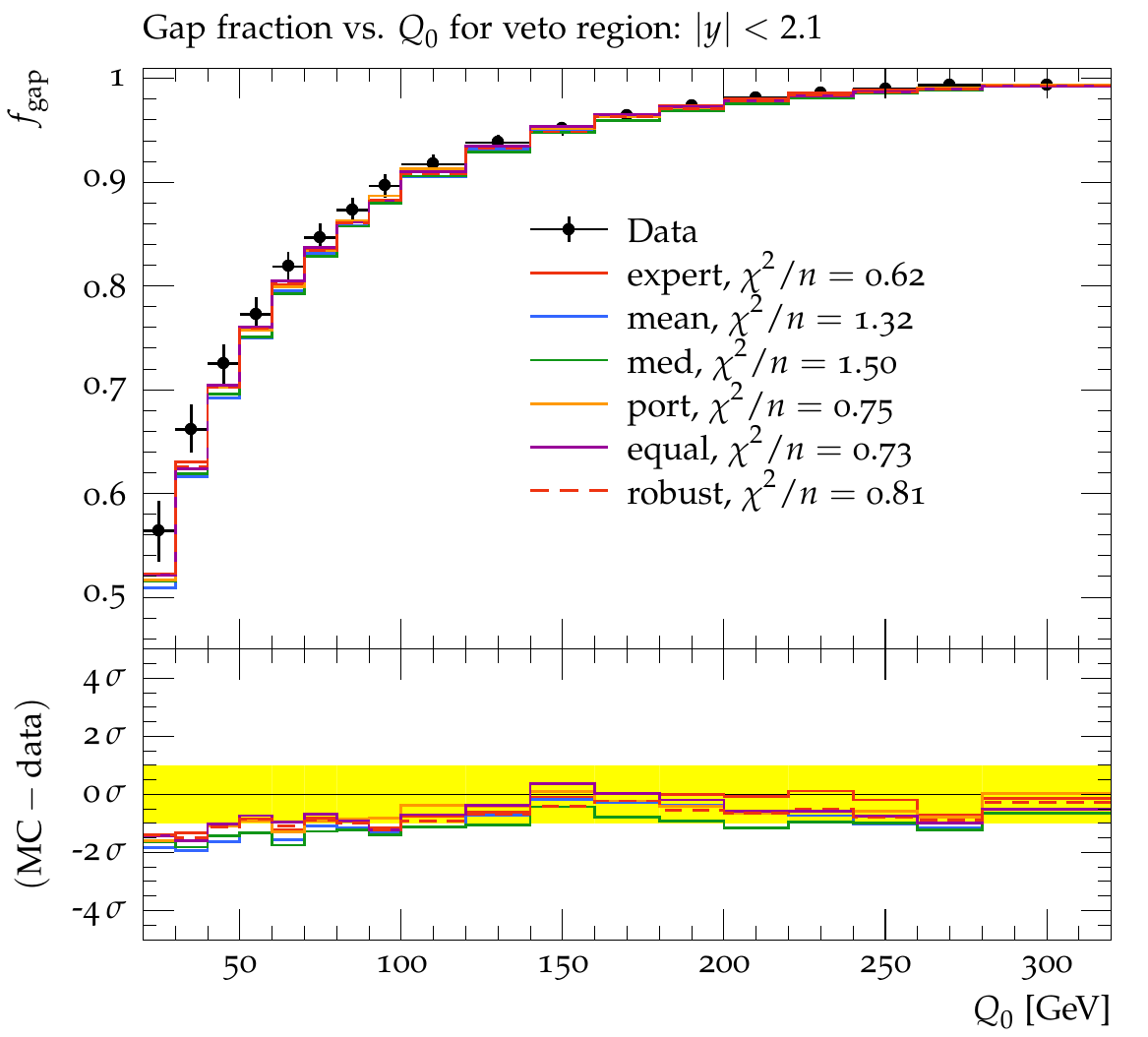}
         \caption{\it $t\bar{t}$ gap}
         \label{fig:worse-a14-ex2}
       \end{subfigure}
                \begin{subfigure}[b]{0.5\textwidth}
           \centering
           \includegraphics[width=\textwidth]{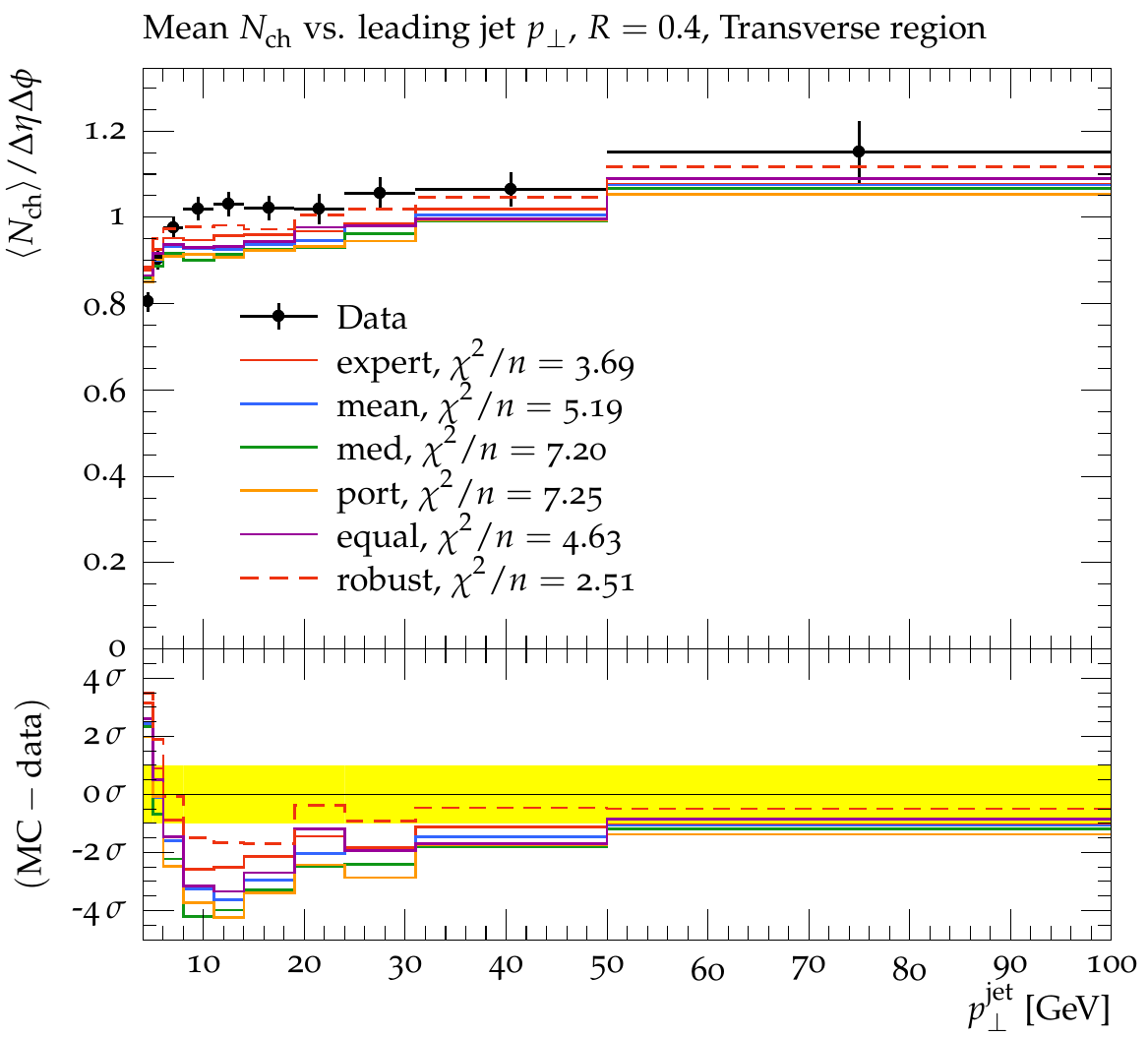}
         \caption{\it Track jet UE}
         \label{fig:worse-a14-ex3}
              \end{subfigure} 
 	\caption{Examples of A14 observables and their  $\chi^2_\calO$ values for which the automated tuning approach performs worse than the  expert's hand tuning. }
\end{figure}

Two categories stand out as being better described by the expert tune.
These are the {\it Multijets} and {\it $t\bar t$ gap} categories that were given a particularly large weight
in the expert tune.    
Some examples can be seen in
Figure~\ref{fig:worse-a14-ex1}-\ref{fig:worse-a14-ex3}.
It is no surprise that these categories are not
described well as the expert tune.
It is surprising  that the parameters sensitive to this observable, namely
{\tt TimeShower:alphaSvalue} and {\tt SpaceShower:alphaSvalue} are actually somewhat larger
than the expert tune values, see Table~\ref{tab:optimal_prms_31}.    Larger values for these parameters should mean forcing the prediction
to look {\it more} like a higher-order calculation.  Clearly, other data, such as {\it Dijet decorr} and $p_T^Z$ prefer
larger values for these parameters than the {\it Multijets} category alone.

It is worth noting that without expert input our automated methods do not emphasize these observables.
However, the results of the eigentunes are that the ``error bands'' for these parameters contain the values for
the A14 tune.   It is understood {\it a priori} that the physics models inside
\pythia behind the description of   {\it Multijets} and {\it $t\bar t$ gap}
need theory corrections.     In fact, even at the time of the A14 tune, methods existed to better describe
the {\it Multijets} category using \pythia, but it requires a process-dependent correction.    The explicit aim of the A14 tune was to
be process-agnostic and applied to physics predictions for which the corrections were not readily available or easily
applicable.   However, if the goal is to provide a tune that can be used in association with process-dependent corrections,
then those provided in this study are more appropriate.

\subsection{Results for \sherpa tuning}

\begin{figure}[!htbp]
  \begin{subfigure}[b]{0.5\textwidth}
    \centering
    \includegraphics[width=\textwidth]{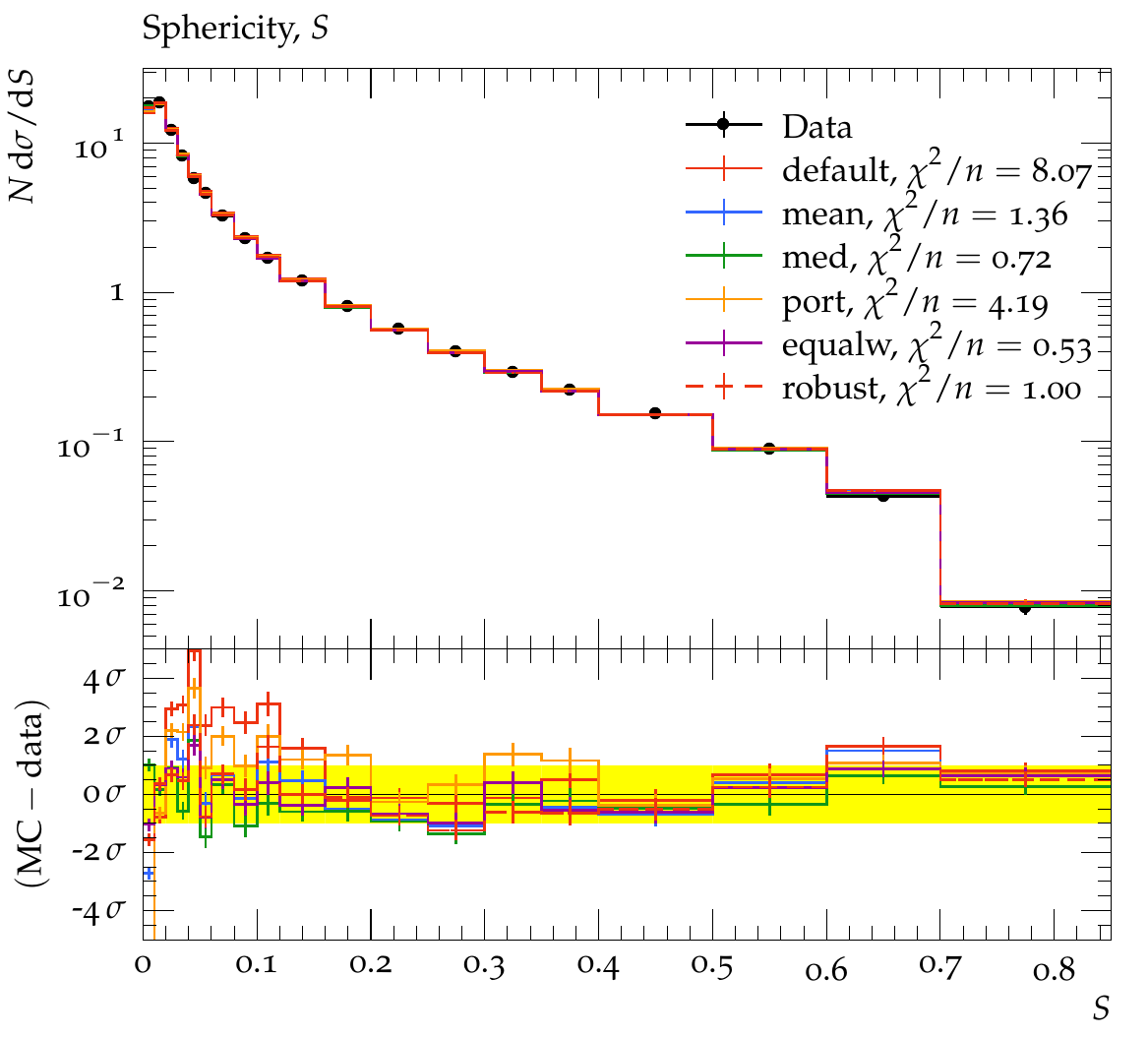}
    \caption{\it Jet shapes}
    \label{fig:sherpa-ex1}
  \end{subfigure}
  \begin{subfigure}[b]{0.5\textwidth}
    \centering
    \includegraphics[width=\textwidth]{\detokenize{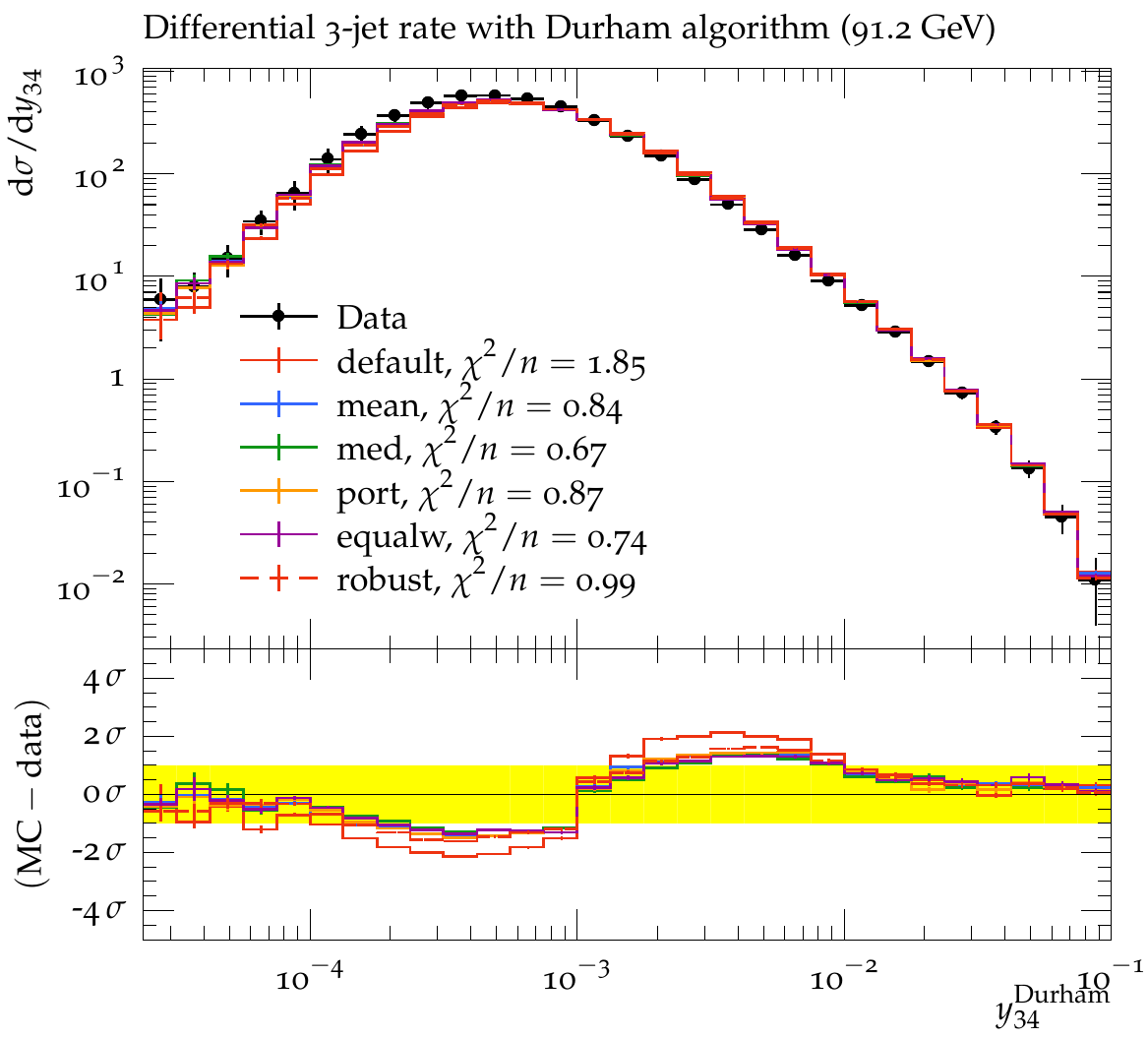}}
    \caption{\it Jet rates}
    \label{fig:sherpa-ex2}
  \end{subfigure}
  \begin{subfigure}[b]{0.5\textwidth}
    \centering
    \includegraphics[width=\textwidth]{\detokenize{plots/PDG_HADRON_MULTIPLICITIES_d03_x01_y03}}
    \caption{\it Particle count}
    \label{fig:sherpa-ex3}
  \end{subfigure}
   \begin{subfigure}[b]{0.5\textwidth}
    \centering
    \includegraphics[width=\textwidth]{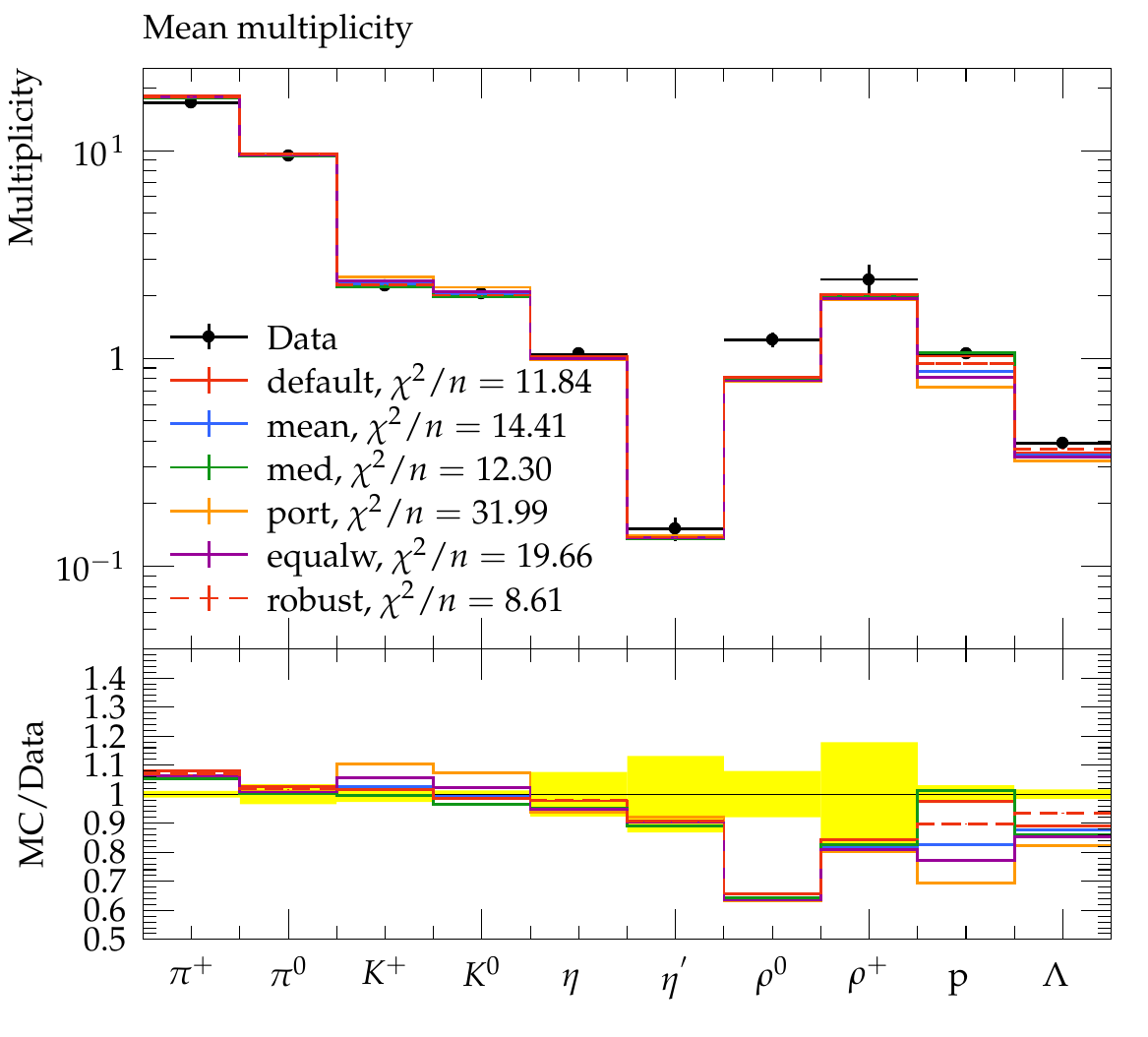}
    \caption{\it Several particle counts side-by-side}
    \label{fig:sherpa-ex4}
  \end{subfigure}
  \caption{Examples of histogram plots of the $\chi^2_\calO$
    values for the \sherpa tune.}
  \label{fig:sherpa}
\end{figure}

Some of the results of the \sherpa tuning are shown in
Figure~\ref{fig:sherpa}.   In general, all of the parameter selection
methods applied
here yield an improved global $\chi^2_\calO$ over the default values.
The parameters varied in this tuning exercise are all related to the
formation of physical particles.  This is a phenomenon that occurs at
a low-energy scale and cannot be described realistically (currently,
at least) from theory.    The model employed in \sherpa is a
cluster model that fissions colorless blobs of energy into particles
using a parameterized probability distribution.   Despite the fact that
hadronization
occurs at a low-energy scale, it has an impact on observables that are
used to test perturbative predictions at relatively high-energy
scales.
For these observables, it is impossible to entirely disentangle the
perturbative prediction from the non-perturbative hadronization model
prediction.   Figures~\ref{fig:sherpa-ex1}-\ref{fig:sherpa-ex2} show comparisons of
our tunes to the default, demonstrating a significant improvement in
most cases.  Figure~\ref{fig:sherpa-ex3} shows mixed results for the
production
of one particular species of particle.   Figure~\ref{fig:sherpa-ex4} is an example of an inclusive
observable that counts the number of particles produced without any
direct reference to their energy or position in the detector.

All of these results are for a certain precision of perturbation
theory.
There are both technical and mathematical reasons to truncate
perturbation theory in a certain order.    These calculations were
based on the lowest order perturbation theory with an improved parton
shower approximation to simulate additional perturbative effects.
The lowest order prediction produces 2 jets using exact perturbation
theory
and any additional jets using the parton shower approximation.
Figure~\ref{fig:sherpa-ex2} is an observable that counts the number of
3-jet events as a function of the jet definition.   While our results
are improved over the default, this indicates higher-order
perturbative
calculations might improve the description even more (e.g., 3 jets
calculated in exact perturbation theory and 4 or more jets from the
parton shower approximation).

Table~\ref{tab:optimal_prms_sherpa_RA} shows the parameters values for the
various tunes.    The simplest comparison is between the default
values and ``All-weights-equal.''  The all-weights-equal method yields the tune that
would result if only the data considered in this study were used.
One result is that several of the parameters take on the extremum of
the values considered here.    Without any additional direction to
choose
the range for our parameter scan, we chose $1/2$ of the default value
to define our sampling window.

One surprising result is that the parameter {\tt P\_QQ1\_by\_P\_QQ0},
which represents the ratio of spin-1 to spin-0 diquarks,
is driven to a value $>1$.  While there is no obvious reason that the
cluster model breaks down, spin-1 diquark production is usually
expected to be suppressed.   The fact that the parameter {\tt BARYON\_FRACTION}
is driven to its minimal value compensates for this large value.

While the type of large scale parameter tuning we have in mind here
can only be performed practically using surrogate models, the
fact that some tuned parameters are pushed to the boundaries suggests
another direction of algorithmic development.   In particular, we
would like our algorithm to have the capability to recognize a trust
region and update the surrogate model with dedicated simulations when
necessary.

\section{Conclusions}

In this paper, we propose several algorithms for automating weighting the importance of data used in the tuning process for Monte Carlo event generators.
We performed two studies.  The first used particle collider data and predictions are from the Large Hadron Collider (LHC) and had an {\it expert}
selection of analysis weights as a benchmark.   The second used data and predictions are from the Large Electron-Positron (LEP) Collider and had only the default parameter choices as a reference.   The algorithms considered included a bilevel optimization based on several scoring procedures and a single-level robust optimization.    We find that our automatic methods produce parameter tunes that are comparable to labor-intensive, by-hand tunes.
For the LHC tuning, filtering of hard-to-describe observables can lead to tunes of superior quality by identifying observables or subsets of observables that cannot be described by the event generator. For the LEP tuning, many of the tuned parameters were driven to the extremum of our sampling range, suggesting that the current models are missing some important physics.

First, the results show that the parameter values we found agree with and have the potential to improve the physicists' hand-tuned results. Second, since we automate the weight adjustment for the tune-relevant observables, physicists do not need to hand-tune the weights for observables anymore; we propose several methods for adjusting the  weights, so physicists are not involved in the subjective re-weighting anymore. Third, by filtering out and excluding observables and bins, we can save computational time during optimization and improve the parameter values. Fourth, we derived new metrics to easily compare different tunes, and it shows that our methods can perform better than the physicists' hand-tuning approach. 

For the \sherpa data, most of the optimal parameters are on the boundaries of the parameter space, indicating that we might need to change the size of the parameter domain to avoid model extrapolation. One possible solution to this problem is to build an outer loop for moving the center of the parameter search space and apply the trust region method. We leave this to future research.

In this work, we assume that each bin is completely independent of all the other bins. Taking into account bin correlations will be left for future research work.

\section*{Acknowledgements}
SM thanks Stefan Hoeche for discussions about our {\sc Sherpa} results.

\paragraph{Funding information}
This work was supported by the U.S. Department of Energy, Office of Science, Advanced Scientific Computing Research, under Contract DE-AC02-06CH11357.
Support for this work was provided through the Scientific Discovery through Advanced Computing (SciDAC) program funded by U.S. Department of Energy, Office of Science, Advanced Scientific Computing Research.
This work was also supported by
the U.S. Department of Energy through grant DE-FG02-05ER25694, and
by Fermi Research Alliance, LLC under Contract No. DE-AC02-07CH11359 with the U.S. Department of Energy, Office of Science, Office of High Energy Physics.
This work was supported in part by the U.S. Department of Energy, Office of Science, Office of Advanced Scientific
Computing Research and Office of Nuclear Physics, SciDAC program through the FASTMath Institute under Contract No. DE-AC02-05CH11231 at Lawrence Berkeley National Laboratory.
\bibliography{biblio,MPEC,ref}

\begin{thebibliography}{10}
\providecommand{\url}[1]{\texttt{#1}}
\providecommand{\urlprefix}{URL }
\expandafter\ifx\csname urlstyle\endcsname\relax
  \providecommand{\doi}[1]{doi:\discretionary{}{}{}#1}\else
  \providecommand{\doi}{doi:\discretionary{}{}{}\begingroup
  \urlstyle{rm}\Url}\fi
\providecommand{\eprint}[2][]{\url{#2}}

\bibitem{Buckley:2011ms}
A.~Buckley \emph{et~al.},
\newblock \emph{{General-purpose event generators for LHC physics}},
\newblock Phys. Rept. \textbf{504}, 145 (2011),
\newblock \doi{10.1016/j.physrep.2011.03.005},
\newblock \eprint{1101.2599}.

\bibitem{Skands:2014pea}
P.~Skands, S.~Carrazza and J.~Rojo,
\newblock \emph{{Tuning PYTHIA 8.1: the Monash 2013 Tune}},
\newblock Eur. Phys. J. \textbf{C74}(8), 3024 (2014),
\newblock \doi{10.1140/epjc/s10052-014-3024-y},
\newblock \eprint{1404.5630}.

\bibitem{ATL-PHYS-PUB-2014-021}
{ATLAS Collaboration},
\newblock \emph{{ATLAS PYTHIA 8 tunes to 7 TeV data}},
\newblock Tech. Rep. ATL-PHYS-PUB-2014-021, CERN, Geneva (2014).

\bibitem{Buckley:2009bj}
A.~Buckley, H.~Hoeth, H.~Lacker, H.~Schulz and J.~{von Seggern},
\newblock \emph{Systematic event generator tuning for the {LHC}},
\newblock The European Physical Journal C \textbf{65}, 331 (2010),
\newblock \doi{10.1140/epjc/s10052-009-1196-7}.

\bibitem{Bellm:2019owc}
J.~Bellm and L.~Gellersen,
\newblock \emph{{High dimensional parameter tuning for event generators}},
\newblock Eur. Phys. J. C \textbf{80}(1), 54 (2020),
\newblock \doi{10.1140/epjc/s10052-019-7579-5},
\newblock \eprint{1908.10811}.

\bibitem{Ilten:2016csi}
P.~Ilten, M.~Williams and Y.~Yang,
\newblock \emph{{Event generator tuning using Bayesian optimization}},
\newblock JINST \textbf{12}(04), P04028 (2017),
\newblock \doi{10.1088/1748-0221/12/04/P04028},
\newblock \eprint{1610.08328}.

\bibitem{ChenFlor:95}
Y.~Chen and M.~Florian,
\newblock \emph{The nonlinear bilevel programming problem: Formulations,
  regularity and optimality conditions},
\newblock Optimization \textbf{32}, 193 (1995).

\bibitem{MarZhu:96}
P.~Marcotte and D.~L. Zhu,
\newblock \emph{Exact and inexact penalty methods for the generalized bilevel
  programming problem},
\newblock Mathematical Programming \textbf{74}(2), 141 (1996).

\bibitem{YeZhuZhu:97}
J.~J. Ye, D.~L. Zhu and Q.~J. Zhu,
\newblock \emph{Exact penalization and necessary optimality conditions for
  generalized bilevel programming problems},
\newblock SIAM J. Optimization \textbf{7}(2), 481 (1997).

\bibitem{colson2007overview}
B.~Colson, P.~Marcotte and G.~Savard,
\newblock \emph{An overview of bilevel optimization},
\newblock Annals of operations research \textbf{153}(1), 235 (2007).

\bibitem{bard2013practical}
J.~F. Bard,
\newblock \emph{Practical bilevel optimization: algorithms and applications},
  vol.~30,
\newblock Springer Science \& Business Media (2013).

\bibitem{portfolio}
H.~Markowitz,
\newblock \emph{Portfolio selection},
\newblock The Journal of Finance \textbf{7}, 77 (1952),
\newblock \doi{10.2307/2975974}.

\bibitem{gneiting2007strictly}
T.~Gneiting and A.~Raftery,
\newblock \emph{Strictly proper scoring rules, prediction, and estimation},
\newblock Journal of the American Statistical Association \textbf{102}(477),
  359 (2007).

\bibitem{Powell1992}
M.~Powell,
\newblock \emph{Advances in Numerical Analysis, vol. 2: wavelets, subdivision
  algorithms and radial basis functions. Oxford University Press, Oxford, pp.
  105-210}, chap. The Theory of Radial Basis Function Approximation in 1990,
\newblock Oxford University Press, London (1992).

\bibitem{Muller2014}
J.~M\"uller and C.~Shoemaker,
\newblock \emph{Influence of ensemble surrogate models and sampling strategy on
  the solution quality of algorithms for computationally expensive black-box
  global optimization problems},
\newblock Journal of Global Optimization \textbf{60}, 123 (2014).

\bibitem{Mueller2017}
J.~M\"uller and J.~Woodbury,
\newblock \emph{{GOSAC}: global optimization with surrogate approximation of
  constraints},
\newblock Journal of Global Optimization \textbf{doi:10.1007/s10898-017-0496-y}
  (2017).

\bibitem{hawkins1980identification}
D.~M. Hawkins,
\newblock \emph{Identification of outliers}, vol.~11,
\newblock Springer (1980).

\bibitem{croarkin2006nist}
C.~Croarkin, P.~Tobias, J.~Filliben, B.~Hembree, W.~Guthrie \emph{et~al.},
\newblock \emph{{NIST/SEMATECH} e-handbook of statistical methods},
\newblock NIST/SEMATECH, July. Available online: http://www. itl. nist.
  gov/div898/handbook  (2006).

\bibitem{upton1996understanding}
G.~Upton and I.~Cook,
\newblock \emph{Understanding statistics},
\newblock Oxford University Press (1996).

\bibitem{rosner1983percentage}
B.~Rosner,
\newblock \emph{Percentage points for a generalized esd many-outlier
  procedure},
\newblock Technometrics \textbf{25}(2), 165 (1983).

\bibitem{grubbs1969procedures}
F.~E. Grubbs,
\newblock \emph{Procedures for detecting outlying observations in samples},
\newblock Technometrics \textbf{11}(1), 1 (1969).

\bibitem{stefansky1972rejecting}
W.~Stefansky,
\newblock \emph{Rejecting outliers in factorial designs},
\newblock Technometrics \textbf{14}(2), 469 (1972).

\bibitem{dixon1953processing}
W.~Dixon,
\newblock \emph{Processing data for outliers},
\newblock Biometrics \textbf{9}(1), 74 (1953).

\bibitem{thompson1985note}
R.~Thompson,
\newblock \emph{A note on restricted maximum likelihood estimation with an
  alternative outlier model},
\newblock Journal of the Royal Statistical Society: Series B (Methodological)
  \textbf{47}(1), 53 (1985).

\bibitem{dardis2004peirce}
C.~Dardis,
\newblock \emph{Peirce's criterion for the rejection of non-normal outliers;
  defining the range of applicability},
\newblock J Stat Softw \textbf{10}, 1 (2004).

\bibitem{tietjen1972some}
G.~L. Tietjen and R.~H. Moore,
\newblock \emph{Some {Grubbs}-type statistics for the detection of several
  outliers},
\newblock Technometrics \textbf{14}(3), 583 (1972).

\bibitem{ETZIONI1988416}
O.~Etzioni,
\newblock \emph{Hypothesis filtering: A practical approach to reliable
  learning},
\newblock In J.~Laird, ed., \emph{Machine Learning Proceedings 1988}, pp. 416
  -- 429. Morgan Kaufmann, San Francisco (CA),
\newblock ISBN 978-0-934613-64-4,
\newblock \doi{https://doi.org/10.1016/B978-0-934613-64-4.50047-5} (1988).

\bibitem{cochran1952}
W.~G. Cochran,
\newblock \emph{The $\chi^2$ test of goodness of fit},
\newblock Ann. Math. Statist. \textbf{23}(3), 315 (1952),
\newblock \doi{10.1214/aoms/1177729380}.

\bibitem{alex2019goodnessoffit}
A.~Shapiro, Y.~Xie and R.~Zhang,
\newblock \emph{Goodness-of-fit tests on manifolds} (2019),
  \eprint{1909.05229}.

\bibitem{10.1145/358234.381162}
J.~Bentley,
\newblock \emph{Programming pearls: Algorithm design techniques},
\newblock Commun. ACM \textbf{27}(9), 865–873 (1984),
\newblock \doi{10.1145/358234.381162}.

\bibitem{pronzato2013design}
L.~Pronzato and A.~P{\'a}zman,
\newblock \emph{Design of experiments in nonlinear models},
\newblock Lecture notes in statistics \textbf{212}, 1 (2013).

\bibitem{crestel2017optimal}
B.~Crestel, A.~Alexanderian, G.~Stadler and O.~Ghattas,
\newblock \emph{A-optimal encoding weights for nonlinear inverse problems, with
  application to the helmholtz inverse problem},
\newblock Inverse problems \textbf{33}(7), 074008 (2017).

\bibitem{KURAM2013159}
E.~Kuram, B.~Ozcelik, M.~Bayramoglu, E.~Demirbas and B.~T. Simsek,
\newblock \emph{Optimization of cutting fluids and cutting parameters during
  end milling by using d-optimal design of experiments},
\newblock Journal of Cleaner Production \textbf{42}, 159  (2013),
\newblock \doi{https://doi.org/10.1016/j.jclepro.2012.11.003}.

\bibitem{Sjostrand:2014zea}
T.~Sj\"ostrand, S.~Ask, J.~R. Christiansen, R.~Corke, N.~Desai, P.~Ilten,
  S.~Mrenna, S.~Prestel, C.~O. Rasmussen and P.~Z. Skands,
\newblock \emph{{An introduction to PYTHIA 8.2}},
\newblock Comput. Phys. Commun. \textbf{191}, 159 (2015),
\newblock \doi{10.1016/j.cpc.2015.01.024},
\newblock \eprint{1410.3012}.

\bibitem{Buckley:2010ar}
A.~Buckley, J.~Butterworth, D.~Grellscheid, H.~Hoeth, L.~Lonnblad, J.~Monk,
  H.~Schulz and F.~Siegert,
\newblock \emph{{{RIVET} user manual}},
\newblock Comput. Phys. Commun. \textbf{184}, 2803 (2013),
\newblock \doi{10.1016/j.cpc.2013.05.021},
\newblock \eprint{1003.0694}.

\bibitem{pope2008algorithms}
S.~B. Pope,
\newblock \emph{Algorithms for ellipsoids},
\newblock Cornell University Report No. FDA pp. 08--01 (2008).

\bibitem{Bothmann:2019yzt}
E.~Bothmann \emph{et~al.},
\newblock \emph{{Event Generation with {Sherpa} 2.2}},
\newblock SciPost Phys. \textbf{7}(3), 034 (2019),
\newblock \doi{10.21468/SciPostPhys.7.3.034},
\newblock \eprint{1905.09127}.

\bibitem{Pfeifenschneider:1999rz}
P.~Pfeifenschneider \emph{et~al.},
\newblock \emph{{QCD analyses and determinations of $\alpha_s$ in $e^+ e^-$
  annihilation at energies between 35-GeV and 189-GeV}},
\newblock Eur. Phys. J. C \textbf{17}, 19 (2000),
\newblock \doi{10.1007/s100520000432},
\newblock \eprint{hep-ex/0001055}.

\bibitem{Abreu:1996na}
P.~Abreu \emph{et~al.},
\newblock \emph{{Tuning and test of fragmentation models based on identified
  particles and precision event shape data}},
\newblock Z. Phys. C \textbf{73}, 11 (1996),
\newblock \doi{10.1007/s002880050295}.

\bibitem{Abe:2002iq}
K.~Abe \emph{et~al.},
\newblock \emph{{Measurement of the b quark fragmentation function in $Z^0$
  decays}},
\newblock Phys. Rev. D \textbf{65}, 092006 (2002),
\newblock \doi{10.1103/PhysRevD.66.079905},
\newblock [Erratum: Phys.Rev.D 66, 079905 (2002)],
\newblock \eprint{hep-ex/0202031}.

\bibitem{Amsler:2008zzb}
C.~Amsler \emph{et~al.},
\newblock \emph{{Review of Particle Physics}},
\newblock Phys. Lett. B \textbf{667}, 1 (2008),
\newblock \doi{10.1016/j.physletb.2008.07.018}.

\bibitem{Pumplin:2001ct}
J.~Pumplin, D.~Stump, R.~Brock, D.~Casey, J.~Huston, J.~Kalk, H.~L. Lai and
  W.~K. Tung,
\newblock \emph{{Uncertainties of predictions from parton distribution
  functions. 2. The Hessian method}},
\newblock Phys. Rev. D \textbf{65}, 014013 (2001),
\newblock \doi{10.1103/PhysRevD.65.014013},
\newblock \eprint{hep-ph/0101032}.

\bibitem{Booker1999}
A.~Booker, J.~{Dennis Jr}, P.~Frank, D.~Serafini, V.~Torczon and M.~Trosset,
\newblock \emph{A rigorous framework for optimization of expensive functions by
  surrogates},
\newblock Structural Multidisciplinary Optimization \textbf{17}, 1 (1999).

\bibitem{Powell1999}
M.~Powell,
\newblock \emph{Recent Research at Cambridge on Radial Basis Functions},
\newblock New Developments in Approximation Theory, pp. 215-232. Birkh\"auser,
  Basel (1999).

\bibitem{kotz2004continuous}
S.~Kotz, N.~Balakrishnan and N.~L. Johnson,
\newblock \emph{Continuous multivariate distributions, Volume 1: Models and
  applications},
\newblock John Wiley \& Sons (2004).

\bibitem{blei2003latent}
D.~M. Blei, A.~Y. Ng and M.~I. Jordan,
\newblock \emph{Latent {Dirichlet} allocation},
\newblock the Journal of machine Learning research \textbf{3}, 993 (2003).

\bibitem{ng2011dirichlet}
K.~W. Ng, G.-L. Tian and M.-L. Tang,
\newblock \emph{Dirichlet and related distributions: Theory, methods and
  applications}  (2011).

\bibitem{Regis2007b}
R.~Regis and C.~Shoemaker,
\newblock \emph{{A} stochastic radial basis function method for the global
  optimization of expensive functions},
\newblock INFORMS Journal on Computing \textbf{19}, 497 (2007).

\end{thebibliography}

\clearpage

\section{Online Supplement}
Online supplement for ``BROOD:   Bilevel  and Robust Optimization and
      Outlier Detection for Efficient Tuning of High-Energy Physics Event Generators''.

\subsection{Solving the outer problem with derivative-free surrogate optimization}
\label{sec:surrogate}
Solving the inner optimization problem can become computationally demanding as it depends on the number of observables involved,   the number of bins per observable (and therefore the number of  parameters), and the starting guess (and therefore the number of iterations needed). Thus, the goal is to  determine the  optimal weights $\w^*$ within as few iterations of the outer loop as possible since this number determines how often we have to solve the inner optimization problem. We  do not have a full analytic expression of $g(\w, \widehat{\p}_{\w})$ (black box) since computing this  value involves solving the inner optimization problem. Thus, also derivatives of $g(\w, \widehat{\p}_{\w})$ are not available. A widely used approach for optimizing computationally expensive black-box functions is to use computationally cheap approximations (surrogates, metamodels) of the expensive function and to use the approximation throughout the optimization to make iterative sampling decisions~\cite{Booker1999}. Here, we approximate  $g(\w, \widehat{\p}_{\w})$ with a radial basis function (RBF) \cite{Powell1999}, although in general any approximation model could be used. An RBF interpolant is defined as follows:
\begin{equation}
s(\w) = \sum_{i = 1}^n \gamma_i \phi(\|\w-\w_i\|_2)+q(\w),\label{eq:rbf}
\end{equation}
where $s:\mathbb{R}^{|\calS_\calO|}\mapsto\mathbb{R}$, $\w_i, \ i =1,\ldots,n$,  are the weight vectors for which we have already evaluated the objective function of the outer optimization problem,  $\gamma_i$ are parameters that must be determined, $\phi(\cdot)$ is the radial basis function (here, we choose the cubic, $\phi(r)=r^3$, but other options are possible), $\|\cdot\|_2$ denotes the Euclidean norm, and $q(\cdot)$ is a polynomial tail whose order depends on the choice of  $\phi$. When using the cubic RBF, the polynomial tail must be at least linear ($q(\w) =\beta_0 +\boldsymbol\beta^{\top}\w$) in order to uniquely determine the RBF parameters ($\gamma_i, i = 1,\ldots,n, \beta_0, \boldsymbol\beta=[\beta_1, \ldots, \beta_{|\calS_\calO|}]^{\top})$. The RBF interpolant $s(\w)$ then predicts the value of the objective function at the point $\w$. It is interpolating, and thus the prediction at an already evaluated point $\w_i$ will agree with the observed function value. Using the RBF, we thus have  $g(\w, \widehat{\p}_{\w}) = s(\w)+e(\w)$, where $e(\w)$ denotes the difference between the RBF and the true function value and it is 0 at already evaluated vectors $\w_i$. The values of the RBF parameters are determined by solving a linear system of equations: 

\begin{equation}
\begin{bmatrix}
\boldsymbol{\Phi} & \mathbf{W} \\
\mathbf{W}^{\top} & \mathbf{0}
\end{bmatrix}
\begin{bmatrix}
\boldsymbol{\gamma}\\
\boldsymbol{\beta}^\prime
\end{bmatrix}=
\begin{bmatrix}
\mathbf{G}\\
\mathbf{0}
\end{bmatrix},\label{eq: RBFsystem} 
\end{equation}
where the elements of $\boldsymbol{\Phi}$ are  $\Phi_{{\iota\nu}}=\phi(\|\w_{\iota}-\w_{\nu}\|_2)$, $\iota,\nu=1\ldots n$, $\mathbf{0}$ is a matrix with all entries  0 of appropriate dimension, and
\begin{equation}
\mathbf{W} = \begin{bmatrix}
\w_{1}^{\top} & 1\\
\vdots &\vdots\\
\w_{{n}}^{\top} & 1\\
\end{bmatrix}
\quad
\boldsymbol{\gamma}=
\begin{bmatrix}
\gamma_1\\
\gamma_2\\
\vdots\\
\gamma_{n}
\end{bmatrix}
\quad
\boldsymbol{\beta}^\prime=
\begin{bmatrix}
\beta_1\\
\beta_2\\
\vdots\\
\beta_{|\calS_\calO|}\\
\beta_0
\end{bmatrix}
\quad
\mathbf{G}=\begin{bmatrix}
g(\w_1,\widehat\p_{\w_1})\\
 g(\w_2,\widehat\p_{\w_2})\\
 \vdots\\
  g(\w_n,\widehat\p_{\w_n})
\end{bmatrix}.
\label{eq: RBFmatrices}
\end{equation}

The linear system in \equaref{eq: RBFsystem} has a solution if and only if rank$(\mathbf{W}) = |\calS_\calO|+1$ \cite{Powell1992}. During the optimization, we use the RBF prediction at unsampled points to determine a new vector $\w$ for which we solve the inner optimization problem. It is important that at this step only weights that sum up to 1 are chosen. The steps of the iterative sampling algorithm are summarized in~\algref{alg:rbfmethod}.

\begin{algorithm2e}
\caption{Derivative-free optimization of the outer equality-constrained optimization problem}
\label{alg:rbfmethod}
\begin{algorithmic}[1]
	\REQUIRE Number of initial experimental design points $n_0$; the maximum number of evaluations $n_{\tmax}$ 
	\ENSURE The best weight vector $\w^*$ and corresponding $\widehat\p^*_{\w^*}$
	\STATE Create an initial experimental design with $n_0$ points; ensure that \equaref{eq:sum1_bilevel} is satisfied for all points;\label{alg:init}
	\STATE Compute the value of the outer optimization objective function at all points in the initial design;
	\STATE Fit an RBF model to the sample data pairs $\{(\w_i, g(\w_i,\widehat\p_{\w_i}))\}_{i = 1}^{n_0}$
	\STATE Set $n=n_0$
	\WHILE {$n<n_\tmax$}\label{alg:forstart}
	\STATE Use the RBF to determine a new point $\w_\tnew$ and ensure that \equaref{eq:sum1_bilevel} is satisfied;
	\STATE Solve the inner optimization problem for $\w_\tnew$ and obtain $\widehat\p_{\w_\tnew}$;
	\STATE Compute the value of the outer optimization objective function for $(\w_\tnew, \widehat\p_{\w_\tnew})$;
	\STATE Update the RBF model with the new data;
	\STATE $n\leftarrow n+1$;
	\ENDWHILE\label{alg:forend}
	\RETURN the best parameter values $(\w^*, \widehat\p^*_{\w^*})$;
\end{algorithmic}
\end{algorithm2e}

The inputs that must be supplied to the algorithm are the number of points $n_0$ to be used in the initial experimental design and the maximum number $n_{\tmax}$ of outer objective function evaluations (i.e., the number of inner optimizations) one is willing to allow. The number $n_0$ should in our case be at least $|\calS_\calO|+1$, since this is the minimum number of points we need to fit the RBF model. $n_{\tmax}$ should depend on how long the inner optimization takes and the time budget of the user. 

When creating the initial experimental design in Step \ref{alg:init}, we have to ensure that the constraint~(\ref{eq:sum1_bilevel}) is satisfied. Also, we have the condition that the weights lie in $[0,1]$ and are uniform in their support. This means that the weights follow the Dirichlet distribution, i.e., the set of points are uniformly distributed over the open standard $(|\calS_\calO| - 1)$-simplex. To achieve this, we generate an initial design where all weights are drawn from the symmetric Dirichlet distribution, $Dir(\alpha_1=\alpha_2=\ldots=\alpha_{|\calS_\calO|}=1)$ \cite{kotz2004continuous,blei2003latent,ng2011dirichlet}.

We evaluate the outer objective function at these points, i.e., we solve the inner optimization problem at each point and we obtain $\mathbf{G}$ in \equaref{eq: RBFmatrices}. 
With the sum-one-scaled initial experimental design, however, the rank of the matrix $\boldsymbol{W}$ is now only $|\calS_\calO|$ (and not the required $|\calS_\calO|+1$). Thus, we solve the problem as one of dimension $|\calS_\calO|-1$, i.e., for  fitting the RBF model, we only use the first $|\calS_\calO|-1$ values of each sample point (the ``reduced'' sample points). Thus, we use
\begin{equation}
\boldsymbol{W}=\begin{bmatrix}
w_{1,1} & w_{1,2} & \ldots &w_{1,|\calS_\calO|-1}  & 1\\
w_{2,1} & w_{2,2} & \ldots &w_{2,|\calS_\calO|-1}  & 1\\
\vdots &\vdots &\vdots &\vdots &\vdots\\
w_{n,1} & w_{n,2} & \ldots &w_{n,|\calS_\calO|-1}  & 1\\
\end{bmatrix}
\end{equation}
and the coefficient vector for the polynomial tail thus becomes $[\beta_1, \ldots, \beta_{|\calS_\calO|-1}, \beta_0]^{\top}$. The vector $\boldsymbol{\gamma}$ and the matrix $\boldsymbol{G}$ do not change. The elements of $\boldsymbol\Phi$ are computed from the $(|\calS_\calO|-1)$-dimensional sample  vectors. Note, however, that when we evaluate the objective function in \equaref{eq:outer_gen}, we always evaluate it for the full-dimensional vectors, as we can simply compute $w_{j, |\calS_\calO|}=1-\sum_{i = 1}^{|\calS_\calO|-1}w_i$ for each $j = 1,\ldots, n$.

In the iterative sampling procedure (Steps~\ref{alg:forstart}-\ref{alg:forend}), we use the RBF model to determine a new vector $\w_{\tnew}$ at which we will do the next evaluation of   \equaref{eq:outer_gen}. Since we do not know whether the objective function is multimodal, we have to balance local and global search steps, i.e., we have to balance our sample point selection such that we select points with low predicted function values but also points that are far away from already evaluated points. 
Moreover, the new sample point must satisfy  \equaref{eq:sum1_bilevel}. 
In order to do so, %
we generate a large set of candidate points from the Dirichlet distribution. %
We   use the RBF to predict the function values at the candidate points. Since the RBF is defined over the $(|\calS_\calO|-1)$-dimensional space, we use only the first $|\calS_\calO|-1$ parameter values of each candidate point.  We denote the  $(|\calS_\calO|-1)$-dimensional candidate points by $\x_1, \ldots, \x_{N_{\tcand}}$, where we choose $N_{\tcand}$ large (for example, $500|\calS_\calO|$). For each candidate point, we use the RBF to predict its function value using (\ref{eq:rbf}) and we obtain $s(\x_k), k=1,\ldots, N_{\tcand}$. We scale these values to [0,1] according to 
\begin{equation}
V_s(\x_k) = \frac{s(\x_k)-s_{\tmin}}{s_{\tmax}-s_{\tmin}}, k = 1, \ldots, N_{\tcand},
\end{equation}
 where 
 \begin{equation}
 s_{\tmin} = \min\{s(\x_k), k = 1,\ldots, N_{\tcand}\} \text{ and } s_{\tmax} = \max\{s(\x_k), k = 1,\ldots, N_{\tcand}\}.
 \end{equation}
 We also compute the distances $d(\x_k, S)$ of each candidate point to the set of already evaluated points $S$ (in the $(|\calS_\calO|-1)$-dimensional Euclidean space),  and we scale these distances to [0,1] according to
 \begin{equation}
V_d(\x_k) = \frac{d_{\tmax}-d(\x_k)}{d_{\tmax}-d_{\tmin}}, k = 1, \ldots, N_{\tcand},
\end{equation}
 where 
 \begin{equation}
 d_{\tmin} = \min\{d(\x_k, S), k = 1,\ldots, N_{\tcand}\} \text{ and } d_{\tmax} = \max\{d(\x_k, S), k = 1,\ldots, N_{\tcand}\}.
 \end{equation}
The ideal new sample point $\w_{\tnew}$ will have a large distance to the set of already evaluated points $S$ and a low predicted objective function value. Using the two criteria defined above, we compute a weighted sum of both (following \cite{Regis2007b})
\begin{equation}
V(\x_k) = \nu V_s(\x_k) +(1-\nu)V_d(\x_k), k = 1, \ldots, N_{\tcand},
\end{equation}
where $\nu \in[0,1]$ is a parameter that determines how much emphasis we put on either criterion. If $\nu$ is large, it means we put most emphasis on $V_s$, and we favor candidate points that have low predicted objective function values. This also means that the search is more local as low function values are usually predicted around the best point found so far. If $\nu$ is small, we put more emphasis on $V_d$ and we favor points that are far away from the set of already evaluated points, and thus the search is more global. By varying the weights $\nu$ between  different values in [0,1], we  can achieve a repeated transition between local and global search, and therefore we can avoid becoming stuck in a local optimum. The candidate point with the lowest $V$ value will become the new sample point $\w_{\tnew}$. We evaluate the objective function (inner optimization) at the new point (augmented with the missing parameter value), and given the new data, we update the RBF model. The algorithm iterates until the maximum number of function evaluations $n_{\tmax}$ has been reached.

\subsection{Polynomial-time algorithm for filtering bins by hypothesis testing}
\label{sec:polyalghyptest}

In this section, we describe the polynomial-time algorithm to solve the problem of finding the largest contiguous subset of bins $\calB\subset\calO$ to be kept for tuning, i.e., finding the largest contiguous subset of bins $\calB\subset\calO$ such that $\chi^2_\calB \le \chi^2_{c,\calB}$, where $\chi^2_{c,\calB}$ is the critical value for bins in $\calB$. This algorithm is described in Algorithm~\ref{A:goodbins} and it is based on the maximum subarray problem~\cite{10.1145/358234.381162}.  

In this algorithm, we first find the critical value for each bin in line 1 as described in Section~\ref{sec:binfilter}. The degree of the freedom is given by $\rho_\calB =  |\calB| - d$ and since $\rho_\calB$ cannot be negative, the critical values for only the bin index $b >d$ is calculated in line 1. Then the $\chisq$ test statistic is computed for each bin in $\calO$ in lines 2-3.
Then, while iterating through the bins in $\calO$, in lines 6, we check whether the current bin $b$ can be added to $\calB$ and if so, we update the counters and add the current bin $b$ to the end of $\calB$ (via $e$) in lines 7-10. If the current bin $b$ cannot be added to $\calB$, then in lines 12-13 we shift the start $s$ of $\calB$ (through $\tau$)  such that the start is now at the bin index where the condition in line 6 could be satisfied in future iterations. 
Finally, in lines 14-19, we perform a sanity check to make sure that $\calB$ contains the set of bins that yield the lowest $\chisq_\calB$ test statistic.

\begin{algorithm2e}[htb!]
	\DontPrintSemicolon \setcounter{AlgoLine}{0}
	\caption{Algorithm to find bins $\calB$  in observable $\calO$ to keep for tuning}\label{A:goodbins}
	\SetKwInOut{Input}{Input}
	\SetKwInOut{Output}{Output}
	\Input{$f_b, \calR_b, \Delta f_b, \Delta\calR_b, \forall b \in \calO;$ significance level $\alpha$}
	\Output{start index $s$ and end index $e$ of bins, i.e., $\calB = \{s,\dots,e\}$ to keep in $\mathcal{O}$}
	Calculate the critical values for each bin:
	\[
	k_b = 
	\begin{cases}
	\chi^2_{c,b} = f(\rho_b,\alpha), & \text{if } b > d  \\
	\infty,& \text{otherwise}  
	\end{cases},
	\quad \forall b\in \{1,2,\dots,\abscalO\},\p\in\Omega\subset\mathbb{R}^d
	\]
	\;
	Find  $\p^*$ by minimizing $\chisq_\calO$ in Eq.~\eqref{eq:obschi2}  \;
	Calculate the test statistic values for each bin:
	$\chisq_b(\p^*) =  \frac{(f_b(\p^*)-\calR_b)^2}{\Delta f_b(\p^*)^2+\Delta\calR_b^2}, \quad \forall b\in \{1,2,\dots,\abscalO\}$\;
	Initialize $\Sigma \gets 0, \widehat{b} \gets 0, s \gets 0, e \gets 0, \tau \gets 0$\;
	\For{b $\in  \{1,2,\dots,\abscalO\}$ }{
		\uIf{$\Sigma + \chi^2_b \le k_{\widehat{b}+ 1}$}{
			$\Sigma \gets \Sigma+ \chi^2_b$\;
			$\widehat{b} \gets \widehat{b}+1$\;
			$s \gets \tau$\;
			$e \gets b$\;
		}
		\uElseIf{$\Sigma \neq 0$}{
			$\Sigma \gets \Sigma - \chi^2_{b - \widehat{b}} + \chi^2_{b}$\;
			$\tau \gets b-\widehat{b}+1$\;
		}
	}
	\uIf{$s > 0$ and $\chi^2_{s-1} < \chi^2_{e}$}{
		$e \gets e - 1$\;
		$s \gets s - 1$\;
	}
	\uElseIf{$e < \abscalO $ and $\chi^2_{s} > \chi^2_{e+1}$}{
		$e \gets e + 1$\;
		$s \gets s + 1$\;
	}
	\Return $\calB = \{s,\dots,e\}$.

\end{algorithm2e}
 
\subsection{A14 and \sherpa physics parameters}
\label{sec:physicsparam}
The A14 tunable physics parameters, their definitions and tuning ranges are shown in Table~\ref{tab:prm-ranges}. The \sherpa parameters, their definitions and  tuning ranges are shown in  Table~\ref{tab:prm-ranges-sherpa}.

\begin{table}[ht!]
	\centering
	\small
	\caption{A14 physics parameters, their definitions and tuning ranges (min, max).}\label{tab:prm-ranges}
	\begin{adjustbox}{width=\textwidth}
		\begin{tabular}{llll}
			\hline
			Parameters                          & Definition & min &  max \\ \hline
\tt			SigmaProcess:alphaSvalue            & The $\alpha_S$ value at scale $Q^2 = M^2_Z$           & 0.12      & 0.15      \\
\tt			BeamRemnants:primordialKThard       & Hard interaction primordial $k_\perp$           & 1.5       & 2.0         \\
\tt			SpaceShower:pT0Ref                  & ISR $p_T$ cutoff           & 0.75      & 2.5       \\
\tt			SpaceShower:pTmaxFudge              & Mult. factor on max ISR evolution scale           & 0.5       & 1.5       \\
\tt			SpaceShower:pTdampFudge             & Factorization/renorm scale damping           & 1.0         & 1.5       \\
\tt			SpaceShower:alphaSvalue             & ISR $\alpha_S$           & 0.10       & 0.15      \\
\tt			TimeShower:alphaSvalue              & FSR $\alpha_S$            & 0.10       & 0.15      \\
\tt			MultipartonInteractions:pT0Ref      &  MPI $p_T$ cutoff          & 1.5       & 3.0         \\
\tt			MultipartonInteractions:alphaSvalue & MPI $\alpha_S$            & 0.10       & 0.15      \\
\tt			BeamRemnants:reconnectRange         & CR strength           & 1.0         & 10.0        \\ \hline
		\end{tabular}
	\end{adjustbox}
\end{table}

\begin{table}[!ht]
	\centering
	\small
	\caption{\sherpa physics parameters, their definitions and tuning ranges (min, max).}\label{tab:prm-ranges-sherpa}
		\begin{tabular}{llll}
			\hline
			Parameters                          & Definition & min &  max \\ \hline
\tt		KT\_0             & generic parameter for non-perturbative transverse momentum    & 0.5  & 1.5   \\
\tt		ALPHA\_G          & gluon fragmentation    & 0.62  & 1.88  \\
\tt		ALPHA\_L          & light quark fragmentation $z$ power    & 1.25  & 3.75  \\
\tt		BETA\_L           & light quark fragmentation  $1-z$ power  & 0.05  & 0.15  \\
\tt		GAMMA\_L          & light quark fragmentation  $\exp$ power  & 0.25  & 0.75  \\
\tt		ALPHA\_H          & heavy quark fragmentation $z$ power   & 1.25  & 3.75  \\
\tt		BETA\_H           & heavy quark fragmentation $1-z$ power   & 0.375  & 1.125  \\
\tt		GAMMA\_H          & heavy quark fragmentation  $\exp$ power   & 0.05  & 0.15  \\
\tt		STRANGE\_FRACTION & suppression of $s$ quarks %
		& 0.25  & 0.75  \\
\tt		BARYON\_FRACTION  & suppression of baryons & 0.09  & 0.27  \\
\tt		P\_QS\_by\_P\_QQ\_norm & fraction of di-quarks with one strange quark & 0.24  & 0.72  \\
\tt		P\_SS\_by\_P\_QQ\_norm & fraction of di-quarks with two strange quarks  & 0.01  & 0.03  \\
\tt		P\_QQ1\_by\_P\_QQ0     & fraction of di-quarks with spin-1 to spin-0  & 0.5  & 1.5 \\ \hline
		\end{tabular}
\end{table}

\subsection{Selection of the best hyperparameter in robust optimization}
\label{sec:ROmuselection}

\begin{figure}[ht!]
	\centering
	\includegraphics[width=0.8\textwidth]{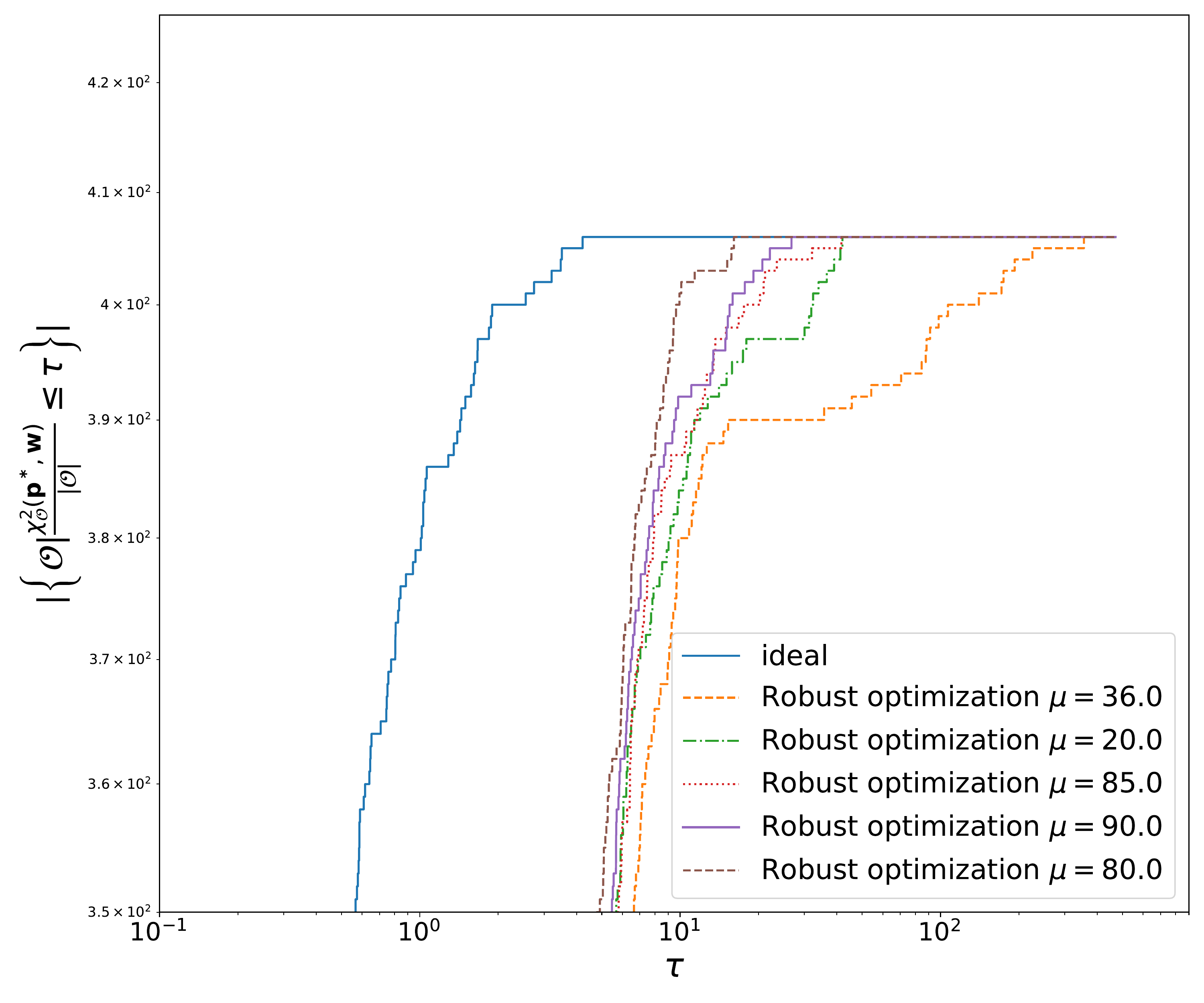}
	\caption{Ideal cumulative density curve and the cumulative density curves of robust optimization runs with selected hyperparameter values $\mu$ for the A14 dataset.}
	\label{fig:romuperformance}
\end{figure}

In order to find the best value for  $\mu$ in the  robust optimization,  we first build for each run (each $\mu$)  a cumulative density curve of the number of observables for which  $\frac{\chisq_\calO (\p^*,\w)}{|\calO|} \le \tau$, where $\p^*$ is the optimal parameter obtained from the robust optimization run, $\w=\bm{1}$, $\tau \in \mathbb{R}^+$ and $\calO \in \calS_\calO$. 
Then, we construct the ``ideal'' cumulative density curve, for which  $\p^*$ in $\frac{\chisq_\calO (\p^*,\w)}{|\calO|} \le \tau$ is obtained by optimizing   for each  observable $\calO$ separately. %
An example plot showing the cumulative density curve from the ideal case to some of the robust optimization runs is shown in Figure~\ref{fig:romuperformance}.

Then, the area between the cumulative density curve for each robust optimization run and the ideal cumulative density curve is computed. For the A14 dataset and all runs completed for robust optimization, the area between the curve is given in Table~\ref{tab:A14ROmuABC} (smaller values are better). %
Finally, for completeness, the best values of $\mu$ found for both the A14 and \sherpa datasets are given in Table~\ref{tab:allbestmu}.

\begin{table}[htp!]
	\centering
	\caption{Area between ideal cumulative density curve and the cumulative density curves of the robust optimization runs for various hyperparameters $\mu$ for the A14 full dataset (smaller values are better). The data are organized in ascending order of the area between the curves.} \label{tab:A14ROmuABC}
	\begin{tabular}{|c|c|c|c|c|c|c|c|c|c|c|c|}\hline
		rank &	$\mathbf{\mu}$ & Area & rank &	$\mathbf{\mu}$ & Area & rank &	$\mathbf{\mu}$ & Area & rank &	$\mathbf{\mu}$ & Area\\\hline
		
		1&	80 & 7.51e+02 &26&	79 & 1.05e+03 & 51&	24 & 1.32e+03 & 76&		38 & 1.66e+03\\\hline
		2&	78 & 7.93e+02 &27&	81 & 1.12e+03 & 52&	35 & 1.32e+03 & 77&		29 & 1.72e+03\\\hline
		3&	76 & 7.95e+02 &28&	71 & 1.13e+03 & 53&	93 & 1.35e+03 & 78&		51 & 1.72e+03\\\hline
		4&	77 & 8.15e+02 &29&	95 & 1.14e+03 & 54&	89 & 1.39e+03 & 79&     49 & 1.73e+03\\\hline
		5&	73 & 8.53e+02 &30&	10 & 1.15e+03 & 55&	45 & 1.40e+03 &  80&    57 & 1.73e+03\\\hline
		6&	70 & 8.91e+02 &31&	11 & 1.17e+03 & 56&	42 & 1.41e+03 & 81&    50 & 1.74e+03\\\hline
		7&	90 & 8.96e+02 &32&	12 & 1.17e+03 & 57&	41 & 1.43e+03 & 82&  43 & 1.74e+03 \\\hline
		8&	26 & 9.09e+02 &33&	18 & 1.18e+03 & 58&	68 & 1.43e+03 & 83& 44 & 1.76e+03 \\\hline
		9&	88 & 9.11e+02 &34&	3& 1.19e+03 & 59&	67 & 1.44e+03 & 84& 55 & 1.79e+03 \\\hline
		10&	74 & 9.14e+02 &35&	21 & 1.19e+03 & 60&	46 & 1.44e+03 & 85&  47 & 1.87e+03\\\hline
		11&	86 & 9.42e+02 &36&	20 & 1.19e+03 & 61&	39 & 1.46e+03 & 86&  60 & 1.93e+03\\\hline
		12&	72 & 9.43e+02 &37&	16 & 1.19e+03 & 62&	30 & 1.48e+03 & 87&  37 & 1.94e+03\\\hline
		13&	27 & 9.44e+02 &38&	69 & 1.20e+03 & 63&	63 & 1.48e+03 & 88&  59 & 1.95e+03  \\\hline
		14&	83 & 9.47e+02 &39&	22 & 1.20e+03 & 64&	40 & 1.51e+03 & 89&  33 & 1.96e+03 \\\hline
		15&	75 & 9.53e+02 &40&	23 & 1.21e+03 & 65&	64 & 1.52e+03 & 90&  54 & 1.97e+03 \\\hline
		16&	87 & 9.59e+02 &41&	19 & 1.21e+03 & 66&	28 & 1.55e+03 & 91&  58 & 1.99e+03\\\hline
		17&	8 & 9.61e+02 &42&	13 & 1.22e+03 & 67&	61 & 1.56e+03 & 92&  53 & 2.08e+03\\\hline
		18&	82 & 9.72e+02 &43&	25 & 1.23e+03 & 68&	98 & 1.56e+03 & 93&  94 & 2.13e+03\\\hline
		19&	2 & 9.77e+02 &44&	97 & 1.24e+03 & 69&	62 & 1.58e+03 & 94&  56 & 2.14e+03\\\hline
		20&	84 & 9.80e+02 &45&	7 & 1.25e+03 & 70&	66 & 1.58e+03 & 95&  52 & 2.23e+03\\\hline
		21&	5 & 9.90e+02 &46&	15 & 1.27e+03 & 71&	32 & 1.58e+03 & 96&  99 & 2.74e+03\\\hline
		22&	85 & 9.99e+02 &47&	17 & 1.28e+03 & 72&	48 & 1.59e+03 & 97&  36 & 3.00e+03 \\\hline
		23&	1 & 1.01e+03 &48&	14 & 1.29e+03 & 73&	31 & 1.60e+03 & 98&34& 3.02e+03\\\hline
		24&	9 & 1.03e+03 &49&	92 & 1.30e+03 & 74&	65 & 1.61e+03 & 99&91&3.05e+03\\\hline
		25&	 6 & 1.04e+03 & 50 & 4 & 1.31e+03 & 75& 96& 1.64e+03& 100&100&3.89e+03\\\hline
	\end{tabular}
\end{table}

\begin{table}[h!]
	\caption{Best $\mu$ obtained for A14 and \sherpa datasets when unfiltered (\textit{All data}), bin filtered and observable filtered data are used for parameter tuning.
	}
	\label{tab:allbestmu}
	
	\centering
	\begin{tabular}{c|c|c|c}
		\hline
		Dataset & All data & Bin filtered & Observable filtered\\
		\hline
		A14 &80&76&80\\\hline
		\sherpa &82&71&73\\\hline
	\end{tabular}
\end{table}

\subsection{Outlier observables in the A14 dataset}
\label{sec:A14outliers}
There are 12 outlier observables using the cubic polynomial approximation and 9 outlier observables using the rational approximation in the A14 dataset.

\begin{tabularx}{\textwidth}{l|l}
\hline
  Cubic Polynomial Model & Rational Approximation Model \\ \hline
\tt    /ATLAS\_2011\_I919017/d01-x02-y02 & \tt /ATLAS\_2011\_I919017/d01-x02-y02	\\
\tt	/ATLAS\_2011\_I919017/d01-x02-y03  & \tt  /ATLAS\_2011\_I919017/d01-x04-y04 	 \\
\tt	/ATLAS\_2011\_I919017/d01-x03-y02 & \tt       /ATLAS\_2011\_I919017/d02-x04-y03\\	
\tt	/ATLAS\_2011\_I919017/d01-x03-y07  & \tt /ATLAS\_2011\_I919017/d02-x04-y04  \\
\tt	/ATLAS\_2011\_I919017/d01-x04-y07 &  \tt       /ATLAS\_2011\_I919017/d02-x04-y05\\
\tt	/ATLAS\_2011\_I919017/d01-x04-y08 &  \tt  /ATLAS\_2011\_I919017/d02-x04-y09  \\
\tt	/ATLAS\_2011\_I919017/d01-x04-y09  & \tt  /ATLAS\_2011\_I919017/d02-x04-y10 \\
\tt	/ATLAS\_2011\_I919017/d02-x04-y04  &  \tt 	/ATLAS\_2011\_I919017/d02-x04-y14\\
\tt	/ATLAS\_2011\_I919017/d02-x04-y10 &   \tt     	/ATLAS\_2011\_I919017/d02-x04-y15  \\
\tt	/ATLAS\_2011\_I919017/d02-x04-y13 & \\
\tt	/ATLAS\_2011\_I919017/d02-x04-y14 &   \tt \\
\tt    /ATLAS\_2011\_I919017/d02-x04-y15  &  \tt  \\\hline
  \end{tabularx}

\subsection{Outlier observables in the \sherpa dataset}
\label{sec:sherpaoutliers}	
There are 2 outlier observables using the cubic polynomial approximation and 3 outlier observables using the rational approximation in the \sherpa dataset.
\begin{tabularx}{\textwidth}{l|l}
\hline
  Cubic Polynomial Model & Rational Approximation Model \\ \hline
\tt   /DELPHI\_1996\_S3430090/d07-x01-y01 & \tt /DELPHI\_1996\_S3430090/d02-x01-y01	\\
\tt	/DELPHI\_1996\_S3430090/d08-x01-y01 & \tt  /DELPHI\_1996\_S3430090/d07-x01-y01 	 \\
\tt	 & \tt       /DELPHI\_1996\_S3430090/d08-x01-y01\\	\hline
  \end{tabularx}

\subsection{Bin filtered data for A14 dataset}
\label{sec:binfiltA14}
In Table~\ref{tab:A14Binfiltdata}, we give the names of the A14 observables from which bins have been filtered, the number of bins filtered out, critical $\chisq$ value, and $\chisq$ test statistic before and after filtering the bins.

\begin{table}[ht!]
	\centering
	\tiny
	\caption{ \tiny Bin filtering of  A14 data: Shown are the observables from which bins were removed and the number of bins removed. We also show the critical $\chi^2$ values and the $\chi^2$ test statistic before and after bin filtering. If all the bins were removed from the observable then the number of bins removed is shown in bold font and the $\chi^2$ test statistic before and after bin filtering is the same.} \label{tab:A14Binfiltdata}
	\begin{tabular}{|c|c|c|c|c|}\hline
		Observable Name&\shortstack{No. of \\filtered bins}& $\chi^2_{c,\calB} $ & \shortstack{$\chi^2_\calB$ before \\filtering bins} & \shortstack{$\chi^2_\calB$ after \\filtering bins}\\\hline
		/ATLAS\_2011\_I919017/d01-x02-y04 & \textbf{11} & 3.84 & 9.77 & 9.77\\\hline
		/ATLAS\_2011\_I919017/d01-x02-y05 & \textbf{13} & 7.81 & 13.51 & 13.51\\\hline
		/ATLAS\_2011\_I919017/d01-x02-y13 & \textbf{11} & 3.84 & 9.43 & 9.43\\\hline
		/ATLAS\_2011\_I919017/d01-x02-y18 & \textbf{11} & 3.84 & 6.20 & 6.20\\\hline
		/ATLAS\_2011\_I919017/d01-x03-y01 & 11 & 21.03 & 24.00 & 3.57\\\hline
		/ATLAS\_2011\_I919017/d01-x03-y02 & 4 & 21.03 & 48.57 & 19.72\\\hline
		/ATLAS\_2011\_I919017/d01-x03-y03 & 2 & 25.00 & 28.99 & 24.72\\\hline
		/ATLAS\_2011\_I919017/d01-x03-y04 & 2 & 32.67 & 35.36 & 32.19\\\hline
		/ATLAS\_2011\_I919017/d01-x03-y06 & 10 & 26.30 & 59.81 & 26.13\\\hline
		/ATLAS\_2011\_I919017/d01-x03-y07 & 7 & 25.00 & 58.78 & 23.91\\\hline
		/ATLAS\_2011\_I919017/d01-x03-y08 & 5 & 28.87 & 36.98 & 28.27\\\hline
		/ATLAS\_2011\_I919017/d01-x03-y09 & 6 & 35.17 & 41.10 & 35.10\\\hline
		/ATLAS\_2011\_I919017/d01-x03-y12 & 3 & 23.68 & 31.51 & 21.33\\\hline
		/ATLAS\_2011\_I919017/d01-x03-y13 & 15 & 31.41 & 58.77 & 26.18\\\hline
		/ATLAS\_2011\_I919017/d01-x03-y14 & 12 & 33.92 & 69.51 & 32.38\\\hline
		/ATLAS\_2011\_I919017/d01-x03-y17 & 3 & 23.68 & 30.48 & 22.60\\\hline
		/ATLAS\_2011\_I919017/d01-x03-y18 & 1 & 30.14 & 30.65 & 26.75\\\hline
		/ATLAS\_2011\_I919017/d01-x03-y19 & 12 & 33.92 & 43.45 & 6.13\\\hline
		/ATLAS\_2011\_I919017/d01-x04-y03 & \textbf{22} & 21.03 & 45.11 & 45.11\\\hline
		/ATLAS\_2011\_I919017/d01-x04-y04 & \textbf{21} & 19.68 & 93.99 & 93.99\\\hline
		/ATLAS\_2011\_I919017/d01-x04-y05 & 4 & 16.92 & 22.81 & 16.74\\\hline
		/ATLAS\_2011\_I919017/d01-x04-y08 & \textbf{22} & 21.03 & 65.21 & 65.21\\\hline
		/ATLAS\_2011\_I919017/d01-x04-y09 & \textbf{22} & 21.03 & 71.99 & 71.99\\\hline
		/ATLAS\_2011\_I919017/d01-x04-y10 & 2 & 18.31 & 25.18 & 18.27\\\hline
		/ATLAS\_2011\_I919017/d01-x04-y13 & 12 & 22.36 & 49.09 & 2.36\\\hline
		/ATLAS\_2011\_I919017/d01-x04-y14 & \textbf{24} & 23.68 & 71.30 & 71.30\\\hline
		/ATLAS\_2011\_I919017/d01-x04-y15 & 4 & 21.03 & 27.53 & 20.80\\\hline
		/ATLAS\_2011\_I919017/d01-x04-y18 & 2 & 23.68 & 23.77 & 22.57\\\hline
		/ATLAS\_2011\_I919017/d01-x04-y19 & 8 & 22.36 & 36.75 & 16.78\\\hline
		/ATLAS\_2011\_I919017/d01-x04-y25 & 3 & 26.30 & 29.14 & 24.98\\\hline
		/ATLAS\_2011\_I919017/d02-x02-y05 & 1 & 11.07 & 13.84 & 9.87\\\hline
		/ATLAS\_2011\_I919017/d02-x02-y09 & 1 & 9.49 & 12.32 & 8.52\\\hline
		/ATLAS\_2011\_I919017/d02-x02-y14 & 1 & 9.49 & 12.19 & 9.18\\\hline
		/ATLAS\_2011\_I919017/d02-x03-y02 & 15 & 30.14 & 40.31 & 7.63\\\hline
		/ATLAS\_2011\_I919017/d02-x03-y06 & 3 & 31.41 & 36.64 & 28.59\\\hline
		/ATLAS\_2011\_I919017/d02-x03-y07 & 4 & 31.41 & 55.51 & 29.12\\\hline
		/ATLAS\_2011\_I919017/d02-x03-y12 & 7 & 31.41 & 45.41 & 30.04\\\hline
		/ATLAS\_2011\_I919017/d02-x03-y17 & 1 & 30.14 & 30.64 & 28.11\\\hline
		/ATLAS\_2011\_I919017/d02-x04-y03 & 10 & 26.30 & 46.87 & 19.20\\\hline
		/ATLAS\_2011\_I919017/d02-x04-y04 & \textbf{25} & 25.00 & 136.83 & 136.83\\\hline
		/ATLAS\_2011\_I919017/d02-x04-y05 & \textbf{28} & 28.87 & 74.75 & 74.75\\\hline
		/ATLAS\_2011\_I919017/d02-x04-y08 & 16 & 27.59 & 82.23 & 25.29\\\hline
		/ATLAS\_2011\_I919017/d02-x04-y09 & \textbf{27} & 27.59 & 156.13 & 156.13\\\hline
		/ATLAS\_2011\_I919017/d02-x04-y10 & \textbf{30} & 31.41 & 126.00 & 126.00\\\hline
		/ATLAS\_2011\_I919017/d02-x04-y13 & 14 & 26.30 & 71.23 & 23.47\\\hline
		/ATLAS\_2011\_I919017/d02-x04-y14 & \textbf{27} & 27.59 & 103.20 & 103.20\\\hline
		/ATLAS\_2011\_I919017/d02-x04-y15 & 9 & 28.87 & 70.47 & 26.49\\\hline
		/ATLAS\_2011\_I919017/d02-x04-y18 & 3 & 26.30 & 32.01 & 25.80\\\hline
		/ATLAS\_2011\_I919017/d02-x04-y19 & 13 & 28.87 & 67.53 & 23.91\\\hline
		/ATLAS\_2011\_I919017/d02-x04-y20 & 10 & 28.87 & 57.69 & 28.31\\\hline
		/ATLAS\_2011\_I919017/d02-x04-y24 & 10 & 28.87 & 43.46 & 28.24\\\hline
		/ATLAS\_2011\_I919017/d02-x04-y25 & 3 & 31.41 & 39.98 & 28.20\\\hline
		/ATLAS\_2011\_ZPT/d02-x01-y01 & 1 & 14.07 & 15.77 & 14.06\\\hline
		/ATLAS\_2011\_ZPT/d02-x02-y02 & 2 & 14.07 & 16.94 & 13.93\\\hline
		/ATLAS\_2011\_ZPT/d03-x01-y01 & 1 & 14.07 & 15.32 & 13.85\\\hline
		/ATLAS\_2013\_JETUE/d08-x01-y03 & 1 & 12.59 & 19.97 & 11.51\\\hline
	\end{tabular}
\end{table}

\subsection{Bin filtered data for \sherpa dataset}
\label{sec:binfiltSherpa}
In Table~\ref{tab:SherpaBinfiltdata}, we give the names of the \sherpa observables from which bins have been filtered, the number of bins filtered out, critical $\chisq$ value, and $\chisq$ test statistic before and after filtering the bins.

\begin{table}[h!]
	\centering
	\caption{Bin filtering of \sherpa data: Shown are the observables from which bins were removed and the number of bins removed. We also show the critical $\chi^2$ values and the $\chi^2$ test statistic before and after bin filtering. If all the bins were removed from the observable then the number of bins removed is shown in bold font and the $\chi^2$ test statistic before and after bin filtering is the same.} \label{tab:SherpaBinfiltdata}
	\begin{tabular}{|c|c|c|c|c|}\hline
		Observable Name&\shortstack{No. of \\filtered bins}& $\chi^2_{c,\calB} $ & \shortstack{$\chi^2_\calB$ before \\filtering bins} & \shortstack{$\chi^2_\calB$ after \\filtering bins}\\\hline
		/DELPHI\_1996\_S3430090/d02-x01-y01 & \textbf{17} & 9.49 & 35.19 & 35.19\\\hline
		/DELPHI\_1996\_S3430090/d04-x01-y01 & \textbf{17} & 9.49 & 24.89 & 24.89\\\hline
		/DELPHI\_1996\_S3430090/d06-x01-y01 & \textbf{21} & 15.51 & 41.59 & 41.59\\\hline
		/DELPHI\_1996\_S3430090/d07-x01-y01 & \textbf{22} & 16.92 & 83.91 & 83.91\\\hline
		/DELPHI\_1996\_S3430090/d08-x01-y01 & \textbf{26} & 22.36 & 80.11 & 80.11\\\hline
		/DELPHI\_1996\_S3430090/d10-x01-y01 & 2 & 11.07 & 15.90 & 10.59\\\hline
		/DELPHI\_1996\_S3430090/d11-x01-y01 & \textbf{20} & 14.07 & 94.30 & 94.30\\\hline
		/DELPHI\_1996\_S3430090/d16-x01-y01 & \textbf{14} & 3.84 & 17.63 & 17.63\\\hline
		/DELPHI\_1996\_S3430090/d18-x01-y01 & \textbf{23} & 18.31 & 101.31 & 101.31\\\hline
		/DELPHI\_1996\_S3430090/d19-x01-y01 & \textbf{21} & 15.51 & 59.12 & 59.12\\\hline
		/DELPHI\_1996\_S3430090/d20-x01-y01 & \textbf{16} & 7.81 & 20.48 & 20.48\\\hline
		/DELPHI\_1996\_S3430090/d33-x01-y01 & 5 & 52.19 & 75.18 & 50.31\\\hline
	\end{tabular}
\end{table}

\subsection{Complete results from filtering out observables and bins}
\label{sec:filterresonlinesup}
\begin{figure}[ht!]
	\centering
	\includegraphics[width=\textwidth]{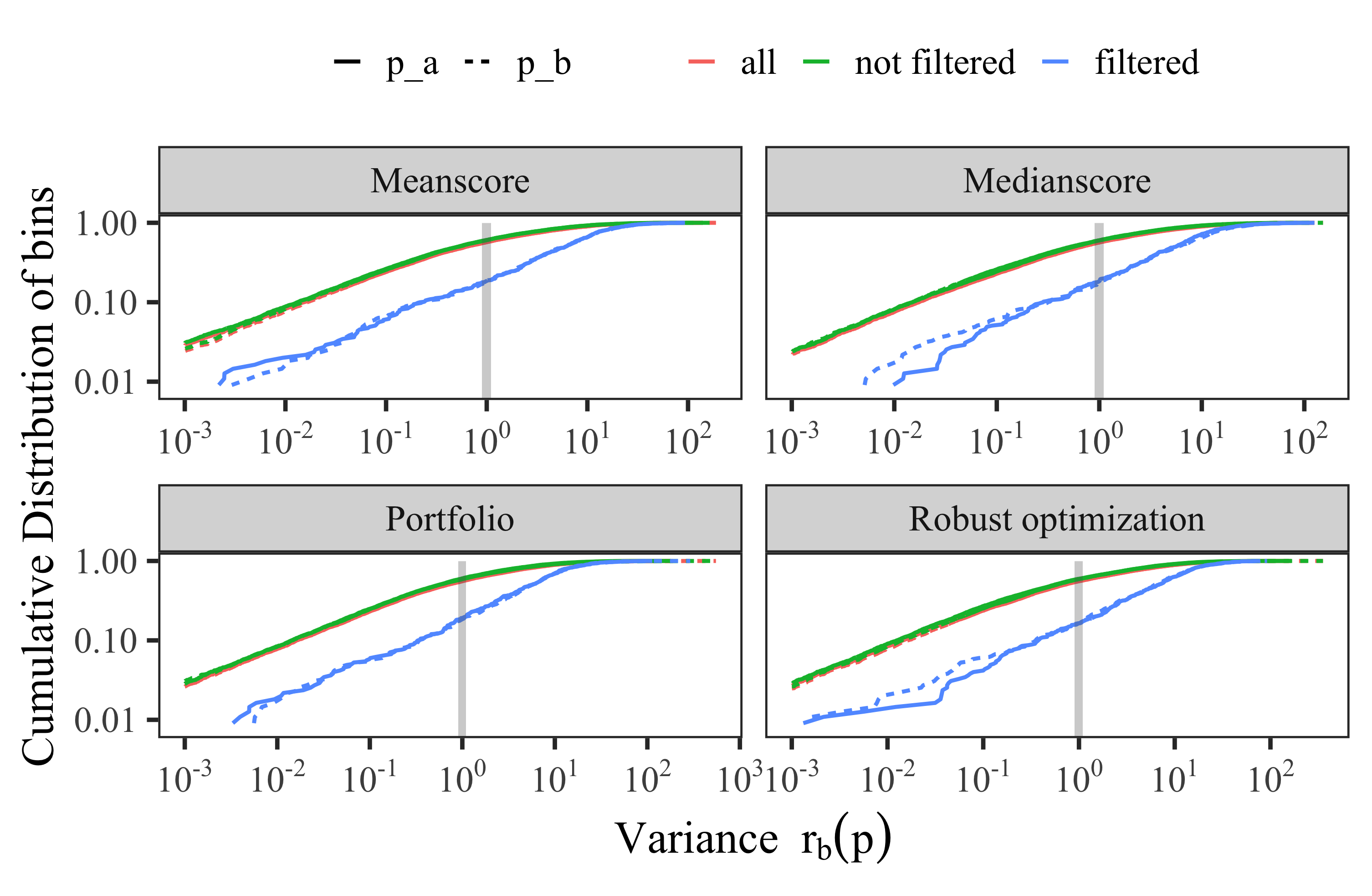}
	\caption{Cumulative distribution of bins for the A14 dataset
          at different bands of variance levels using different  approaches.
          Results are shown using the parameters $\p^*_a$ obtained
          using all bins during optimization, and the parameters
          $\p^*_b$ obtained when only the bin filtered bins are used during optimization.
        }
	\label{Fig:chi2perf2binFilterA14}        
\end{figure}

\begin{figure}[ht!]
	\centering 
	\includegraphics[width=\textwidth]{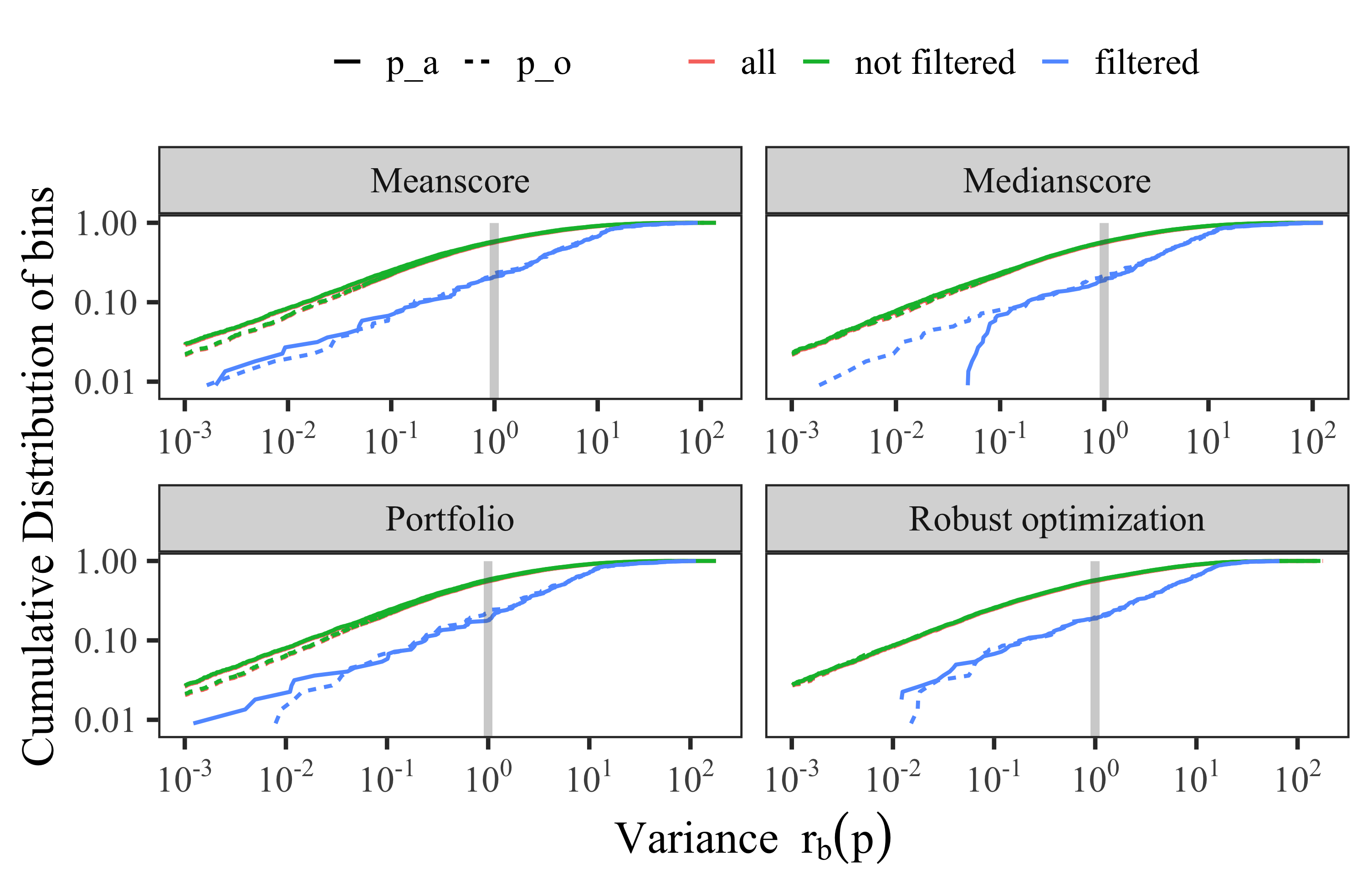}
        \caption{Same as Figure~\ref{Fig:chi2perf2binFilterA14}, but
          using observable filtering.}
	\label{Fig:chi2perf2obsFilterA14}
\end{figure}

\begin{figure}[ht!]
	\centering 
	\includegraphics[width=\textwidth]{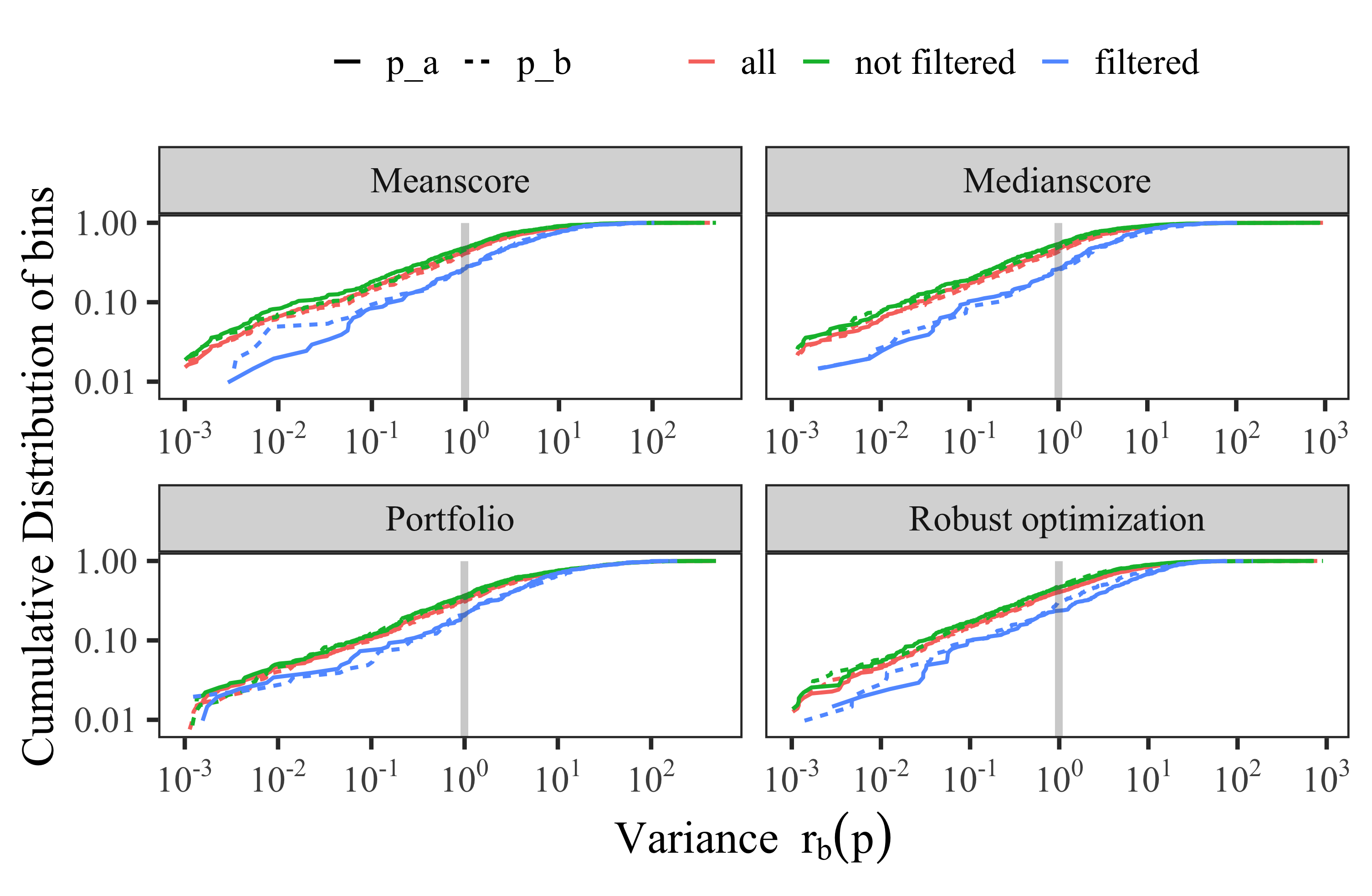}
        \caption{Same as Figure~\ref{Fig:chi2perf2binFilterA14} , but
          for the \sherpa dataset.}          
	\label{Fig:chi2perf2binFilterSherpa}
	
\end{figure}

\begin{figure}[ht!]
	\centering 
	\includegraphics[width=\textwidth]{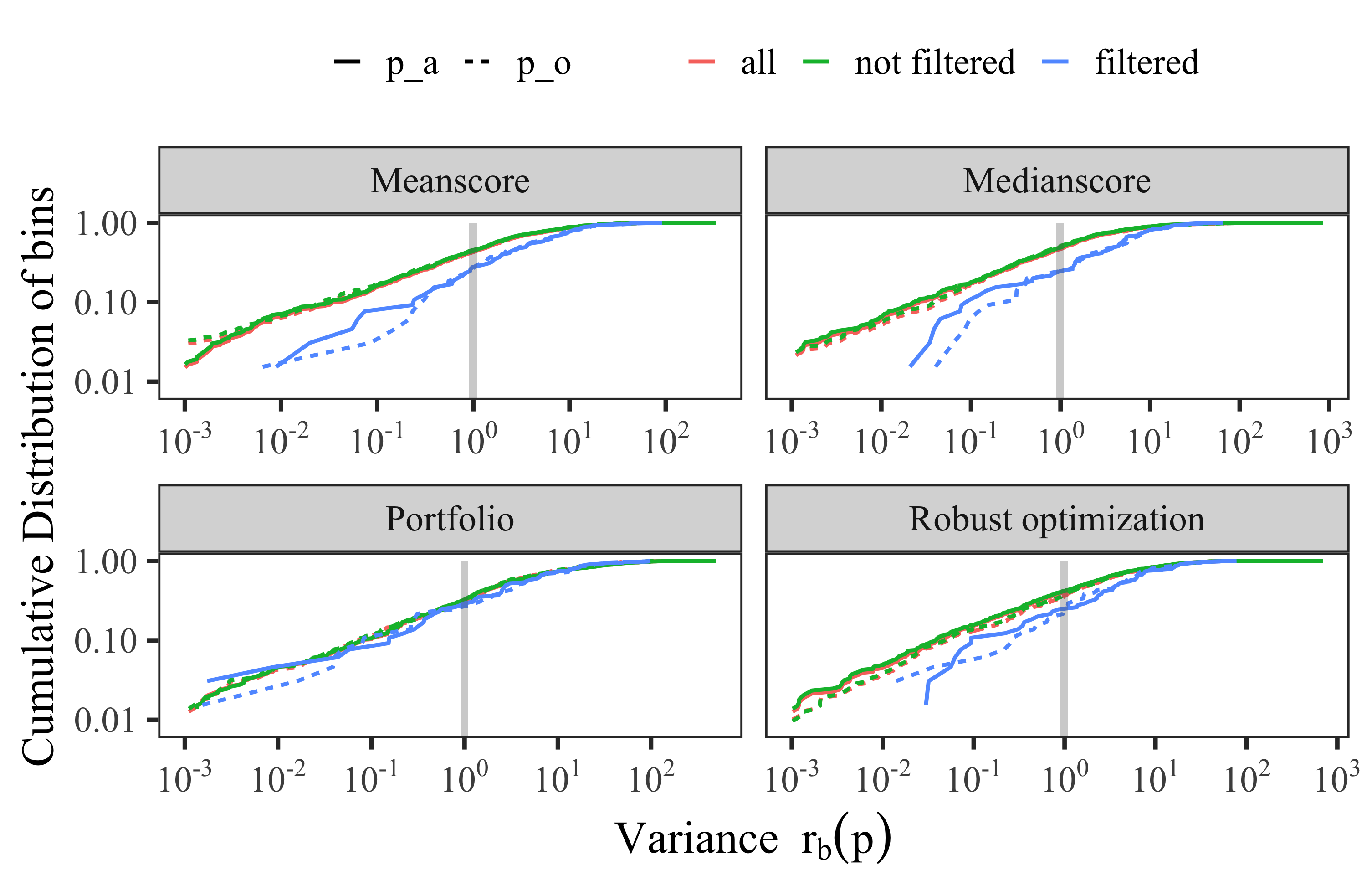}
        \caption{Same as Figure~\ref{Fig:chi2perf2binFilterSherpa},
          but using observable filtering.}
	\label{Fig:chi2perf2obsFilterSherpa}	
\end{figure}

In Figures~\ref{Fig:chi2perf2binFilterA14} and~\ref{Fig:chi2perf2obsFilterA14}, the cumulative distribution plots for parameters obtained after bin filtering and observable filtering for the A14 data are presented.  In Figures~\ref{Fig:chi2perf2binFilterSherpa} and~\ref{Fig:chi2perf2obsFilterSherpa}, the cumulative distribution plots for parameters obtained after bin filtering and observable filtering for the \sherpa data are presented.  
From these figures, we observe that there is no significant difference in the number of bins within the 1 $\sigma$ variance level between the optimal parameters $\p^*_a$ obtained when all bins were used for tuning and the optimal parameters $\p^*_b$ and $\p^*_o$ obtained when only the bin filtered and observable filtered bins are used for tuning, respectively.
This indicates that the MC generator cannot explain the bins removed by the filtering approaches very well. Hence, removing these bins from the tuning process does not reduce the information required to achieve a good tune as it performs very similarly to when all bins are used for tuning.

\subsection{Comparison of the rational approximation with the MC generator}
\label{sec:mcvsra2}
\begin{figure}[!ht]
	\centering
	\includegraphics[width=\textwidth]{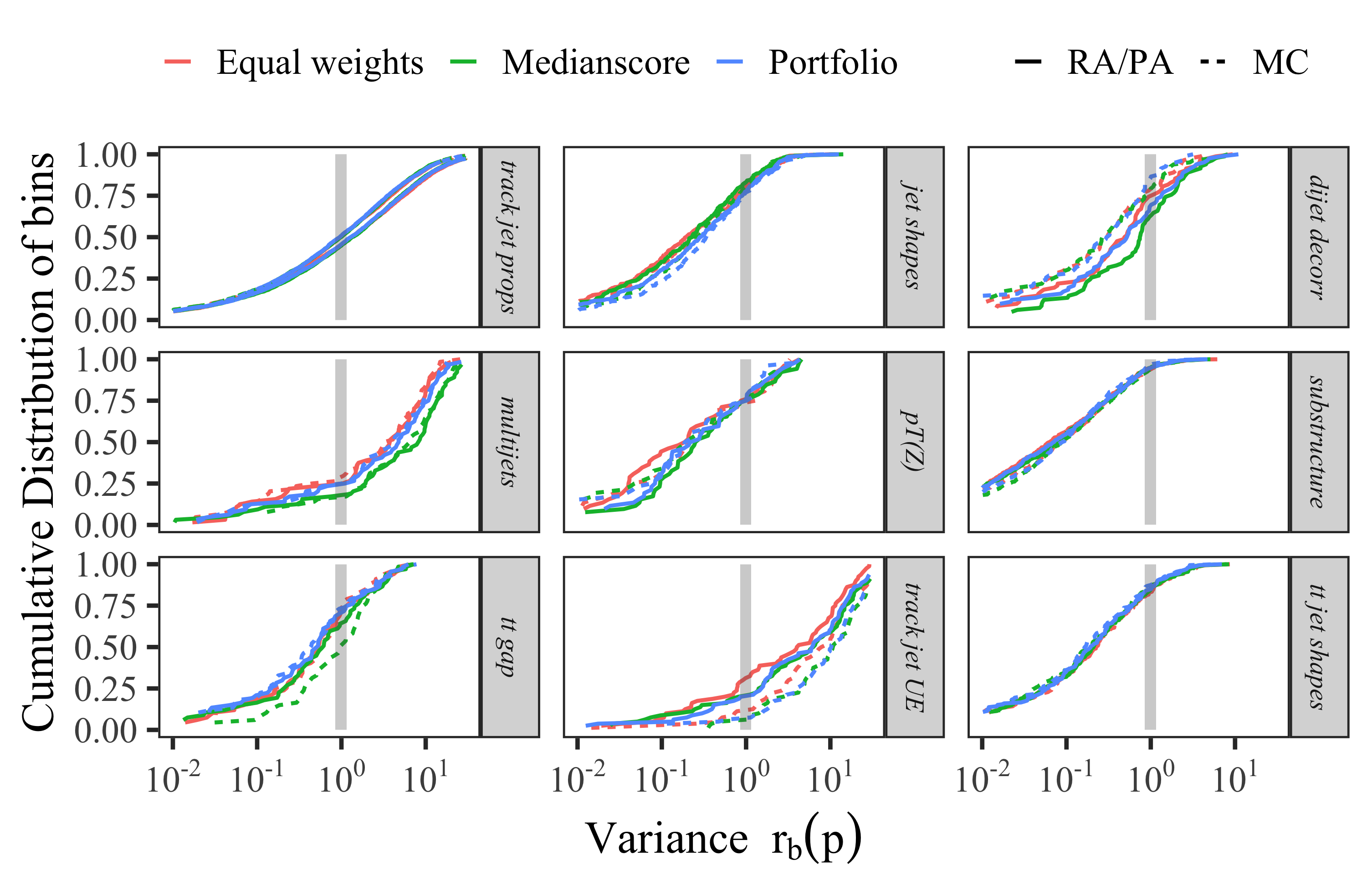}
	\caption{Cumulative distribution of bins in each category of the A14 dataset at different bands of variance levels computed with the rational approximation (RA) given by, $r_b(\p)=\frac{\left(f_b(\p) - \calR_b\right)^2}{\Delta f_b(\p)^2+\Delta\calR_b^2}$ and the MC simulation given by $\widetilde{r_b(\p)}=\frac{\left(\MC_b(\p) - \calR_b\right)^2}{\Delta \MC_b(\p)^2+\Delta\calR_b^2}$\\ 
	}
	\label{Fig:RAvsMCchi2perf2binsCat2}
	
\end{figure}

We compare the cumulative distribution of bins at different bands of variance levels computed using the rational approximation (RA) model as $r_b(\p)=\frac{\left(f_b(\p) - \calR_b\right)^2}{\Delta f_b(\p)^2+\Delta\calR_b^2}$ and the MC generator model as $\widetilde{r_b(\p)}=\frac{\left(\MC_b(\p) - \calR_b\right)^2}{\Delta \MC_b(\p)^2+\Delta\calR_b^2}$, where $\p$ are the parameters obtained from the different tuning approaches. In Figure~\ref{Fig:RAvsMCchi2perf2binsCat}, we showed the plot of this comparison for bins in each category of the A14 dataset using the parameters from three approaches. For completeness, in Figure~\ref{Fig:RAvsMCchi2perf2binsCat2}, we show the plot of this comparison for the remaining  three approaches. 

We observe in Figure~\ref{Fig:RAvsMCchi2perf2binsCat2} that around the variance boundary, except for in the \textit{Track-jet UE} and \textit{Multijets} categories, there is no significant difference in performance  between $r_b(\p)$ and $\widetilde{r_b(\p)}$ for each approach. 
In the case of \textit{Track-jet UE} and  \textit{Multijets} categories, the number of bins that lie within the variance boundary is quite low compared to other categories. This suggests that many bins in these categories cannot be explained well by either the MC generator or the approximation for the optimal tuning parameters reported by the approaches. Additionally, we observe in these categories that the approximations are not able to capture the MC generator  perfectly.

\subsection{Optimal parameter values for the A14 dataset with the rational approximation}
\label{sec:pstar_a14_31}

\begin{figure}[ht!]
	\centering
	\includegraphics[width=\textwidth]{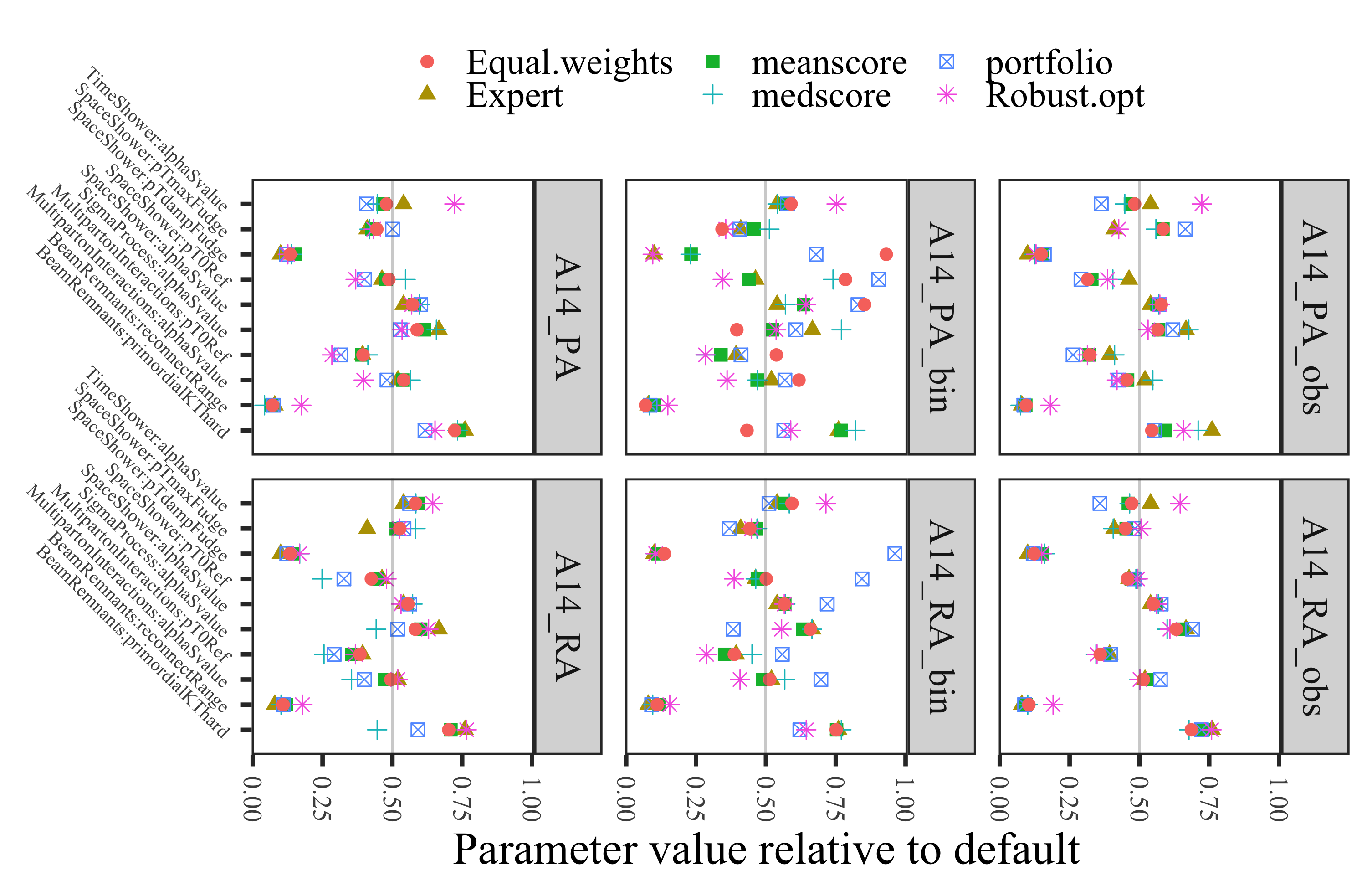}                        
	\caption{Optimal parameter values for the A14 dataset obtained
          when using all, bin-filtered (\_bin) and observable-filtered (\_obs) data in the optimization
          and the polynomial approximation (PA) and rational approximation (RA). Values are normalized to [0,1].}
	\label{fig:prms_all_31}
	\label{fig:prms_obsfilt_31}
	\label{fig:prms_binfilt_31}
 	\label{fig:prms_all_c}
 	\label{fig:prms__obsfilt_c}
 	\label{fig:prms__binfilt_c}                
\end{figure}

To better visually compare the different solutions obtained with our optimization methods, we show the [0,1]-scaled optimal parameter values in Figure~\ref{fig:prms_all_31}.%

\subsection{Results for using the cubic polynomial to approximate the MC simulation}
In the main paper, we showed the numerical results when using a rational approximation of the MC simulation during tuning. In the A14 publication~\cite{ATL-PHYS-PUB-2014-021}, a cubic polynomial was used. Thus, in this section, we present the results obtained with our optimization methods when using a cubic polynomial instead of a rational approximation.

\subsubsection{Comparison metric outcomes for the A14 dataset using the cubic polynomial approximation}
\label{sec:a14_30_metric}
Tables~\ref{tab:comparison_full}-\ref{tab:comparison_binf} show the comparison metrics we introduced in the main paper in Section~\ref{sec:metrics} when using the cubic polynomial approximation for the full data, the observable-filtered data, and the bin-filtered data, respectively. We see that for all three cases and most criteria (except for the D-optimality in the observable-filtered case), our automated methods for adjusting the observable weights perform better than the expert solution (i.e., using the parameters published in~\cite{ATL-PHYS-PUB-2014-021}).

\begin{table}[htbp]
	\caption{A14 results for the \textit{full data set} when using the \textit{cubic polynomial} approximation. Smaller numbers are better. The best results are indicated in bold.}%
	\centering
	\begin{tabular}{l|rrr}
		\hline
		Method         & Weighted $\chi^2$ & A-optimality & D-optimality (log) \\\hline
		Bilevel-meanscore  & 0.1290 & 	0.5358 & 	-66.0364
 \\
		Bilevel-medianscore    & 0.1645	 & \textbf{0.4114} & 	\textbf{-70.0545}
         \\
		Bilevel-portfolio   & 0.1900 & 	0.6590	 & -63.0378
\\
		Expert &0.1306	 & 0.5466 & 	-68.6511
 \\
		All-weights-equal  & 0.1034	 & 0.5553 & 	-65.6099
          \\
		Robust optimization    & \textbf{0.0697}	 & 0.9749 & 	-66.7931

  \\\hline
	\end{tabular}
	\label{tab:comparison_full}
\end{table}

\begin{table}[htbp]
	\caption{A14 results for the \textit{observable-filtered data} when using the \textit{cubic polynomial} approximation (12 observables were filtered out and not used during the optimization). Smaller numbers are better. Best results are indicated in bold. All values are computed over all 406 observables.}%
	\centering
	\begin{tabular}{l|rrr}
		\hline
		Method         & Weighted $\chi^2$ & A-optimality & D-optimality (log) \\\hline
		Bilevel-meanscore   & 0.1079  & 	0.8082  & 	-61.9210
\\
		Bilevel-medianscore    & 0.1702	  & \textbf{0.4955}  & 	-66.6920
  \\
		Bilevel-portfolio   & 0.1764  & 	0.7408	  & -61.3839
\\
		Expert &0.1306	  & 0.5466	  & \textbf{-68.6511}
 \\
		All-weights-equal  & 0.1049	  & 0.6689  & 	-63.6502
 \\
		Robust optimization    & \textbf{0.0829}	  & 1.0574	  & -66.3665

		\\\hline
	\end{tabular}
	\label{tab:comparison_obsf}
\end{table}

\begin{table}[htbp]
	\caption{A14 results for the \textit{bin-filtered data} when using the \textit{cubic polynomial} approximation. %
	1811 out of 7010 total bins were filtered out and not used during the optimization. Smaller numbers are better. Best results are indicated in bold. All results are computed over all bins.}
	\centering
	\begin{tabular}{l|rrr}
		\hline
		Method         & Weighted $\chi^2$ & A-optimality & D-optimality (log) \\\hline
		Bilevel-meanscore   & 0.1244 & 	0.7147 & 	-64.7848
  \\
  Bilevel-medianscore    & 0.2171 & 	0.5433	 & -69.7202
\\
		Bilevel-portfolio   & \textbf{0.1159} & 	0.5205	 & \textbf{-70.1573}
 \\
		Expert &0.1306 & 	0.5466	 & -68.6511
 \\
		All-weights-equal  & 0.1406	 & \textbf{0.4122} & 	-69.2732

\\
		Robust optimization    & 0.1234	 & 0.8075 & 	-67.1015

\\\hline
	\end{tabular}
	\label{tab:comparison_binf}
\end{table}

\subsubsection{Optimal parameter values for the A14 dataset using the cubic polynomial approximation}
Table~\ref{tab:optimal_prms}  shows the optimal values for the tuned parameters obtained by all methods for the A14 dataset when using all observables in the tune.  For Bilevel-meanscore, -medianscore and -portfolio, we repeated the experiments three times using different random number seeds and we report the best results among the three trials based on their respective objective functions. From these results, we can see that the Bilevel-medianscore method leads to a solution that is closest to the expert's solution. 

To better visually compare the different solutions obtained with our methods, we show the [0,1]-scaled optimal values in Figure~\ref{fig:prms_all_c}.
We can see that there are differences between the optimal parameters obtained with the different methods, in particular, the results of the robust optimization method tend to be further away from the expert's solution for parameters 1, 2, 3, 7, 8, 9 and 10. The results of the portfolio optimization differ from the expert tune in particular for parameters 1, 2, 3, 4 and 7. The mean- and medianscore results are very similar to each other as well as to the expert's solution.

\begin{table}[htbp]
	\centering
	\small
	\caption{Optimal parameter values for the A14 dataset obtained when using all observables in the optimization and the \textit{cubic polynomial} approximation.}\label{tab:optimal_prms}
	\begin{adjustbox}{width=\textwidth}
	\begin{tabular}{llllllll}
		\hline
		ID & Parameter name & Expert & Bil.-meanscore & Bil.-medianscore & Bil.-portfolio & Robust opt    & All-weights-equal \\\hline
1  & \tt SigmaProcess:alphaSvalue            & 0.143  & 0.139 & 0.141 & 0.140 & 0.136 & 0.138 \\
2  & \tt BeamRemnants:primordialKThard       & 1.904  & 1.867 & 1.884 & 1.866 & 1.826 & 1.862 \\
3  & \tt SpaceShower:pT0Ref                  & 1.643  & 1.632 & 1.735 & 1.651 & 1.395 & 1.603 \\
4  & \tt SpaceShower:pTmaxFudge              & 0.908  & 0.939 & 0.904 & 0.988 & 0.933 & 0.944 \\
5  & \tt SpaceShower:pTdampFudge             & 1.046  & 1.079 & 1.069 & 1.047 & 1.063 & 1.067 \\
6  & \tt SpaceShower:alphaSvalue             & 0.123 & 0.129 & 0.130 & 0.130 & 0.128 & 0.129 \\
7  & \tt TimeShower:alphaSvalue              & 0.128 & 0.123 & 0.124 & 0.121 & 0.136 & 0.124 \\
8  & \tt MultipartonInteractions:pT0Ref      & 2.149  & 2.083 & 2.065 & 2.039 & 1.925 & 2.092 \\
9  & \tt MultipartonInteractions:alphaSvalue & 0.128
 & 0.127 & 0.127 & 0.126 & 0.120 & 0.127 \\
10 & \tt BeamRemnants:reconnectRange         & 1.792  & 1.531 & 1.405 & 1.591 & 2.567 & 1.636
\\\hline\hline
		& Euclidean distance from the expert solution & &0.246&	0.235&	0.428&	0.451&	0.259

	\\\hline %
	\end{tabular}
	\end{adjustbox}
\end{table}

We conducted a similar analysis on the observable- and bin-filtered data. Table~\ref{tab:optimal_prms_obsfilter} shows the optimal parameter values that we obtain with the automated optimization methods after filtering out the 12 observables that the model cannot explain (see also Section~\ref{sec:obsfilter}). The expert solution is the same as before and based on all observables. We include it for easier comparison. With only a few exceptions, all parameters obtained with the automated optimizations change (as compared to using the full dataset). 
Figure~\ref{fig:prms__obsfilt_c} shows the optimal parameter values obtained with each method scaled to [0,1]. In comparison to when using the full dataset, we see that the results of the robust optimization now agree better with the expert's tune for parameters 3, 4, and 8, but less agreement is achieved for parameter 10. Of the three bilevel methods, the medianscore objective function leads to optimal parameters that are most similar to the expert tune. %

\begin{table}[htbp]
	\small
	\centering
	\caption{Optimal parameter values for A14 when using the \textit{cubic polynomial} approximation with all methods after outlier detection to filter out observables that cannot be approximated well by the model.}\label{tab:optimal_prms_obsfilter}
	\begin{adjustbox}{width=\textwidth}
	\begin{tabular}{llllllll}
		\hline
	ID	& Parameter name & Expert & Bilevel-meanscore & Bilevel-medianscore & Bilevel-portfolio & Robust opt    & All-weights-equal \\\hline
1  & \tt SigmaProcess:alphaSvalue            & 0.143  & 0.136 & 0.141 & 0.137 & 0.136 & 0.137 \\
2  & \tt BeamRemnants:primordialKThard       & 1.904  & 1.793 & 1.853 & 1.754 & 1.829 & 1.772 \\
3  & \tt SpaceShower:pT0Ref                  & 1.643  & 1.329 & 1.369 & 1.218 & 1.425 & 1.301 \\
4  & \tt SpaceShower:pTmaxFudge              & 0.908  & 1.079 & 1.088 & 1.223 & 0.926 & 1.085 \\
5  & \tt SpaceShower:pTdampFudge             & 1.046  & 1.069 & 1.053 & 1.101 & 1.065 & 1.074 \\
6  & \tt SpaceShower:alphaSvalue             & 0.123 & 0.129 & 0.128 & 0.129 & 0.129 & 0.129 \\
7  & \tt TimeShower:alphaSvalue              & 0.128 & 0.124 & 0.123 & 0.116 & 0.136 & 0.124 \\
8  & \tt MultipartonInteractions:pT0Ref      & 2.149  & 1.971 & 2.098 & 1.870 & 1.971 & 1.983 \\
9  & \tt MultipartonInteractions:alphaSvalue & 0.128 & 0.122 & 0.126 & 0.120 & 0.121 & 0.123 \\
10 & \tt BeamRemnants:reconnectRange         & 1.792  & 1.812 & 1.614 & 1.714 & 2.632 & 1.851
\\\hline 		\hline
& Euclidean distance from the expert solution & &0.447&	0.279&	0.553&	0.432&	0.480

		\\\hline %
	\end{tabular}
	\end{adjustbox}
\end{table}

In Table~\ref{tab:optimal_prms_binfilter} and Figure~\ref{fig:prms__binfilt_c} we show the optimal parameter values obtained with  our  methods after  applying the bin-filtering approach described in Section~\ref{sec:binfilter} in the main document.  In comparison to our results that do not use any filtering, we can see a much larger disagreement in the optimal parameters for all methods. In fact, all methods yield optimal parameters that are significantly further away from the expert's solution, except for parameters 7 and 10. The Euclidean distance between the optimal parameters obtained by our proposed methods and the expert solution shows that the bilevel-medianscore method leads to the most similar parameter values  while all the other methods lead to very different tunes.

\begin{table}[htbp]
	\centering
	\small
	\caption{Optimal  parameter  values  obtained  for  A14  with  the  \textit{cubic polynomial} approximation with  all methods after using the bin-filtering approach that excludes individual bins from the optimization.}\label{tab:optimal_prms_binfilter}
	\begin{adjustbox}{width=\textwidth}
	\begin{tabular}{llllllll}
		\hline
	ID	& Parameter name &  Expert & Bilevel-meanscore & Bilevel-medianscore & Bilevel-portfolio & Robust opt    & All-weights-equal \\\hline
1  & \tt SigmaProcess:alphaSvalue            & 0.143  & 0.141 & 0.143 & 0.136 & 0.136 & 0.132 \\
2  & \tt BeamRemnants:primordialKThard       & 1.904  & 1.919 & 1.918 & 1.575 & 1.794 & 1.716 \\
3  & \tt SpaceShower:pT0Ref                  & 1.643  & 1.802 & 2.284 & 2.300 & 1.355 & 2.123 \\
4  & \tt SpaceShower:pTmaxFudge              & 0.908  & 0.968 & 1.014 & 0.920 & 0.856 & 0.843 \\
5  & \tt SpaceShower:pTdampFudge             & 1.046  & 1.071 & 1.147 & 1.442 & 1.047 & 1.465 \\
6  & \tt SpaceShower:alphaSvalue             & 0.123 & 0.130 & 0.130 & 0.144 & 0.132 & 0.143 \\
7  & \tt TimeShower:alphaSvalue              & 0.128 & 0.129 & 0.127 & 0.131 & 0.138 & 0.130 \\
8  & \tt MultipartonInteractions:pT0Ref      & 2.149  & 2.059 & 1.800 & 2.228 & 1.925 & 2.306 \\
9  & \tt MultipartonInteractions:alphaSvalue & 0.128 & 0.126 & 0.120 & 0.131 & 0.118 & 0.131 \\
10 & \tt BeamRemnants:reconnectRange         & 1.792  & 1.860 & 1.922 & 1.807 & 2.340 & 1.622
\\\hline	\hline
		& Euclidean distance from the expert solution & & 0.376	&0.354&	0.848&	0.525&	1.111

\\\hline %
 \end{tabular}
	\end{adjustbox}
\end{table}

	\subsubsection{Comparison of optimal weights for the A14 dataset with cubic polynomial approximation}
In Table~\ref{tab:weights_c} we present the optimal weights assigned to each observable group by each method following the presentation style in \cite{ATL-PHYS-PUB-2014-021}. The weights reported for our  method are  averages of the  weights over all observables that belong to the same group. We scaled the weights such that they are on equal footing (all add up to 4580). %

The largest differences between the expert-adjusted  values and the values determined by our methods are for \textit{Multijets}, \textit{$t \Bar{t}$ gap} and \textit{Jet UE}, while for the remaining  groups, the values are very similar. These results, together with our analysis above let us conclude that an automated method  for adjusting the weights of observables  for tuning parameters is a viable approach and can  lead to better results than hand-tuning.
  \begin{table}[htbp]
		\centering
		\small
		\caption{Comparison of the optimal weights obtained by each method using the \textit{cubic polynomial} approximation. The observable grouping corresponds to the same grouping used in~\cite{ATL-PHYS-PUB-2014-021}.}\label{tab:weights_c}
		\begin{adjustbox}{width=\textwidth}
			\begin{tabular}{p{9cm}p{1.5cm}p{1.5cm}p{1.5cm}p{1.5cm}p{1.5cm}}
				\hline
				& expert & Bilevel-meanscore & Bilevel-medianscore  & Bilevel-portfolio & robustopt \\\hline
			\textbf{Track jet properties}                                        &     &       &       &       &       \\
			Charged jet multiplicity (50   distributions)                                           & 10  & 10.74
 & 14.98  & 10.64  & 19.38   \\
			Charged jet $z$ (50 distributions)                                                      & 10  & 11.29 & 8.66 & 13.71  & 0.00   \\
			Charged jet $p^{rel}_T$ (50   distributions)                                            & 10  & 11.20 & 10.39  & 10.99  &  0.00  \\
			Charged jet $\rho_{ch}(r)$ (50   distributions)                                         & 10  & 11.57 & 10.58 & 12.55 & 0.00 \\\hline
			\textbf{Jet shapes}                                                    &     &       &       &       &       \\
			Jet shape $\rho$ (59 distributions)                                                     & 10  & 11.57 & 11.06 & 10.20 & 19.38 \\\hline
			\textbf{Dijet decorr}                                                  &     &       &       &       &       \\
			Decorrelation $\Delta \phi$ (Fit range:   $\Delta \phi>0.75$) (9 distributions)         & 20  & 12.39 & 8.37  & 9.39  & 15.07 \\\hline
			\textbf{Multijets}                                                     &     &       &       &       &       \\
			3-to-2 jet ratios (8 distributions)                                                     & 100 & 12.99 & 27.19 & 5.88  & 19.38 \\\hline
			\textbf{$p^Z_T$} (Fit range: $p^Z_T<50   \text{GeV}$)                  &     &       &       &       &       \\
			Z-boson $p_T$ (20 distributions)                                                        & 10  & 12.78 & 14.53 & 6.71  & 19.38 \\\hline
			\textbf{Substructure}                                                  &     &       &       &       &       \\
			Jet mass, $\sqrt{d_{12}},  \sqrt{d_{23}},  \tau_{21}, \tau_{23}$ (36 distributions)     & 5   & 10.55 & 9.91  & 9.74  & 15.61 \\\hline
			\textbf{$t \Bar{t}$ gap}                                               &     &       &       &       &       \\
			Gap fraction vs $Q_0$,  $Q_{\text{sum}}$ for $|y|<0.8$                                  & 100 & 0.18  & 2.10  & 3.88  & 19.38 \\
			Gap fraction vs $Q_0$,  $Q_{\text{sum}}$ for $0.8<|y|<1.5$                              & 80  & 0.75  & 9.52  & 5.71  & 19.38 \\
			Gap fraction vs $Q_0$,  $Q_{\text{sum}}$ for $1.5<|y|<2.1$                              & 40  & 7.93  & 8.31  & 39.20 & 19.38 \\
			Gap fraction vs $Q_0$,  $Q_{\text{sum}}$ for $|y|<2.1$                                  & 10  & 18.19 & 13.43 & 11.05 & 19.38 \\\hline
			\textbf{Track-jet UE}                                                  &     &       &       &       &       \\
			Transverse region $N_{ch}$ profiles (5   distributions)                                 & 10  & 15.87 & 13.45 & 13.53 & 19.38 \\
			Transverse region mean $p_T$ profiles for   $R=0.2, 0.4, 0.6$ (3 distributions)         & 10  & 7.56  & 11.72 & 10.30 & 19.38 \\\hline
			\textbf{$t \Bar{t}$ jet shapes}                                       &     &       &       &       &       \\
			Jet shapes $\rho(r),  \psi(r)$ (20 distributions)                                       & 5   & 10.86 & 10.91 & 12.25 & 10.66  \\\hline
			\textbf{Jet UE}                                                        &     &       &       &       &       \\
			Transverse, trans-max, trans-min sum   $p_T$ incl. profiles (3 distributions)           & 20  & 12.76 & 22.51 & 9.65  & 19.38 \\
			Transverse, trans-max, trans-min $N_{ch}$   incl. profiles (3 distributions)            & 20  & 15.57 & 9.65  & 6.01  & 19.38 \\
			Transverse sum $E_T$ incl. profiles (2   distributions)                                 & 20  & 12.71 & 12.75 & 25.03 & 3.73  \\
			Transverse sum $ET/$sum $p_T$ ratio incl.,   excl. profiles (2 distributions)            & 5   & 7.53  & 18.29 & 28.35 & 19.38 \\
			Transverse mean $p_T$ incl. profiles (2   distributions)                                & 10  & 7.65  & 7.45  & 13.34 & 19.38 \\
			Transverse, trans-max, trans-min sum   $p_T$ incl. distributions (15 distributions)     & 1   & 9.39  & 5.50  & 11.04 & 19.38 \\
			Transverse,  trans-max, trans-min sum $N_{ch}$ incl.   distributions (15 distributions) & 1   & 11.92 & 9.85  & 14.52 & 19.38 \\\hline
			\end{tabular}
		\end{adjustbox}
	\end{table}

\subsubsection{Optimal parameter values for the \sherpa dataset with rational approximation}
\label{sec:pstar_sherpa_31}

\begin{figure}[ht!]
	\centering
        \includegraphics[width=\textwidth]{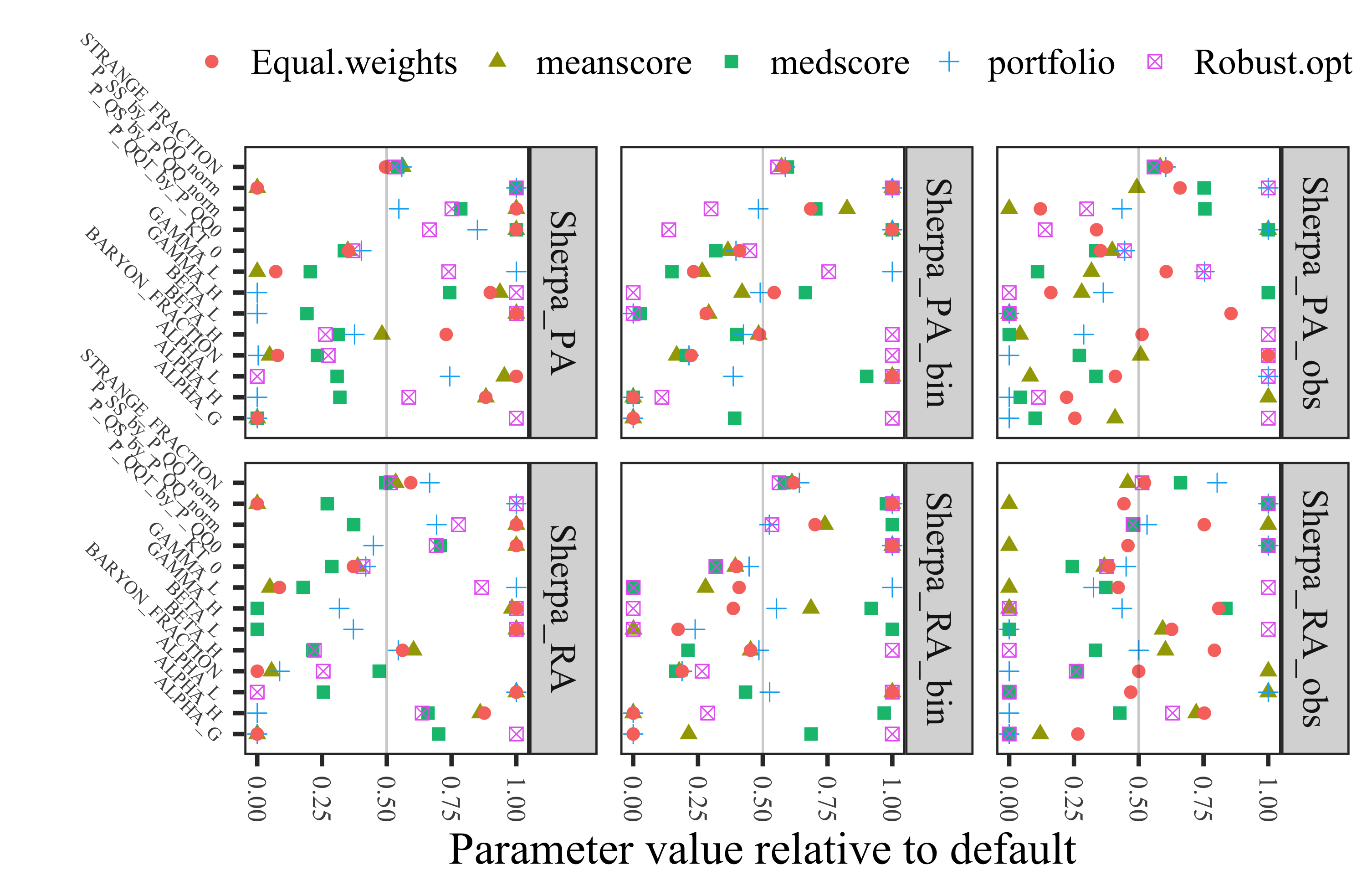}                
	\caption{Comparison of the optimal parameter values for
          \sherpa obtained with the different optimization methods
          when no, observable, and bin data filtering was applied and
          the \textit{rational} and polynomial  approximation was used. Values are normalized to [0,1].}
	\label{fig:prms_all_sherpa_31}
	\label{fig:prms_obsfil_sherpa_31}
        \label{fig:prms_binfilt_sherpa_31}
 	\label{fig:prms_all_sherpa}
 	\label{fig:prms_obsfil_sherpa}
 	\label{fig:prms_binfilt_sherpa}                        
\end{figure}

For a better visual comparison of the different solutions obtained with our methods, we show the [0,1]-scaled optimal values in Figure~\ref{fig:prms_all_sherpa_31}. %
Compared to the results for the A14 dataset, we see that there are significant differences between the optimal parameters obtained with the different  methods.

\subsubsection{Optimal parameter values for the \sherpa dataset with the cubic polynomial approximation}
The physics parameters $\p$ and their optimization ranges  are shown in Table~\ref{tab:prm-ranges-sherpa}. %
Tables~\ref{tab:optimal_prms_sherpa}, \ref{tab:optimal_prms_obsfilter_sherpa} and \ref{tab:optimal_prms_binfilter_sherpa} shows 
 the optimal values for the physics parameters obtained by all methods when no filtering was applied before optimization, after using outlier detection to remove observables from the optimization, and after using the bin-filtering approach that excludes individual bins from the optimization, respectively. For an illustrative comparison, we show the [0,1]-scaled optimal parameter values in Figure~\ref{fig:prms_all_sherpa}. %
 The default values lie right in the middle of the parameter range.

\begin{table}[htp!]
	\centering
	\small
	\caption{Optimal parameter values for the \sherpa dataset obtained with all methods using the \textit{cubic polynomial} approximation when no filtering was applied before optimization (88 observables).}\label{tab:optimal_prms_sherpa}
	\begin{adjustbox}{width=\textwidth}
	\begin{tabular}{lllllllll}
		\hline
		ID &	Parameter name	& Default & Bilevel-meanscore & Bilevel-medscore & Bilevel-portfolio & Robust opt & All-weights-equal   \\\hline
1  & \tt KT\_0                  & 1.00 & 0.850 & 0.837 & 0.903 & 0.870 & 0.853 \\
2  & \tt ALPHA\_G               & 1.25 & 0.626 & 0.626 & 0.626 & 1.874 & 0.626 \\
3  & \tt ALPHA\_L               & 2.50 & 3.634 & 2.022 & 3.108 & 1.252 & 3.749 \\
4  & \tt BETA\_L                & 0.10 & 0.150 & 0.069 & 0.050 & 0.150 & 0.150 \\
5  & \tt GAMMA\_L               & 0.50 & 0.250 & 0.353 & 0.750 & 0.619 & 0.286 \\
6  & \tt ALPHA\_H               & 2.50 & 3.455 & 2.047 & 1.251 & 2.712 & 3.454 \\
7  & \tt BETA\_H                & 0.75 & 0.736 & 0.610 & 0.657 & 0.573 & 0.922 \\
8  & \tt GAMMA\_H               & 0.10 & 0.144 & 0.124 & 0.050 & 0.150 & 0.140 \\
9  & \tt STRANGE\_FRACTION      & 0.50 & 0.531 & 0.521 & 0.529 & 0.514 & 0.497 \\
10 & \tt BARYON\_FRACTION       & 0.18 & 0.099 & 0.132 & 0.091 & 0.139 & 0.104 \\
11 & \tt P\_QS\_by\_P\_QQ\_norm & 0.48 & 0.720 & 0.617 & 0.502 & 0.601 & 0.720 \\
12 & \tt P\_SS\_by\_P\_QQ\_norm & 0.02 & 0.010 & 0.030 & 0.030 & 0.030 & 0.010 \\
13 & \tt P\_QQ1\_by\_P\_QQ0     & 1.00 & 1.499 & 1.499 & 1.349 & 1.164 & 1.499
	 \\\hline\hline
	& Euclidean distance from the default solution & & 1.508 & 	1.130 & 	1.400 & 	1.236 & 	1.497

	\\\hline %
		\end{tabular}

              \end{adjustbox}
\label{tab:sherpa_full}
            \end{table}

\begin{table}[htbp]
	\small
	\centering
	\caption{Optimal parameter values for the \sherpa dataset  obtained with all methods using the \textit{cubic polynomial} approximation after using outlier detection to  remove observables from the optimization (3 observables removed).}\label{tab:optimal_prms_obsfilter_sherpa}
	\begin{adjustbox}{width=\textwidth}
	\begin{tabular}{llllllll}
		\hline
	ID	&	Parameter name	& Default & Bilevel-meanscore & Bilevel-medscore & Bilevel-portfolio & Robust opt & All-weights-equal   \\\hline
1  & \tt KT\_0                  & 1.00 & 0.898 & 0.834 & 0.946 & 0.945 & 0.853 \\
2  & \tt ALPHA\_G               & 1.25 & 1.136 & 0.751 & 0.626 & 1.874 & 0.942 \\
3  & \tt ALPHA\_L               & 2.50 & 1.454 & 2.088 & 3.749 & 3.749 & 2.275 \\
4  & \tt BETA\_L                & 0.10 & 0.050 & 0.050 & 0.050 & 0.050 & 0.136 \\
5  & \tt GAMMA\_L               & 0.50 & 0.409 & 0.305 & 0.627 & 0.626 & 0.553 \\
6  & \tt ALPHA\_H               & 2.50 & 3.748 & 1.358 & 1.251 & 1.533 & 1.804 \\
7  & \tt BETA\_H                & 0.75 & 0.406 & 0.375 & 0.591 & 1.125 & 0.760 \\
8  & \tt GAMMA\_H               & 0.10 & 0.078 & 0.150 & 0.086 & 0.050 & 0.066 \\
9  & \tt STRANGE\_FRACTION      & 0.50 & 0.541 & 0.528 & 0.552 & 0.529 & 0.553 \\
10 & \tt BARYON\_FRACTION       & 0.18 & 0.181 & 0.139 & 0.090 & 0.270 & 0.270 \\
11 & \tt P\_QS\_by\_P\_QQ\_norm & 0.48 & 0.240 & 0.602 & 0.449 & 0.384 & 0.298 \\
12 & \tt P\_SS\_by\_P\_QQ\_norm & 0.02 & 0.020 & 0.025 & 0.030 & 0.030 & 0.023 \\
13 & \tt P\_QQ1\_by\_P\_QQ0     & 1.00 & 1.499 & 1.499 & 1.499 & 0.639 & 0.837
 \\\hline\hline
		& Euclidean distance from the default solution & & 1.222 & 	1.327 & 	1.378 & 	1.463 & 	0.937

		\\\hline %
		
	\end{tabular}
	\end{adjustbox}
\end{table}

\begin{table}[htbp]
	\centering
	\small
	\caption{Optimal parameter values for the \sherpa dataset obtained with all methods using the \textit{cubic polynomial} approximation after using the bin-filtering approach that excludes individual bins from the optimization  (204 bins out of 5246 total bins were removed).}\label{tab:optimal_prms_binfilter_sherpa}
	\begin{adjustbox}{width=\textwidth}
	\begin{tabular}{llllllll}
		\hline
	ID	& Parameter name		& Default & Bilevel-meanscore & Bilevel-medscore & Bilevel-portfolio & Robust opt & All-weights-equal   \\\hline
1  & \tt KT\_0                  & 1.00 & 0.866 & 0.820 & 0.897 & 0.950 & 0.911 \\
2  & \tt ALPHA\_G               & 1.25 & 0.626 & 1.114 & 0.626 & 1.874 & 0.626 \\
3  & \tt ALPHA\_L               & 2.50 & 3.749 & 3.502 & 2.216 & 3.749 & 3.749 \\
4  & \tt BETA\_L                & 0.10 & 0.079 & 0.053 & 0.050 & 0.050 & 0.078 \\
5  & \tt GAMMA\_L               & 0.50 & 0.383 & 0.325 & 0.750 & 0.627 & 0.367 \\
6  & \tt ALPHA\_H               & 2.50 & 1.251 & 1.251 & 1.251 & 1.527 & 1.251 \\
7  & \tt BETA\_H                & 0.75 & 0.738 & 0.675 & 0.694 & 1.125 & 0.741 \\
8  & \tt GAMMA\_H               & 0.10 & 0.092 & 0.116 & 0.099 & 0.050 & 0.104 \\
9  & \tt STRANGE\_FRACTION      & 0.50 & 0.536 & 0.547 & 0.543 & 0.529 & 0.541 \\
10 & \tt BARYON\_FRACTION       & 0.18 & 0.120 & 0.127 & 0.129 & 0.270 & 0.130 \\
11 & \tt P\_QS\_by\_P\_QQ\_norm & 0.48 & 0.636 & 0.578 & 0.472 & 0.384 & 0.569 \\
12 & \tt P\_SS\_by\_P\_QQ\_norm & 0.02 & 0.030 & 0.030 & 0.030 & 0.030 & 0.030 \\
13 & \tt P\_QQ1\_by\_P\_QQ0     & 1.00 & 1.499 & 1.499 & 1.499 & 0.637 & 1.499
 \\\hline
\hline
		& Euclidean distance from the default solution & & 1.263 & 	1.215 & 	1.272 & 	1.464 & 	1.224

 \\\hline %
							\end{tabular}

	\end{adjustbox}
\end{table}

\subsubsection{Comparison metric outcomes for the \sherpa dataset with the cubic polynomial approximation}
\label{sec:sherpa_30_metric}

Tables~\ref{tab:comparison_full_sherpa},~\ref{tab:comparison_obsf_sherpa}, and~\ref{tab:comparison_binf_sherpa} show the comparison metrics of our experiments when using the cubic polynomial approximation for the full data, the observable-filtered data, and the bin-filtered data, respectively.
Smaller numbers indicate better performance. The smallest number of each metric is bold for better visualization.
\begin{table}[]
	\caption{Comparison metrics for the full \sherpa  dataset when using the  \textit{cubic polynomial} approximation: Weighted $\chi^2$, A-optimality, D-optimality of different tuning methods. Smaller numbers are better and the best results are in bold. Note that we do not have an expert solution for this dataset.}
	\centering
	\begin{tabular}{l|rrr}
		\hline
		method         & Weighted $\chi^2$ & A-optimality    & D-optimality (log) \\\hline
		Bilevel-meanscore & 0.1777 &	9.0959 &	-39.9863
 \\
		Bilevel-medscore  & 0.2370 &	13.3943	& -37.1420
          \\
		Bilevel-portfolio & 0.3409 &	8.7863 &	-39.6956
         \\
		All-weights-equal         & 0.2305	& \textbf{6.8732} &	\textbf{-42.0678}
      \\
		Robust optimization          &  \textbf{0.0507} &	56.9168 &	-21.9561
       
\\\hline 
	\end{tabular}
	\label{tab:comparison_full_sherpa}
\end{table}

\begin{table}[h]
	\caption{Comparison metrics for the \sherpa dataset after filtering out  outlier observables and  using the \textit{cubic polynomial} approximation: Weighted $\chi^2$, A-optimality, D-optimality of different tuning methods.  3 observables were filtered out and not used during the optimization and results are computed based on all observables. Smaller numbers are better and the best results are indicated in bold. Note that we do not have an expert solution for this dataset.}
	\centering
	\begin{tabular}{l|rrr}
		\hline
		method         & Weighted $\chi^2$ & A-optimality & D-optimality(log) \\\hline
		Bilevel-meanscore & 0.4740&	14.3374 &	-35.2608
      \\
		Bilevel-medscore  & 0.4786 &	13.6299	& -36.6594
         \\
		Bilevel-portfolio & 0.2139	& \textbf{10.4481} &	\textbf{-36.8254}
         \\
		All-weights-equal         & 0.4789	& 28.2419	& -28.1536

\\
		Robust optimization     &    \textbf{0.0093} &	94.7811	& -23.5723
 \\\hline        
	\end{tabular}
	\label{tab:comparison_obsf_sherpa}
\end{table}

\begin{table}[]
	\caption{Comparison metrics for the \sherpa dataset after  bin-filtering and using the \textit{cubic polynomial} approximation: Weighted $\chi^2$, A-optimality, D-optimality of the different tuning methods.  204 out of 5426 bins were filtered out and not used during the optimization, and results were computed over all bins. Smaller numbers are better and the best results are indicated in bold. Note that we do not have an expert solution for this dataset.}
	\centering
	\begin{tabular}{l|rrr}
		\hline
		method         & Weighted $\chi^2$ & A-optimality & D-optimality(log) \\\hline
Bilevel-meanscore & 0.2504	& 17.2334	& -32.7683
          \\
		Bilevel-medscore  & 0.1835&	16.9248	&-32.1289
       \\
		Bilevel-portfolio & 0.2906 &	13.3500 &	-36.3598
 \\
		All-weights-equal         & 0.1928	& \textbf{10.4897} &	\textbf{-37.0305}

         \\
		Robust optimization            & \textbf{0.0364} &	72.5601	& -26.8516
     	\\\hline
	\end{tabular}
	\label{tab:comparison_binf_sherpa}
\end{table}

Based on these results, we can see that the all-weights-equal method (i.e. not adjusting any weights)  has the best performance for the full dataset under the A- and D-optimality. The bilevel-portfolio method performs best under the A- and D-optimality for both the observable- and bin-filtered datasets. The robust optimization method performs best in all three cases under the Weighted $\chi^2$ criterion.

\subsection{Weights assigned by different fitting methods}
\begin{figure}[h]
\label{fig:sherpa_weights}
  \includegraphics[width=\textwidth]{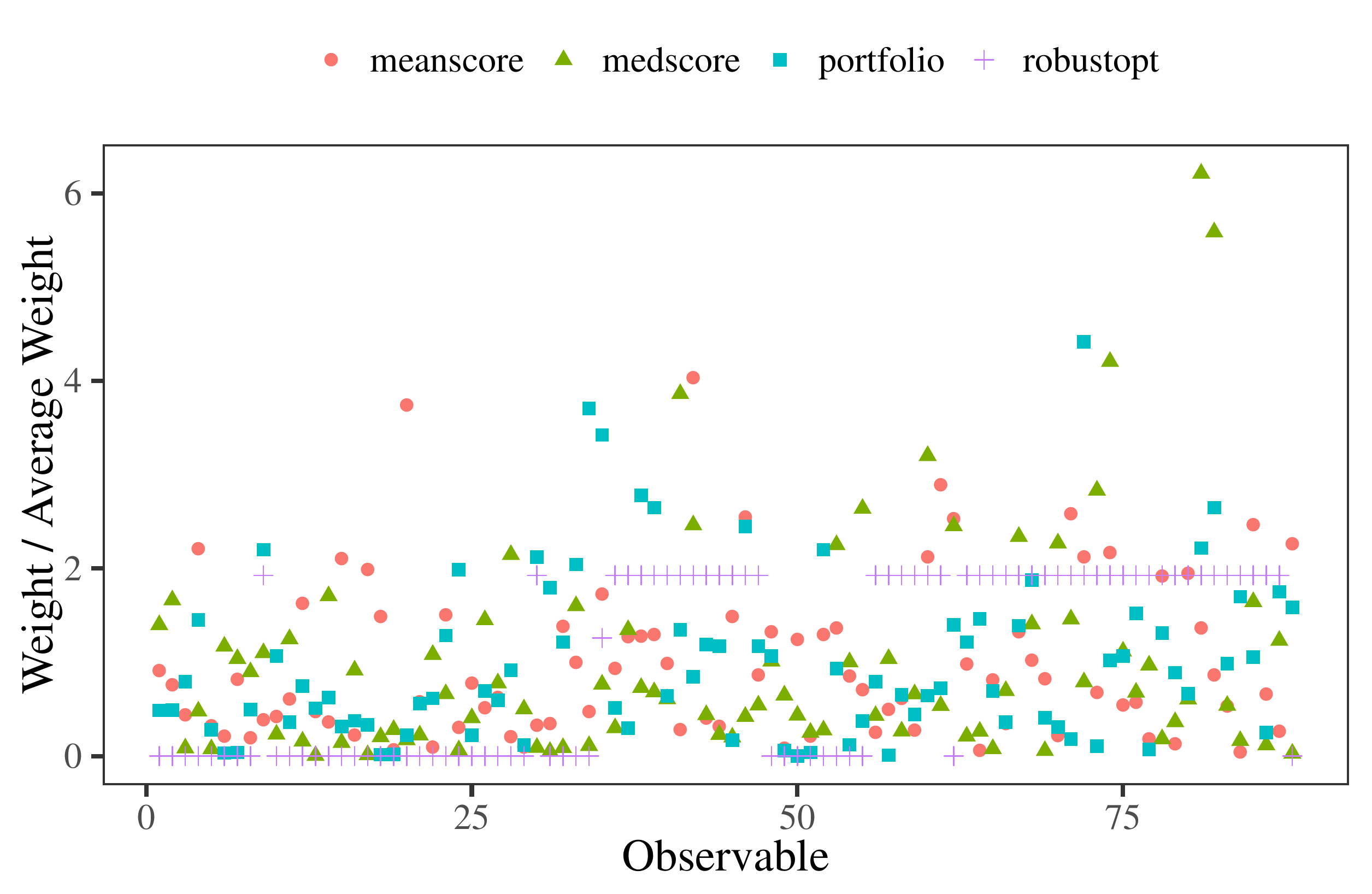}
  \caption{Distribution of weights assigned to observables for the
    different fitting methods described in the paper.
   Observables to the left are based on kinematic properties of
   events, while those to the right are particle multiplicities.}
\end{figure}

Figure~\ref{fig:sherpa_weights} shows the weights per observable
obtained from the tune to \sherpa using the methods described in this paper.

\subsection{Eigentunes for the results obtained with the cubic polynomial approximation}
Tables~\ref{tab:eigentune_a14_poly} and~\ref{tab:eigentune_sherpa_poly} shows the eigentune results for the A14 and \sherpa datasets, respectively, when using the cubic polynomial approximation.

\begin{table}[htbp]
\caption{Eigentune results for the A14 dataset using the optimal physics parameters $\p^*$ obtained with the different optimization methods when using the \textit{cubic polynomial} approximation.%
}
\label{tab:eigentune_a14_poly}
\begin{adjustbox}{width=\textwidth}
\begin{tabular}{l|rr|rr|rr|rr|rr}\hline
    Parameters & \multicolumn{2}{c|}{Expert} & \multicolumn{2}{c|}{Bilevel-meanscore} & \multicolumn{2}{c|}{Bilevel-medianscore} & \multicolumn{2}{c|}{Bilevel-portfolio} & \multicolumn{2}{c}{Robust optimization}  \\ 
     &  min & max &  min & max &  min & max &  min & max &  min & max\\  \hline 
     
\tt SigmaProcess:alphaSvalue &  0.072 &  0.196 &  0.071 &  0.197 &  0.079 &  0.190 &  0.076 &  0.191 &  0.079 &  0.187 \\
\tt BeamRemnants:primordialKThard &  1.899 &  1.904 &  1.849 &  1.888 &  1.877 &  1.894 &  1.855 &  1.881 &  1.764 &  1.895 \\
\tt SpaceShower:pT0Ref &  1.616 &  1.633 &  1.622 &  1.640 &  1.733 &  1.737 &  1.631 &  1.667 &  1.377 &  1.411 \\
\tt SpaceShower:pTmaxFudge &  0.904 &  0.914 &  0.938 &  0.940 &  0.884 &  0.923 &  0.986 &  0.990 &  0.932 &  0.935 \\
\tt SpaceShower:pTdampFudge &  1.039 &  1.047 &  1.059 &  1.102 &  1.053 &  1.085 &  1.045 &  1.049 &  1.061 &  1.064 \\
\tt SpaceShower:alphaSvalue &  0.116 &  0.128 &  0.128 &  0.130 &  0.118 &  0.141 &  0.129 &  0.131 &  0.128 &  0.129 \\
\tt TimeShower:alphaSvalue &  0.076 &  0.199 &  0.034 &  0.223 &  0.046 &  0.205 &  0.083 &  0.145 &  0.042 &  0.198 \\
\tt MultipartonInteractions:pT0Ref &  1.749 &  2.666 &  1.533 &  2.707 &  1.536 &  2.621 &  1.989 &  2.116 &  1.866 &  1.965 \\
\tt MultipartonInteractions:alphaSvalue &  0.045 &  0.186 &  0.095 &  0.154 &  0.114 &  0.140 &  0.044 &  0.180 &  0.100 &  0.133 \\
\tt BeamRemnants:reconnectRange &  1.719 &  1.719 &  1.523 &  1.541 &  1.390 &  1.420 &  1.589 &  1.595 &  2.565 &  2.568 \\

\hline
\end{tabular}
\end{adjustbox}
\end{table}

\begin{table}[htbp]
\caption{Eigentune results for the \sherpa dataset using the optimal physics parameters $\p^*$ obtained with the different optimization methods when using the \textit{cubic polynomial} approximation.}%
\label{tab:eigentune_sherpa_poly}
\begin{adjustbox}{width=\textwidth}
\begin{tabular}{l|rr|rr|rr|rr}\hline
    Parameters & \multicolumn{2}{c|}{Bilevel-meanscore} & \multicolumn{2}{c|}{Bilevel-medianscore} & \multicolumn{2}{c|}{Bilevel-portfolio} & \multicolumn{2}{c}{Robust optimization}  \\ 
     &  min & max &  min & max & min & max &  min & max\\  \hline 

\tt KT\_0 &  0.572 &  1.845 &  0.818 &  0.884 &  0.798 &  1.002 &  0.350 &  1.021 \\
\tt ALPHA\_G &  0.113 &  0.769 &  0.472 &  0.690 &  0.612 &  0.639 &  1.288 &  2.044 \\
\tt ALPHA\_L &  3.468 &  4.227 &  1.956 &  2.181 &  2.917 &  3.309 & 0 &  1.697 \\
\tt BETA\_L & 0 &  0.255 & 0 &  0.487 & 0 &  0.305 & 0 &  0.233 \\
\tt GAMMA\_L &  0.064 &  0.915 &  0.226 &  0.405 &  0.746 &  0.755 &  0.328 &  1.625 \\
\tt ALPHA\_H &  2.981 &  3.587 &  2.000 &  2.162 &  1.235 &  1.268 &  2.427 &  2.898 \\
\tt BETA\_H &  0.662 &  0.771 &  0.582 &  0.677 &  0.637 &  0.675 & 0 &  0.741 \\
\tt GAMMA\_H &  0.045 &  0.190 &  0.070 &  0.255 & 0 &  0.134 & 0 &  0.652 \\
\tt STRANGE\_FRACTION &  0.068 &  0.749 &  0.446 &  0.655 &  0.501 &  0.558 &  0.413 &  0.546 \\
\tt BARYON\_FRACTION & 0 &  0.335 &  0.117 &  0.166 & 0 &  0.186 &  0.030 &  0.516 \\
\tt P\_QS\_by\_P\_QQ\_norm &  0.669 &  0.828 &  0.576 &  0.715 &  0.458 &  0.549 &  0.537 &  0.619 \\
\tt P\_SS\_by\_P\_QQ\_norm & 0 &  0.087 & 0 &  0.105 & 0 &  0.076 & 0 &  0.050 \\
\tt P\_QQ1\_by\_P\_QQ0 &  1.496 &  1.508 &  1.498 &  1.500 &  1.348 &  1.349 &  1.153 &  1.200 \\

\hline
\end{tabular}
\end{adjustbox}
\end{table}

\subsection{Generator settings for \pythia and \sherpa}
\label{sec:gen_setup}
Typical run card for A14 studies using \pythia {\sc v8.186}.

\begin{Verbatim}[fontsize=\small]
Tune:pp = 14
Tune:ee = 7

PDF:useLHAPDF = on
PDF:LHAPDFset = NNPDF23_lo_as_0130_qed
PDF:LHAPDFmember = 0
PDF:extrapolateLHAPDF = off

! 3) Beam parameter settings. Values below agree with default ones.
Beams:idA = 2212                   ! first beam, p = 2212, pbar = -2212
Beams:idB = 2212                   ! second beam, p = 2212, pbar = -2212
Beams:eCM = 7000.                ! CM energy of collision

# uncomment for QCD
PhaseSpace:pTHatMin = 10.0
HardQCD:all = on
PhaseSpace:bias2Selection = on
PhaseSpace:bias2SelectionRef = 10.0
# uncomment for t-tbar
#Top:qqbar2ttbar = on
#Top:gg2ttbar = on
#SpaceShower:pTmaxMatch = 2
#SpaceShower:pTmaxFudge = 1
#SpaceShower:pTdampMatch = 1
# uncomment for Z
#WeakSingleBoson:ffbar2gmZ = On
#23:onMode = off
#23:onIfAny = 11 13 15 5 4 3
#SpaceShower:pTmaxMatch = 2
#SpaceShower:pTmaxFudge = 1
#SpaceShower:pTdampMatch = 1

# Example set of tuning parameters
SigmaProcess:alphaSvalue               0.1343
BeamRemnants:primordialKThard           1.711
SpaceShower:pT0Ref                      1.823
SpaceShower:pTmaxFudge                  1.047
SpaceShower:pTdampFudge                 1.492
SpaceShower:alphaSvalue                0.1302
TimeShower:alphaSvalue                 0.1166
MultipartonInteractions:pT0Ref          2.953
MultipartonInteractions:alphaSvalue     0.127
BeamRemnants:reconnectRange             4.747

ParticleDecays:limitTau0 = on
ParticleDecays:tau0Max = 10
\end{Verbatim}
We used these settings to reproduce the original results when necessary and
to make full predictions for parameters selected using the surrogate function.
Some of the original data using in the A14 study was private at that time and
was only made public later.   In a relatively small number of cases, the public
data was in a different form than that used for the original study, so we were unable
to reproduce those predictions.

Typical run card for \sherpa studies using {\sc v3.0.0}.
\begin{Verbatim}[fontsize=\small]
# general settings

SHOWER_GENERATOR: CSS
ANALYSIS: Rivet
FRAGMENTATION: Ahadic
INTEGRATION_ERROR: 0.02

# model parameters

ALPHAS(MZ): 0.1188
ORDER_ALPHAS: 2

# collider setup

BEAMS: [11, -11]
BEAM_ENERGIES: 45.6
  
# hadronization parameters 
AHADIC:
 KT_0  : 0.9088969039427998
 ALPHA_G  : 1.8736652396525728
 ALPHA_L  : 1.2518697247467987
 BETA_L  : 0.14989272155179253
 GAMMA_L  : 0.6832145156132761
 ALPHA_H  : 2.840868263919124
 BETA_H  : 0.5404054759080933
 GAMMA_H  : 0.14984034099619253
 STRANGE_FRACTION  : 0.5075082631730515
 BARYON_FRACTION  : 0.1357479921139296
 P_QS_by_P_QQ_norm  : 0.612797404412154
 P_SS_by_P_QQ_norm  : 0.029994467832440565
 P_QQ1_by_P_QQ0  : 1.1896505751927051

PARTICLE_DATA: 
  4: {Massive: true}
  5: {Massive: true}

PARTICLE_CONTAINER: 
  1098: {Name: C, Flavours: [4, -4]}
  1099: {Name: B, Flavours: [5, -5]}

PROCESSES: 
- 11 -11 -> 93 93: 
    Order: {QCD: 0, EW: 2}
- 11 -11 -> 4 -4: 
    Order: {QCD: 0, EW: 2}
- 11 -11 -> 5 -5: 
    Order: {QCD: 0, EW: 2}

RIVET:
   ANALYSES:
       - SLD_2002_S4869273
       - DELPHI_1996_S3430090
       - JADE_OPAL_2000_S4300807
       - PDG_HADRON_MULTIPLICITIES
\end{Verbatim}
We used these settings to reproduce the data for our surrogate function and 
to make full predictions for parameters selected using the surrogate function.

\nolinenumbers

\vfill
\begin{flushright}
\scriptsize
\framebox{\parbox{\textwidth}{
This manuscript has been authored by an author at Lawrence Berkeley National Laboratory under Contract No. DE-AC02-05CH11231 with the U.S. Department of Energy. The U.S. Government retains, and the publisher, by accepting the article for publication, acknowledges, that the U.S. Government retains a non-exclusive, paid-up, irrevocable, world-wide license to publish or reproduce the published form of this manuscript, or allow others to do so, for U.S. Government purposes.}}

\framebox{\parbox{\textwidth}{
The submitted manuscript has been created by UChicago Argonne, LLC, Operator of Argonne National Laboratory (“Argonne”). 
Argonne, a U.S. Department of Energy Office of Science laboratory, is operated under Contract No. DE-AC02-06CH11357. 
The U.S. Government retains for itself, and others acting on its behalf, a paid-up nonexclusive, irrevocable worldwide 
license in said article to reproduce, prepare derivative works, distribute copies to the public, and perform publicly 
and display publicly, by or on behalf of the Government.  The Department of Energy will provide public access to these 
results of federally sponsored research in accordance with the DOE Public Access Plan. 
\url{http://energy.gov/downloads/doe-public-access-plan}.
}}
\normalsize
\end{flushright}

\end{document}